\documentclass[10pt]{article}

\usepackage{framed}
\usepackage{amssymb}
\usepackage{latexsym}
\usepackage{url}
\usepackage{xcolor}
\usepackage{booktabs,multirow,arydshln}
\usepackage{setspace}
\definecolor{newcolor}{rgb}{.8,.349,.1}
\makeatletter
\def\ps@pprintTitle{%
	\let\@oddhead\@empty
	\let\@evenhead\@empty
	\def\@oddfoot{}%
	\let\@evenfoot\@oddfoot}
\makeatother
\usepackage[margin=1.5cm]{geometry}

\usepackage[titletoc]{appendix}

\usepackage[normalem]{ulem}

\usepackage{float}
\usepackage{amsmath}
\usepackage{amsthm}
\usepackage{psfrag}
\usepackage{graphicx}
\usepackage{fancyhdr}
\usepackage{verbatim}
\usepackage{blindtext}
\usepackage{stmaryrd}
\usepackage{soul}
\usepackage{sidecap}
\usepackage{tikz}
\usepackage{rotating}
\usepackage{comment}
\usepackage{todonotes}
\usepackage{changepage}
\usepackage{wrapfig}
\usepackage{subcaption}  
\usepackage{hyperref}
\usepackage{cleveref}

\crefname{appendix}{Appendix}{Appendices}
\Crefname{appendix}{Appendix}{Appendices}

\crefformat{appendix}{Appendix~#2#1#3}
\Crefformat{appendix}{Appendix~#2#1#3}
\crefmultiformat{appendix}{Appendix~#2#1#3}{ and~#2#1#3}{, #2#1#3}{ and~#2#1#3}
\Crefmultiformat{appendix}{Appendix~#2#1#3}{ and~#2#1#3}{, #2#1#3}{ and~#2#1#3}
\usepackage{natbib}
\setcitestyle{numbers}
\setcitestyle{square}
\definecolor{lightblue}{rgb}{0.32,0.45,0.90}
\definecolor{lightgreen}{rgb}{0.42,0.7,0.40}
\numberwithin{equation}{section}
\numberwithin{figure}{section}

\hypersetup{colorlinks=true,linkcolor = blue,citecolor=blue}

\usetikzlibrary{shapes,arrows}
\usetikzlibrary{arrows.meta}
\usetikzlibrary{decorations.markings}

\def\b{\boldsymbol}

\makeatletter
\newcommand{\vast}{\bBigg@{4}}
\newcommand{\Vast}{\bBigg@{5}}

\newcommand\bv{{\boldsymbol{v}}}

\newcommand\bx{\boldsymbol{x}}

\colorlet{cgray}{gray!20!white}

\theoremstyle{definition}
\newtheorem{definition}{Definition}[section]
\newtheorem{method}{Method}[section]

\newtheorem{theorem}{Theorem}[section]
\theoremstyle{remark}
\newtheorem*{remark}{Remark}
\newtheorem{example}{Example}[section]
\newtheorem{lemma}{Lemma}[section]
\newtheorem{corollary}{Corollary}[section]

\tikzset{every label/.style={font=\footnotesize,inner sep=1pt}}

\usepackage{pgfplots}
\pgfplotsset{compat = newest}
\usepgfplotslibrary{fillbetween,decorations.softclip}

\allowdisplaybreaks

\begin{document}

\title{High Order Multidimensional Slope Limiters with Local Maximum Principles}

\author{James Woodfield
\footnote{Department of Mathematics, Imperial College London, South Kensington Campus, London SW7 2AZ, United Kingdom.}
\footnote{Department of Mathematics and Statistics, University of Reading, Whiteknights, Reading RG6 6AX, United Kingdom.}
}

% #\addtocounter{footnote}{1} 
% \footnotetext{Department of Mathematics, Imperial College London, South Kensington Campus, London SW7 2AZ, United Kingdom.}
% #\addtocounter{footnote}{1}
% \footnotetext{Department of Mathematics and Statistics, University of Reading, Whiteknights, Reading RG6 6AX, United Kingdom.}
%\date{\today}
\date{}
\maketitle
\begin{abstract}
%%%
Higher-order numerical methods are used to find accurate numerical solutions to hyperbolic partial differential equations. Limiting is required to either converge to the correct type of solution or to adhere to physically motivated local maximum principles and less restrictive limiting procedures are required so as to not severely decrease the accuracy. 

In this paper, we adapt the existing slope limiter framework introduced in [Zhang \& Shu, J. Comput. Phys., 229(9):3091–3120, 2010] to achieve distinct $\emph{local}$ boundedness principles. We conclude that quadrature points contributing to numerical fluxes on either side of a face can be limited based on shared face-defined maximum principles and the resulting cell mean at the next timestep satisfies a cell mean maximum principle. Furthermore additional points arising in a decomposition of a cell mean must be limited locally when going beyond piecewise linear reconstructions. 

This allows the design of new multidimensional limiters which at second order can attain the same cell mean maximum principle as existing slope limiters, but allows more of the higher order flux to be used, generalises beyond second order schemes and can be modified for user specified local maximum principles.

% This allows more of the higher-order flux to be used whilst preserving the same cell mean maximum principle. We preliminarily investigate the practical application of the introduced framework to a second-order finite volume scheme as well as a fourth-order finite volume scheme, in the context of the advection equation and a 2D extension of Burgers' equation. 
%%%%
\end{abstract}

% \section{Introduction}
% \begin{document}

% \verso{Woodfield \textit{et al.}}

% \begin{frontmatter}

% \title{High Order Multidimensional Slope Limiters with Local Maximum Principles }

% \author[inst1,inst2]{James Woodfield\corref{cor1}}
% \cortext[cor1]{Corresponding author: 
%    Tel.: +447758899493; }
   
% \address[inst1]{Department of Mathematics, Imperial College London, South Kensington Campus, London SW7 2AZ, United Kingdom}

% \address[inst2]{Department of Mathematics and Statistics, University of Reading, Whiteknights, Reading RG6 6AX,United Kingdom}

% \received{31 May 2024}
% \finalform{-}
% \accepted{-}
% \availableonline{31 May 2024}
% \communicated{James Woodfield}

{\footnotesize\tableofcontents}
\section{Introduction}

\subsection{Historical context and motivation}\label{sec: intro: motivation: chapter 2}

Harten, Hyman, Lax and Keyfitz (HHLK) introduced a notion of monotonicity \cite{HHL_1976} suitable for numerical study of hyperbolic partial differential equations of various types, but also showed that such schemes (including nonlinear ones) must necessarily be first order. Since then, several different nonlinear limiting strategies have been proposed for more general schemes and meshes, typically with more relaxed definitions of monotonicity. For orthogonal meshes one dimensional slope limiting is sufficient for multidimensional incompressible flow to maintain discrete local maximum principles \cite{woodfield2024new,S_1984,S_1987}, however for more general meshes truly multidimensional slope limiting becomes necessary.
The unstructured multidimensional limiter of Barth-Jespersen \cite{BJ_1989} and the Kuzmin limiter \cite{kuzmin2010vertex} (equivalently the MLP limiter \cite{park2010multi}), are widely used and shown to be effective and discrete local maximum principle preserving for second-order methods on both unstructured and structured grids. However, existing BJ/Kuzmin-type limiters, and higher order extensions such as \cite{WANG2004716,yang2026toward}, do not provably enforce discrete local maximum principles once reconstructions go beyond piecewise linear.

By adapting the \emph{\textbf{global}}  boundedness limiter framework in Zhang and Shu \cite{zhang2010maximum} (see also \cite{zhang2010positivity,zhang2012maximum,fan2022positivity}), we derive sufficient conditions under which higher-order finite volume limiters enforce strictly stronger \emph{\textbf{local}} boundedness principles. Following \cite{zhang2010maximum}, we note additional points arising from a cell mean decomposition must be locally bounded when reconstructions go beyond piecewise linear. In addition, flux contributing quadrature points on both sides of a face can be limited by shared face defined maximum principles to reduce correction factors. This leads to new multidimensional limiters, which limit less than existing slope limiters (per defined maximum principle), enforce arbitrary local maximum principles, and can be extended to higher order methods and arbitrary meshes.

The scope of this paper concerns slope-limiting procedures with provable local cell mean maximum principles, rather than other notions of non-physical oscillation control.
There exist several useful multidimensional and one-dimensional limiting frameworks capable of suppressing some aspects of nonphysical oscillatory behaviour without enforcing local maximum principles. WENO (weighted essentially non-oscillatory) and its many variants \cite{jiang1996efficient,shu1999high,shu2003high,gerolymos2009very,balsara2016efficient,levy2000compact} are known examples where strict maximum principles are avoided but oscillations are still controlled locally. Other multidimensional frameworks exist where strict local maximum principles are deliberately avoided \cite{liu2017novel,venkatakrishnan1993accuracy,michalak2006differentiability,ollivier1997quasi,nishikawa2022new} in favour of differentiability. In other works \cite{wang2002spectral,wang2003high,WANG2004716,wang2004spectral,yang2026toward}, the subcell solution at reconstructed quadrature values (within spectral volumes) are bounded in terms of their local neighbour cell mean values in direct analogy to the Barth-Jespersen limiter. Such limiters, have been widely used to help control oscillations but have no formal discrete (local or global) cell mean maximum principle. 

% and there is a tendency in the literature for the ``maximum principle" to be imposed on larger and larger stencils as correction factors are required to be less restrictive. However, these limiters do not bound points arising in an acceptable decomposition of the cell average, consequently a discrete local maximum principle for the cell mean value at the next timestep is never established, and positivity or global boundedness is not proven (a counter example is found in \cite{zhang2010maximum}). 

% There already exists slope limiting procedures that can attain distinct local cell mean maximum principles but only for linear subcell representations (second order schemes). Of particular note is the Barth-Jespersen limiter \cite{BJ_1989} which limits the slope of the cell to ensure flux contributing quadrature points are bounded by a cell mean defined maximum principle. Another notable approach is taken by Kuzmin's vertex limiter equivalent to Park et al's MLP limiter, \cite{park2010multi,kuzmin2010vertex} where corner points are limited based on a local vertex neighbourhood (following the structured work in \cite{yoon2008multi}). 

\subsection{Contributions}
We develop a local maximum principle framework for unstructured finite volume (FV) schemes (\cref{thm: slope limiters}) following \cite{zhang2010maximum}. For arbitrary-order FV schemes on unstructured meshes, we introduce new limiter classes dependent on fixed neighbourhoods $GN(K)$ and $BJ(K)$. One introduced limiter is a generalisation of the BJ limiter for arbitrary-order FV schemes, termed the GBJ--$BJ(K)$ limiter, which differs from the literature by locally bounding non-flux contributing quadrature points arising from a decomposition of the cell average. Another novel limiter is the Sharp--$GN(K)$ limiter based upon sufficient conditions for \cref{thm: slope limiters}, which differs by using face defined local maximum principles for provably less restrictive correction factors. Another new limiter class is the $GN(K)$ limiter based on more interpretable face defined bounds. 

For arbitrary-order schemes on general 2D and 3D unstructured meshes, we prove a discrete local maximum principle (DLMP) for the GBJ--$BJ(K)$ limiter, the Sharp--$GN(K)$ limiter, and the $GN(K)$ limiter. We also prove given any GBJ--$BJ(K)$ limiter (including the original BJ limiter), one can always define a Sharp--$GN(K)$ limiter and a $GN(K)$ limiter that enforce the same cell-mean maximum principle with less severe correction factors.

At second order with a linear subcell representation on general 2D and 3D unstructured meshes, the GBJ-$N(K)\cup K$-limiter becomes the original BJ limiter, and we prove that the Sharp--$N^2(K)\cup N(K)$ limiter, and the $N^2(K)\cup N(K)$ limiter, enforce the same cell-mean maximum principle with less severe correction factors. At second order with a linear subcell representation on general 2D and 3D unstructured meshes, we prove the Sharp--$VN(K)$ limiter, the $VN(K)$ limiter, and the original Kuzmin \cite{kuzmin2010vertex}/MLP \cite{park2010multi} limiter enforce the same vertex-based cell-mean maximum principle, while the Sharp--$VN(K)$ and $VN(K)$ limiters use less severe correction factors.

We introduce a new fourth-order FV scheme (FV4), with a novel symmetry-based decomposition of the cell average that enables all new local boundedness limiters to be implemented. Numerical testing is consistent with analytical correction factor inequalities, and proven DLMP's.

\subsection{Outline of the paper}
The remainder of this paper is organised as follows.
In \cref{sec:equations} we introduce the equations of interest, and the desired monotonicity properties of the scheme. In \cref{sec:fv_update_fluxes} we introduce the numerical scheme and admissible meshes considered. In \cref{sec:mesh} we introduce unstructured mesh notation. 

In \cref{sec: Theory: theorems: chapter 2} we indicate sufficient conditions for an arbitrary order scheme to retain a local boundedness principle on an unstructured mesh (\cref{thm: slope limiters}). In \cref{sec: Theory: theorems: chapter 2} we introduce the higher order generalisation of the BJ-limiter (GBJ-$BJ(K)$), the Sharp-General-Neighbourhood-$GN(K)$-limiter, and its intersect variant $GN(K)$. Local maximum principles are proven in \cref{sec: Theory: theorems: chapter 2}, and correction factor inequalities in \cref{sec:Comparing generalised local limiting procedures}. In \cref{sec:instanciations_of_limiters} we instantiate some limiters for later numerical experiments.

The first example is in \cref{sec: Application: second order method: chapter 2} where a second order finite volume method (called FV2) is employed on a structured grid. The second example (\cref{sec: Application2: fv4 chapter 2}) introduces a fourth order advection algorithm (called FV4), we then demonstrate how \cref{thm: slope limiters} can be used for a new limiting strategy once a novel type of decomposition of the cell average is found and present the numerical results.

In \cref{sec:FV1} we review extended notions of monotonicity developed by Harten, Hyman, Lax, and Keyfitz \cite{hhlk1976} as a self contained \emph{First order} introduction for motivating the desired notions of monotonicity in this paper. In \cref{proof:bj_derivation} we elaborate on the correction factor function. 
In \cref{sec:fv2}, we define a second order finite volume scheme for unstructured meshes in 2D and 3D, derive 2D and 3D Zhang-acceptable decompositions of the cell mean, provide limiters (under the FV2 simplification) and prove correction factor inequalities for the Kuzmin-limiter.

\subsection{Equations}\label{sec:equations}
We consider numerical solutions to the equation
\begin{align}
u_t + \operatorname{div}\left(f(u)\b v\right) = 0,\quad u(\b x,0) = u_0(\b x) \label{eq:eq},
\end{align}
over $\b x \in \Omega\subseteq\mathbb{R}^d$, $t\in \mathbb{R}^{\geq 0}$, $d\in\lbrace 2,3\rbrace$. Subject to the usual assumptions \cite{EGH_2000} of bounded initial data $u_0\in L^{\infty}(\Omega)$, $u_{0}\in [m,M]$, $m,M\in \mathbb{R}$, continuous ($\b v \in C^{1}(\Omega\times \mathbb{R}^{\geq 0};\mathbb{R}^d)$) potentially divergence free ($\operatorname{div}(\b v)= 0$) bounded ($\exists V\in \mathbb{R}$ s.t. $|\b v|<V$) velocity $\b v$, and continuous $f$. With a specific focus on flux form advection where $f(u)=u$, $d=2$. The solution to \cref{eq:eq} is sign-preserving, i.e. a positive solution remains positive for all time
\begin{align}
u_0(\b x)\geq 0 \implies u(\b x,t) \geq 0 ,\quad \forall t\in [0,T],\quad \forall \b x \in \Omega. 
\end{align}
Conditional on the divergence-free property of $\b v$, solutions to \cref{eq:eq} also satisfy global maximum principles of the type 
\begin{align}
u_0(\b x)\in [m,M] \implies u(\b x,t) \in [m,M],\quad \forall t\in [0,T],\quad \forall \b x \in \Omega. 
\end{align} 
Such properties can be expected of unique (weak) entropy solutions, and are often desired in numerical schemes. 

Higher order globally bounded numerical solutions can be constructed using the slope limiter framework in \cite{zhang2010maximum}. In particular, let $\bar{u}_K^{n}$ denote the cell mean within a cell $K$ belonging to the mesh $\mathcal{M}$ at the timestep $n$. Then it can be shown \cite{zhang2010maximum} numerical approximations to \cref{eq:eq} can be made sign-preserving in the sense
\begin{align}
\bar{u}^{0}_{K}\geq 0 \implies \bar{u}^{n}_{K} \geq 0 ,\quad \forall n\in \lbrace 1,...,N\rbrace,\quad \forall K \in \mathcal{M}, 
\end{align}
and conditional on a discrete divergence-free property of $\b v$ satisfy
global
maximum principles of the form 
\begin{align}
\bar{u}^{0}_K \in [m,M]\implies \bar{u}^n_K \in [m,M],\quad  \forall n\in  \lbrace1,...,N\rbrace,\quad \forall K\in \mathcal{M}.
\end{align}
This paper describes the extension or application of this framework to preserve stronger local maximum principles of the type
\begin{align}
\bar{u}^{n+1}_K \in [m^{n}_{K},M^{n}_K],\quad \forall K\in \mathcal{M}
\end{align}
where $m^{n}_{K}, M^{n}_{K}$ are locally defined time-dependent maxima and minima, typically chosen as maxima and minima of cell mean values at the previous time level $t^{n}$ over a locally defined neighbourhood of cell $K$. In the remainder of this work we will omit the superscript $n$ where it is apparent that it refers to the $n$-th time step value.

\subsection{Mesh and Method}\label{sec:fv_update_fluxes}

\begin{figure}[H]
\centering
	\hspace{10pt}	\includegraphics[scale=1.2,trim={10mm 5mm 40mm 10mm},clip]{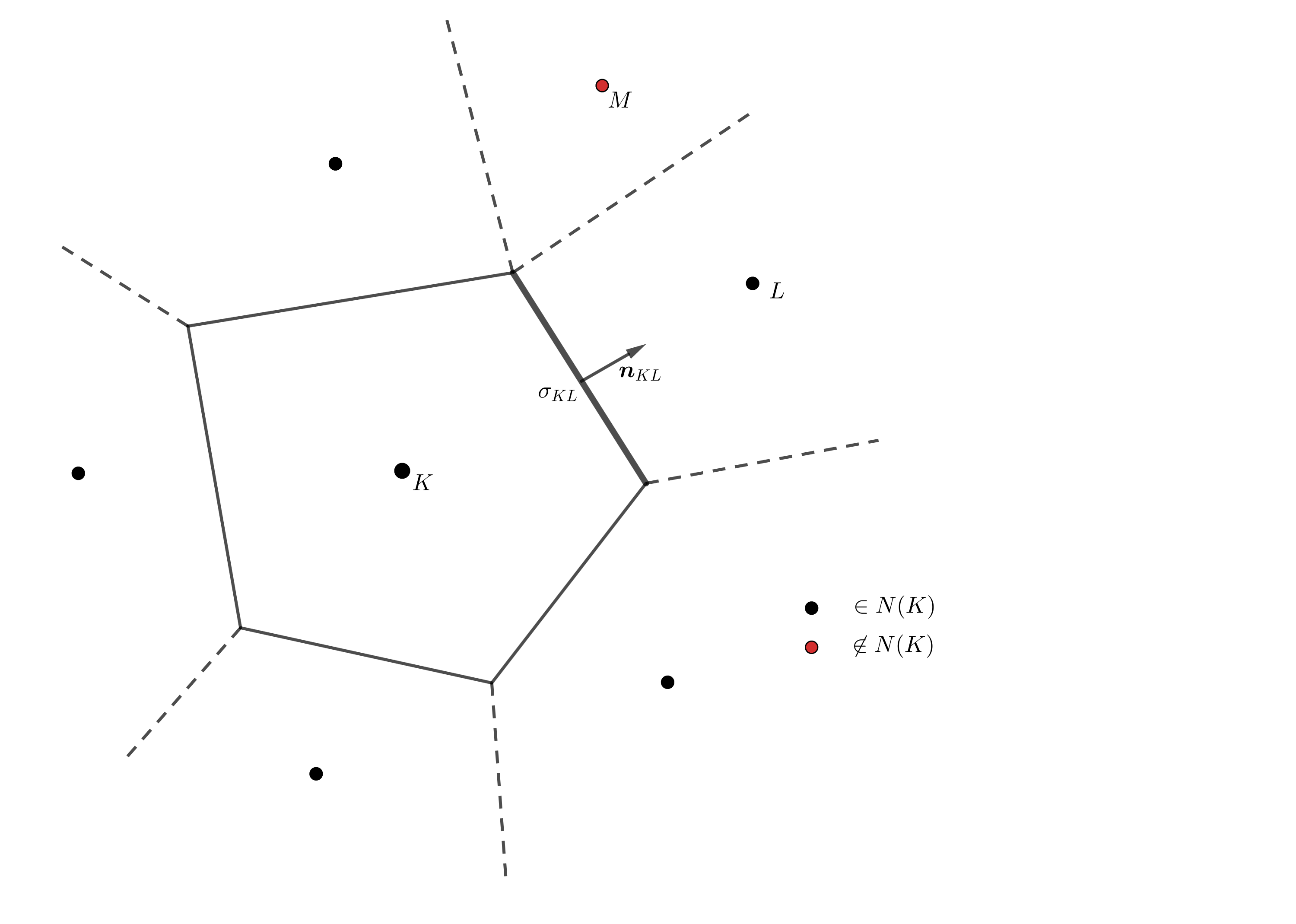} 
	\vspace{-10pt}
	\caption{Cell $K$, and the face $\sigma_{KL}$ of a face sharing neighbour $L \in N(K)$, with outward unit normal $\b n_{KL}$.}
	\label{stencil:unstructured_chap2}
\end{figure}

Let $\mathcal M$ be a partition of $\Omega$ into control volumes $K$. For each face-sharing neighbour $L\in N(K)$, let $\sigma_{KL}=\partial K\cap \partial L$ be the common face, oriented by the outward unit normal $\b n_{KL}$ from $K$ (so $\b n_{LK}=-\b n_{KL}$). Denote the cell measure by $|K|$ and the face measure by $|\sigma_{KL}|$. We call $\mathcal{M}$ an admissible finite volume mesh if: each cell $K$ is a convex polygon (2D) or polyhedron (3D), each face $\sigma_{KL}$ lies in a $(d-1)$-dimensional affine hyperplane, and each face is shared by at most two cells. See \cref{stencil:unstructured_chap2}.

A single forward-Euler step for a finite volume scheme can be written in flux form as
\begin{align}\label{eq:FV_FE_generic}
\bar u^{\,n+1}_K
=
\bar u^{\,n}_K
-\frac{\Delta t}{|K|}\sum_{L\in N(K)}\,
\mathcal F_{KL}^{\,n},
\end{align}
where $\bar u^n_K$ is the cell mean of cell $K$ at time level $n$, $\mathcal F_{KL}$ is the approximation of the flux through the face $\sigma_{KL}$. The flux computation uses a face defined quadrature rule with normalised weights and points (assumed not on vertices in this work) denoted $\lbrace  w_q^{\sigma_{KL}} , \b x^{\sigma_{KL}}_q  \rbrace_{q=1}^{Q_{KL}}$ as follows
\begin{align}
\mathcal F_{KL} := |\sigma_{KL}|\sum_{q\in \sigma_{KL}^{I}} w_q^{\sigma_{KL}} f_{KL}(p_K(\b x^{\sigma_{KL}}_q),p_L(\b x^{\sigma_{KL}}_q);\b v(\b x_q^{\sigma_{KL}}) \cdot \b n_{KL})\approx \int_{\b x\in \sigma_{KL}} f_{KL}(p_K(\b x), p_L(\b x); \b v(\b x) \cdot \b n_{KL})dS
\end{align}
Here we have used the associated index set by $\sigma^{I}_{KL}:=\{1,\dots,Q_{KL}\}$, to account for the possibility of different number of points per face.
Where $p_K(\b x)$ denotes the subcell approximation in cell $K$ and $f_{KL}$ denotes a conservative consistent monotone flux, from cell $K$, to cell $L$, defined as follows. 

\begin{definition}[Consistent, conservative, monotone flux function]\label{def:flux conservative consistent lipshitz monotone}
A two-point numerical flux $f_{KL}(a,b;c)$ is
\emph{consistent} if $f_{KL}(u,u;c)=f(u)c$, $ \forall u\in\mathbb R,\ \forall c\in\mathbb R.$
It is \emph{conservative} if
$f_{KL}(a,b;c)=-f_{LK}(b,a;-c)$, $\forall a,b,c\in\mathbb R$. 
It is \emph{monotone} if it is nondecreasing in its first argument and nonincreasing in its second: $\partial_{a}f_{KL}(a,b;c)\geq 0$, $\partial_{b}f_{KL}(a,b;c)\leq 0$,
for all $a,b,c\in\mathbb R$. Examples can be found in \cref{sec:FV1}.
\end{definition}

\subsection{Mesh neighbourhood notation}\label{sec:mesh}
% Before introducing new limiters, we first remark that several different local maximum principles are already proposed to control non-physical oscillations.
% The one-dimensional limiting procedures in \cite{woodfield2024new} produces a numerical solution with a local maximum principle on the inclusive face sharing neighbourhood, whilst requiring less stringent requirements of the subcell representation. The Barth-Jespersen limiter \cite{BJ_1989}, limits the subcell representation locally, and satisfies a maximum principle based on the ``squared" inclusive face sharing neighbourhood. The Kuzmin limiter \cite{kuzmin2010vertex} satisfies a maximum principle based on the inclusive vertex sharing neighbours. 

We define some neighbourhoods below
\begin{itemize}
\item $GN(K)$ denotes a set of cells in the general neighbourhood of cell $K$ containing the cell $K$.
	\item $N(K) :=\lbrace L\in\mathcal{M}\setminus{\lbrace K \rbrace}   : |\sigma_{KL}|>0\rbrace$ denotes the set of face sharing neighbours of cell $K$.
	\item $N(K)\cup K$ denotes the inclusive face sharing neighbourhood i.e. the set of cells attainable from cell $K$ by at most one face sharing step.
    \item $N^2(K):=\bigcup_{L\in N(K)}N(L)$ denotes the set of cells, which are attainable in exactly two face sharing steps from cell $K$, and contains the cell $K$, but not $N(K)$.
	\item $N^{2}(K) \cup N(K)$ denotes the ``squared" inclusive face sharing neighbourhood, containing all cells attainable in at most two face sharing steps from cell $K$.
	\item $N(v)\;:=\lbrace K\in\mathcal{M}: v\in K\rbrace$ denotes the set of cells sharing vertex $v$.
    \item $V(K)$ denotes the set of vertices of cell $K$.
	\item $VN(K) := \cup_{v\in V(K)} N(v)$ is the inclusive vertex neighbourhood of a cell $K$, i.e $L\in VN(K)$ iff $L$ shares a vertex with $K$.
    \item $V(\sigma_{KL}):=\lbrace v\in \mathcal{V} : v\in K\cap L\rbrace$, denotes those vertices on a face $\sigma_{KL}$.
    \item $VN(\sigma_{KL}):= \left\lbrace K\in \mathcal{M}:\bigcup_{v\in V(\sigma_{KL})}N(v)\right\rbrace$ denotes the set of cells, sharing a vertex with both cell $K$ and cell $L$. (vertex neighbours of face $\sigma_{KL}$).
\end{itemize}

\begin{figure}[H]
	\centering
	\includegraphics[scale=0.45]{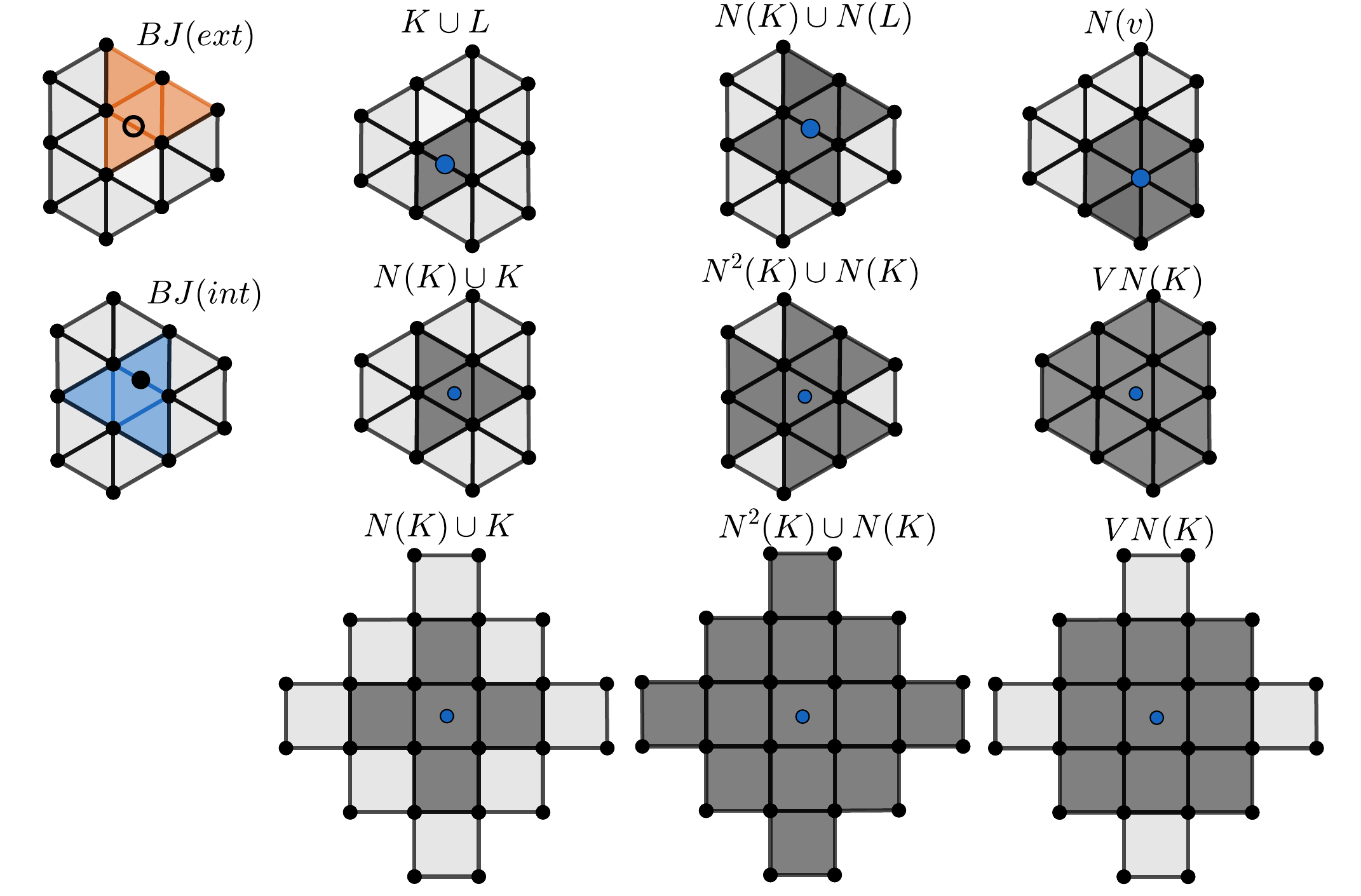}
	\caption{Visualisations of some neighbourhoods for some common meshes. Blue dot is an informal representation of the middle of the neighbourhood. Dark grey denotes the cells contained in the neighbourhood of the midpoint.  Blue and orange regions are specific to the Barth-Jespersen limiter.}
	\label{fig:neighboorhoods}
\end{figure}

\Cref{fig:neighboorhoods} visualises some of these different regions for some different meshes.

\section{High order multidimensional slope limiting}\label{sec: Theory: theorems: chapter 2}
We consider a (higher-order) forward Euler cell-mean evolution equation of the form
\begin{align}
	\bar{u}^{n+1}_{K}
	&= \bar{ u}^{n}_{K}
	- \frac{ \Delta t }{ |K|} \sum_{L\in N(K)} |\sigma_{KL}| \sum_{q \in \sigma^{I}_{KL}} w_{q}^{\sigma_{KL}}
	f_{KL}\!\left( p_{K}(\b x^{\sigma_{KL}}_{q}),\ p_{L}(\b x^{\sigma_{KL}}_{q}),\ \b v(\b x^{\sigma_{KL}}_{q}) \cdot \b n_{KL}\right). \label{eq:dg}
\end{align}
Here $p_{K}(\b x)$ denotes the subcell representation approximating the true solution in cell $K$.
For each face $\sigma_{KL}$ we fix a quadrature rule with $Q_{KL}$ nodes $\{\b x^{\sigma_{KL}}_{q}\}_{q=1}^{Q_{KL}} \subset \sigma_{KL}$ not on vertices and corresponding \emph{normalised} weights $\{w_{q}^{\sigma_{KL}}\}_{q=1}^{Q_{KL}}$.
We denote the associated index set by $\sigma^{I}_{KL}:=\{1,\dots,Q_{KL}\}$, and when $q\in\sigma^I_{KL}$, we write $\b x_q:=\b x^{\sigma_{KL}}_q$. 
The numerical flux in \cref{eq:dg} arises from the quadrature approximation of the flux integral 
\begin{align}
	\int_{\b x\in \sigma_{KL}} f_{KL}\left(p_{K}(\b x),p_{L}(\b x),\b v(\b x)\cdot \b n_{KL}\right)dS
	\approx |\sigma_{KL}| \sum_{q \in \sigma^{I}_{KL}} w_{q}^{\sigma_{KL}}
	f_{KL}\!\left( p_{K}(\b x_{q}), p_{L}(\b x_{q}), \b v(\b x_{q}) \cdot \b n_{KL}\right).
\end{align}
Normalisation and positivity of quadrature weights give
\begin{align}
w_{q}^{\sigma_{KL}}>0, \quad \forall q\in \sigma_{KL}^I,
\qquad
\sum_{q\in \sigma_{KL}^I}  w_{q}^{\sigma_{KL}}= 1,
\qquad
\forall L\in N(K), \quad\forall K\in \mathcal{M} . \label{eq:positivity_face_quadrature}
\end{align}

The scheme \cref{eq:dg} is not a monotone function of surrounding cell mean values in the HHLK \cite{HHL_1976} sense (unless the subcell approximation is piecewise constant \cref{sec:FV1}). However \cite{zhang2010maximum,zhang2012maximum} show that under suitable decompositions of the cell average, the scheme is a monotone function of quadrature point evaluations and non-flux contributing quadrature points.
The key assumption is that the cell mean $\bar{u}_K$ admits a decomposition as a weighted combination of flux-contributing quadrature point evaluations and non-flux contributing point evaluations and all such weights are positive.

To formalise this, we introduce an index set $K^{nfc}:=\{1,\dots,Q^{nfc}_K\}$ for the non-flux contributing points
$\{\b x^K_q\}_{q\in K^{nfc}}\subset \overline K$, and for $q\in K^{nfc}$ we write $\b x_q:=\b x_q^{K}$. In the cell-mean decomposition we assign weights $\gamma^K_q$ to the non-flux points ($q\in K^{nfc}$), and weights
$\gamma^K_{q,L}$ to the face quadrature points ($q\in\sigma^I_{KL}$, $L\in N(K)$). Then the following cell-mean decomposition is assumed:
\begin{align}
	\bar{u}_{K}
	= \frac{1}{|K|}\int_{K} p_{K}(\b x)\, d\b x
	= \sum_{q \in K^{nfc}}  p_{K}(\b x_q)\, \gamma^K_{q}
	+ \sum_{L\in N(K)}\sum_{q\in \sigma^I_{KL}} p_{K}(\b x^{\sigma_{KL}}_q)\, \gamma^K_{q,L}.
	\label{eq: cell mean decomposition}
\end{align}
where the weights satisfy positivity and summation to one must arises if the decomposition holds for a constant,
\begin{align}
(\gamma_{q}^{K}>0, \quad \forall q\in K^{nfc}), \quad (\gamma^K_{q,L}>0,\quad \forall q \in \sigma_{KL}^{I},\quad \forall L\in N(K)),\quad
\forall K\in \mathcal{M}. \label{eq:positivity_cell_mean}
\\
  \sum_{q \in K^{nfc}}  \gamma^K_{q}
	+ \sum_{L\in N(K)}\sum_{q\in \sigma^I_{KL}}  \gamma^K_{q,L}= 1,
\quad
\forall K\in \mathcal{M}. \label{eq:sum_to_one}
\end{align}

We call \eqref{eq: cell mean decomposition} Zhang-acceptable when all weights are positive (i.e. \cref{eq:positivity_cell_mean}) and sum to one (\cref{eq:sum_to_one}). Technically, one can allow a flux contributing point to live on both LHS and RHS of the decomposition.
\Cref{thm: slope limiters} describes the sufficient conditions for a local cell mean boundedness principle. The method of proof is similar to that in \cite{zhang2010maximum,zhang2012maximum}; but treats \emph{local} maximum principles under incompressible flow and arbitrary meshes.

\begin{theorem}\label{thm: slope limiters} The cell mean value at the next time-step $\bar{u}_{K}^{n+1}$ evolved by the cell mean evolution equation \cref{eq:dg} with a conservative consistent monotone flux function (\cref{def:flux conservative consistent lipshitz monotone}) can be expressed as a monotone function of quadrature point evaluations and points arising in a decomposition of the cell mean \cite{zhang2010maximum,zhang2012maximum}. If all quadrature point evaluations arising from a Zhang-acceptable cell mean decomposition are non-negative sum to one, and all the face defined Courant number restrictions are satisfied:
	\begin{align}
		\Delta t \frac{w^{\sigma_{KL}}_{q} |\sigma_{KL}|}{\gamma^{K}_{q,L}|K|}\frac{\partial f_{KL}}{\partial p_{K}(\b x_{q})}\leq 1, \quad \forall q \in \sigma_{KL}^{I},\quad \forall L \in N(K),\quad \forall K \in\mathcal{M}.
	\end{align}
	This implies positivity preservation $\bar{u}_{K}^{n+1}\geq 0$ for compressible and incompressible flow \cite{zhang2010maximum,zhang2012maximum}, with negativity and sign preservation holding similarly. If, in addition, the following discrete divergence free condition holds, 
	\begin{align}
		\frac{1}{|K|}\sum_{L\in N(K)}|\sigma_{KL}| \sum_{q \in \sigma_{KL}^{I}} w_{q}^{\sigma_{KL}}( \b v(\b x^{\sigma_{KL}}_q)\cdot \b n_{KL}) = 0,
	\end{align}
    and the values at quadrature points are locally bounded as follows,
	\begin{align}
		&p_{K}(\b x_q) \in[m_{K},M_{K}],\quad \forall q \in K^{nfc}, \\
		&p_{K}(\b x^{\sigma_{KL}}_q), p_{L}(\b x^{\sigma_{KL}}_q) \in[m_{K},M_{K}],\quad \forall q \in \sigma_{KL}^{I},\quad \forall L \in N(K),
	\end{align}
	then the next time level will satisfy a $\bar{u}_{K}^{n+1} \in [m_K,M_K]$ local boundedness principle.
\end{theorem}

\begin{proof}{(of {\cref{thm: slope limiters}})} Use the Zhang-acceptable abstract cell mean decomposition \cref{eq: cell mean decomposition}, to write the numerical method \cref{eq:dg} as a positive sum of non-flux contributing quadrature point evaluations and HHLK-monotone \cite{hhlk1976} two-point schemes at the flux contributing quadrature points as follows
	\begin{align}
		\bar{u}_{K}^{n+1} &= \bar{u}^{n}_{K}   -  \frac{\Delta t}{ |K|} \sum_{L\in N(K)}  |\sigma_{KL}| \sum_{q \in \sigma_{KL}^{I}} w_{q}^{\sigma_{KL}}f_{KL}\left(  p_{K}(\b x_q),  p_{L}(\b x_q),  \b v(\b x_q,t)\cdot \b n_{KL}\right),\\
		&= \sum_{q \in K^{nfc}} \gamma_{q}^{K}  p_{K}(\b x_{q})  + \sum_{L\in N(K)} \sum_{q \in \sigma_{KL}^{I}} \gamma_{q,L}^{K} \left( p_{K}(\b x_{q}) - \Delta t \frac{w^{\sigma_{KL}}_{q} |\sigma_{KL}|}{\gamma^{K}_{q,L}|K|}f_{KL}( p_{K}(\b x_{q}), p_{L}(\b x_{q}); \b v(\b x_{q},t)\cdot \b n_{KL})\right).
	\end{align}
    The derivative of the numerical method \cref{eq:dg} with respect to each quadrature point evaluation is given by
	\begin{align}
		\frac{\partial \bar{u}_{K}^{n+1}}{\partial p_{K}(\b x_{q})} & = \gamma_{q}^{K}, \quad \forall q \in K^{nfc},\\
		\frac{\partial \bar{u}_{K}^{n+1}}{\partial p_{K}(\b x_{q})} & =\gamma_{q,L}^K \left[ 1-   \frac{\Delta t |\sigma_{KL}| w_{q}^{\sigma_{KL}}}{|K| \gamma^{K}_{q,L}}\frac{\partial f_{KL}}{\partial p_{K}(\b x_{q})} \right], \quad \forall q \in  \sigma^{I}_{KL}, \quad \forall L\in N(K), \\
		\frac{\partial \bar{u}_{K}^{n+1}}{\partial  p_{L}(\b x_{q})} & = -\frac{\Delta t |\sigma_{KL}| w_{q}^{\sigma_{KL}}}{|K|}\frac{\partial f_{KL}}{\partial  p_{L}(\b x_{q})}, \quad \forall q \in  \sigma_{KL}^{I}, \quad \forall L\in N(K).
	\end{align}
The weight positivity of the cell mean decomposition $(\gamma_{q,L}^{K}> 0, \forall q \in \sigma_{KL}^{I},\forall L\in N(K))$, $(\gamma_{q}^{K}> 0$, $\forall q \in K^{nfc})$, the positivity of flux weights $(w_{q}^{\sigma_{KL}} > 0, \forall q \in \sigma_{KL}^{I}, \forall L\in N(K))$, the monotone property of the flux $\partial_{a} f_{KL}(a,b;c) \geq 0, \partial_{b} f_{KL}(a,b;c)\leq 0$, and the flux contributing time-step restrictions $\frac{\Delta t |\sigma_{KL}| w_{q}^{\sigma_{KL}}}{|K|\gamma^K_{q,L}}\frac{\partial f_{KL}}{\partial  p_{K}(\b x_{q})} \leq 1, \quad \forall q \in \sigma_{KL}^{I},\forall L\in N(K)$, imply all above derivatives are non-negative. Proving, the scheme is a monotone function of quadrature point evaluations, and points arising in the decomposition of the cell mean. This implies sign preservation of the scheme. For instance, given an arbitrary velocity field, the following conditions
	\begin{align}
		&p_{K}(\b x_{q}) \geq 0,\quad \forall q \in K^{nfc}, \\
		&p_{K}(\b x_{q}),\quad p_{L}(\b x_{q}) \geq 0,\quad \forall q \in \sigma_{KL}^{I},\quad \forall L \in N(K),\\
		&\Delta t \frac{w^{\sigma_{KL}}_{q} |\sigma_{KL}|}{\gamma^{K}_{q,L}|K|}\frac{\partial f_{KL}}{\partial p_{K}}\leq 1, \quad \forall q \in \sigma_{KL}^{I},\quad \forall L \in N(K),
	\end{align}
	are sufficient for the scheme to be positivity preserving and $\bar{u}^{n+1}_{K}\geq 0$. 
	If, in addition, we suppose that the following discrete divergence free condition holds,
	\begin{align}
		\frac{1}{|K|}\sum_{L\in N(K)}|\sigma_{KL}| \sum_{q \in \sigma_{KL}^{I}} w_{q}^{\sigma_{KL}}( \b v(\b x_q,t) \cdot \b n_{KL}) &= 0,
	\end{align}
	and that the numerical fluxes are consistent, then we can derive equation consistency, from the identity
	\begin{align}
		\frac{1}{|K|}\sum_{L\in N(K)} |\sigma_{KL}| \sum_{q \in \sigma_{KL}^{I}} w_{q}^{\sigma_{KL}} \Big( f_{KL}(c,c;\b v(\b x_{q},t)\cdot \b n_{KL}) \Big)  &= 0, \label{eq:consistency_ho}
	\end{align}
	and it implies the preservation of constants of the scheme. More formally writing \cref{eq:dg} as
    
    \begin{align}
    \bar{u}^{n+1}_K = H\left(\lbrace p_K(\b x_q) \rbrace_{ q \in K^{nfc}}, \left\lbrace \lbrace p_{K}(\b x_q)\rbrace_{ q \in \sigma_{KL}^{I}},\lbrace p_{L}(\b x_q)\rbrace_{ q \in \sigma_{KL}^{I}} ,\lbrace \b v(\b x_{q},t^n) \cdot \b n_{KL} \rbrace_{ q \in \sigma_{KL}^{I}} \right\rbrace_{L\in N(K)}\right),
    \end{align}
    the identity \cref{eq:consistency_ho} implies 
    \begin{align}
    c = H\left(\lbrace c \rbrace_{ q \in K^{nfc}}, \left\lbrace \lbrace c\rbrace_{ q \in \sigma_{KL}^{I}},\lbrace c\rbrace_{ q \in \sigma_{KL}^{I}} ,\lbrace \b v(\b x_{q},t^n) \cdot \b n_{KL} \rbrace_{ q \in \sigma_{KL}^{I}} \right\rbrace_{L\in N(K)}\right).
    \end{align}
If in addition we assume that the internal and boundary quadrature points are locally bounded as follows,
	\begin{align}
		&p_{K}(\b x_{q}) \in[m_{K},M_{K}],\quad \forall  q \in K^{nfc}, \\
		&p_{K}(\b x_{q}), p_{L}(\b x_{q}) \in[m_{K},M_{K}],\quad \forall q \in \sigma_{KL}^I,\quad \forall L \in N(K),
	\end{align}
we can then deduce the following maximum principle, 
 \begin{align}
		m_{K} \leq \bar{u}_{K}^{n+1} \leq M_{K},
	\end{align}
 by the monotonicity and the consistency of the numerical method, by considering $c=m_K$ and $c=M_K$.
 %as follows
 %\begin{align}
 %m_{K} &= H(\lbrace m_{K}\rbrace_{\forall q \in K^{nfc}},\lbrace m_{K}\rbrace_{\forall q \in K^{fc}},\lbrace m_{K}\rbrace_{\forall q \in K^{fc}}; \lbrace \bv(x_q,t)\cdot \b n_{KL}\rbrace_{\forall q \in K^{fc}}) \\
 %&\leq H(\lbrace p_K(\b x_q) \rbrace_{\forall q \in K^{nfc}}, \lbrace p_{K}(\b x_q)\rbrace_{\forall q \in K^{fc}},\lbrace p_{L}(\b x_q)\rbrace_{\forall q \in K^{fc}}; \lbrace \b v(\b x_{q},t^n) \cdot \b n_{KL} \rbrace_{\forall q \in K^{fc}}) %\\
 %&\leq H(\lbrace M_{K}\rbrace_{\forall q \in K^{nfc}},\lbrace M_{K}\rbrace_{\forall q \in K^{fc}},\lbrace M_{K}\rbrace_{\forall q \in K^{fc}}; \lbrace \bv(x_q,t)\cdot \b n_{KL}\rbrace_{\forall q \in K^{fc}})= M_{K}.
 %\end{align}
\end{proof}
\phantom{ok}\newline

The main distinction from the literature is that the flux contributing quadrature points at a face $\sigma_{KL}$ must satisfy two local boundedness principles,
\begin{align}
	p_{L}(\b x_{q}), p_{K}(\b x_{q}) \in[m_{K},M_{K} ], \\
	p_{L}(\b x_{q}), p_{K}(\b x_{q}) \in[m_{L},M_{L}],
\end{align}
when the requirements of \cref{thm: slope limiters} are viewed from the perspective of cells $K,L$ respectively (to ensure $\bar{u}_{K}^{n+1}\in[m_K,M_K]$, $\bar{u}_{L}^{n+1}\in[m_L,M_L]$). This has important consequences on the design of multidimensional limiter functions, since we must require $p_{L}(\b x_{q}), p_{K}(\b x_{q}) \in[m_{K},M_{K}]\cap [m_{L},M_{L}]$, as summarised in the following \cref{cor:limiting proceedure}. 

\phantom{ok}\newline

\begin{corollary}
[Local boundedness limiter]\label{cor:limiting proceedure} Consider a FV scheme \cref{eq:dg} such that the conditions in \cref{thm: slope limiters} hold. Let the flux contributing points $p_{L}(\b x_{q}), p_{K}(\b x_{q})$ be limited based on the shared face defined maximum principle,
\begin{align}
p_{L}(\b x_{q}), p_{K}(\b x_{q}) \in \left[\max \lbrace m_{K}, m_{L} \rbrace , \min \lbrace  M_{K}, M_{L} \rbrace\right] \quad \forall q \in \sigma_{KL}^{I}, \quad \forall L\in N(K).
\end{align}

Let $p(\b x_{q})$ arising from a cell average decomposition, not contributing to a flux be limited via
\begin{align}
 p_{K}(\b x_{q}) \in \left[ m_{K}, M_{K}\right], \quad \forall q \in K^{nfc}.
\end{align}
Then the next timestep will satisfy the maximum principle 
\begin{align}
\bar{u}^{n+1}_{K} \in \left[m_K,M_K\right].
\end{align}
\end{corollary}
\begin{proof}
Follows immediately from \cref{thm: slope limiters}. \Cref{cor:limiting proceedure} defines a class of limiting method, sufficient to ensure that $\bar{u}_{K}^{n+1}$ is bounded by some cell defined bounds $[m_K,M_K]$. 
\end{proof}

Before defining any methods we must first define correction factors used to limit the value of a polynomial evaluated at a particular point following \cite{barth2003finite}.
\begin{definition}
[Barth-Jespersen correction factors \cite{barth2003finite}]\label{BJ_correction_factors}
 Given a subcell representation $p_K(\b x)$ within cell $K$, and slope limited representation $\tilde{p}_{K}(\b x) = \bar{u}_K + \alpha (p_{K}(\b x)-\bar{u}_K)$, the Barth-Jespersen correction factor $\alpha_q$, at the point $\b x_{q}\in K$, is defined as the \emph{\textbf{largest}} $\alpha \in[0,1]$ such that $\tilde{p}_K(\b x_{q})\in [m,M]$, and can be equivalently be defined as
    \begin{align}
  \alpha_q =   \alpha(m,M,\bar{u}_{K},p_K(\b x_q)) = \begin{cases}
			\min \left\lbrace 1, \frac{M - \bar{u}_{K}}{p_K(\b x_{q}) - \bar{u}_{K}}  \right\rbrace \quad \text{if} \quad  p_K(\b x_{q}) - \bar{u}_{K}>0,\\
			\min \left\lbrace 1,\frac{m - \bar{u}_{K}}{p_K(\b x_{q}) -\bar{u}_{K}} \right\rbrace \quad \text{if} \quad  p_{K}(\b x_{q}) - \bar{u}_{K}<0,\\
			1 \quad  \text{if} \quad  p_K(\b x_{q}) - \bar{u}_{K} =0.
		\end{cases}
	\end{align}
\end{definition}
\begin{remark}
    The derivation of this explicit formula, is included in \cref{proof:bj_derivation}, to emphasise this is the largest admissible correction factor to ensure boundedness of $\tilde{p}_K$ at $\b x_q$.
\end{remark}

By combining \cref{cor:limiting proceedure,thm: slope limiters} with \cref{BJ_correction_factors} one attains the following method, capable of achieving different local maximum principles. 

\begin{method}[Sharp--\protect\ensuremath{GN(K)} limiter]\label{method:sharp_m_K_M_K}
\noindent 
Let $GN(K)\subseteq\mathcal{M}$ be a fixed neighbourhood of $K$, with $K\in GN(K)$.
\begin{enumerate}
    \item  Per cell $K$ we compute the desired local maximum principle,
	\begin{align}
		[m_{K} , M_{K}] := \left[ \min_{L\in GN(K)}\bar{u}_{L}^{n}, \max_{L \in GN(K)}\bar{u}_{L}^{n}\right].
	\end{align}
     This is associated to each non-flux contributing quadrature point $\b x_{q} \in K^{nfc}$.
    \item Per face $\sigma_{KL}$, we compute  the local face defined maximum principle bounds
	\begin{align}
		[m_{\sigma_{KL}},M_{\sigma_{KL}}]= \left[\max \left\lbrace m_{K} , m_{L}\right\rbrace ,\min\left\lbrace M_{K} , M_{L}\right\rbrace \right],\quad \forall L\in N(K).
	\end{align}
	This is associated to each flux contributing quadrature point $\b x_{q}$ such that $q\in \sigma_{KL}^{I}$, $\forall L\in N(K)$
	\item 
	We then compute per cell all the Barth-Jespersen quadrature corrections factors $\alpha_{q}$, to ensure $\tilde{p}_{K}(x) = \alpha (p_{K}(x) - \bar{u}_{K}) + \bar{u}_{K}$ satisfies the conditions in \cref{thm: slope limiters} at required quadrature points, i.e.
	\begin{align}
		&\tilde{p}_{K}(\b x_{q}) \in [m_{K},M_{K}],\quad \forall q \in K^{nfc},\\
		&\tilde{p}_{K}(\b x_{q}) \in [m_{\sigma_{KL}},M_{\sigma_{KL}}],\quad \forall q \in \sigma_{KL}^{I} , \quad \forall L\in N(K),
	\end{align}
	by choosing the smallest value $\alpha = \min_{\forall q \in K^{I} } \alpha_{q}$.
    This is done explicitly via \cref{BJ_correction_factors} by computing
    \begin{align}
    \alpha^{\sigma_{KL}}_{q} = \alpha_{q}(m_{\sigma_{KL}},M_{\sigma_{KL}},\bar{u}_{K},p_K(\b x_{q})), \quad \forall q \in \sigma^{I}_{KL} , \quad \forall L\in N(K),\\
    \alpha^{I}_{q}=\alpha_{q}(m_{K},M_{K},\bar{u}_{K},p(\b x_{q})), \quad \forall q \in K^{nfc},
    \end{align}
    then choosing
    \begin{align}
        \alpha_{K}: = \min\left( \min_{\forall L\in N(K)} \left[\min_{q \in \sigma^{I}_{KL} } \alpha^{\sigma_{KL}}_{q}\right] , \min_{q\in K^{nfc}}\alpha^{I}_q \right).
    \end{align}
\end{enumerate}
\end{method}

\begin{remark}For FV methods with greater than two nodes per face, one need only compute two BJ-correction factors per face since one can bound the worst violators of the shared face-defined maximum principle,
\begin{align}
\min_{q \in \sigma^{I}_{KL} } \alpha^{\sigma_{KL}}_{q} = \min \left( \alpha_{q}\left(m_{\sigma_{KL}},M_{\sigma_{KL}},\bar{u}_{K},\min_{q\in \sigma_{KL}^{I}}p(\b x_{q})\right), \alpha_{q}\left(m_{\sigma_{KL}},M_{\sigma_{KL}},\bar{u}_{K},\max_{q\in \sigma_{KL}^{I}}p(\b x_{q})\right) \right). 
\end{align}
\end{remark}

\begin{corollary}
The arbitrary order FV scheme in \cref{eq:dg} with the sharp--$GN(K)$ limiter under the conditions in \cref{thm: slope limiters} attains the discrete local maximum principle:
\begin{align}
   \min_{L\in GN(K)} \bar{u}_{L}^{n}  \leq \bar{u}^{n+1}_K \leq \max_{L\in GN(K)} \bar{u}_{L}^{n}.
\end{align}
\end{corollary}
\begin{proof}\Cref{thm: slope limiters}, \cref{cor:limiting proceedure} where $[m_K,M_{K}] = [\min_{L\in GN(K)} \bar{u}_{L}^{n},\max_{L\in GN(K)} \bar{u}_{L}^{n}]$. 
\end{proof}

\begin{method}[\protect\ensuremath{GN(K)} limiter]\label{method:Dull_m_K_M_K}
\noindent 
Identical to \cref{method:sharp_m_K_M_K}, however per face $\sigma_{KL}$, we instead compute  the local face-defined maximum principle bounds via the intersect bounds,
\begin{align}
	[m_{\sigma_{KL}},M_{\sigma_{KL}}]= \left[\min_{i\in GN(K)\cap GN(L) }\bar{u}^n_{i} , \max_{i\in GN(K)\cap GN(L) }\bar{u}^n_{i}\right],\quad \forall L\in N(K).
\end{align}
This is associated to each flux contributing quadrature point $\b x_{q}$ such that $q\in \sigma_{KL}^{I}$, $\forall L\in N(K)$. The face defined bounds here are less sharp from that in \cref{method:sharp_m_K_M_K} and \cref{thm: slope limiters}, to be slightly more interpretable. For instance, in 1D one can have face-defined bounds like $u_{j+1/2}\in [\min(\bar{u}_{j},\bar{u}_{j+1}),\max(\bar{u}_{j},\bar{u}_{j+1})]$, rather than $u_{j+1/2}\in [\max(\min(\bar{u}_{j-1},\bar{u}_{j},\bar{u}_{j+1}),\min(\bar{u}_{j},\bar{u}_{j+1},\bar{u}_{j+2})),\min(\max(\bar{u}_{j-1},\bar{u}_{j},\bar{u}_{j+1}),\max(\bar{u}_{j},\bar{u}_{j+1},\bar{u}_{j+2}))]$.
\end{method}
\begin{corollary}
The arbitrary order FV scheme in \cref{eq:dg} with the $GN(K)$ limiter under the conditions in \cref{thm: slope limiters} attains the discrete local maximum principle:
\begin{align}
   \min_{L\in GN(K)} \bar{u}_{L}^{n}  \leq \bar{u}^{n+1}_K \leq \max_{L\in GN(K)} \bar{u}_{L}^{n}.
\end{align}
\end{corollary}
\begin{proof}\Cref{thm: slope limiters}, \cref{cor:limiting proceedure}. 
\end{proof}

Two limiting procedures were proposed for a fourth-order spectral volume method in \cite{WANG2004716}. Both limit the flux quadrature points to be cellwise bounded, closely mirroring the philosophy of the second-order Barth–Jespersen limiter. However, these procedures do not establish a positivity result, a discrete local maximum principle, or global boundedness for the updated cell mean. Moreover, \cite{zhang2010maximum} (Remark 2.3) highlights that controlling flux quadrature points alone is generally insufficient: bounding internal quadrature points can be necessary to guarantee a maximum principle for the cell mean. Motivated by this, we introduce a generalised BJ limiter that also constrains non-flux quadrature points and allows arbitrary neighbourhood choices. 

\begin{method}[GBJ--$BJ(K)$ limiter.] \label{method:GBJ}
Generalised Barth-Jespersen--$BJ(K)$ limiter.
Let $BJ(K)\subseteq\mathcal{M}$ denote a fixed neighbourhood of $K$ with $K\in BJ(K)$, typically taken to be $N(K)\cup K$. 
\begin{enumerate}
    \item Per cell we compute the cell defined bounds
    \begin{align}
		[m^{BJ}_{K},M^{BJ}_{K}] := \left[ \min_{L\in BJ(K)}\bar{u}_{L}^{n}, \max_{L\in BJ(K)}\bar{u}_{L}^{n}\right]
	\end{align} 
    associated to all points in the Zhang-acceptable decomposition of the cell mean.
    \item 
    We then compute per cell all the Barth- Jespersen quadrature corrections factors $\alpha_{q}$, to ensure $\tilde{p}_{K}(x) = \alpha (p_{K}(x) - \bar{u}_{K}) + \bar{u}_{K}$, satisfy the following bounds
	\begin{align}
		&\tilde{p}_{K}(\b x_{q}) \in [m^{BJ}_{K},M^{BJ}_{K}],\quad \forall q \in K^{nfc},\\
		&\tilde{p}_{K}(\b x_{q}) \in [m^{BJ}_{K},M^{BJ}_{K}],\quad \forall q \in \sigma_{KL}^{I} , \quad \forall L\in N(K),
	\end{align}
	by choosing the smallest value $\alpha^{BJ} = \min_{\forall q \in K^{I} } \alpha_{q}$.
    This is done explicitly via correction factors \cref{BJ_correction_factors}, and a refactoring speed-up in \cref{eq:bjfast}.
\end{enumerate}
\end{method}
\begin{remark}[Relationship to the original BJ-limiter]
    At second order with a linear representation in each cell, the Generalised Barth-Jespersen--$N(K)\cup K$ limiter becomes the usual Barth-Jespersen limiter. This arises for piecewise linear representations on unstructured meshes since there exists a decomposition of the cell average onto only flux contributing facet midpoints in both 2D and 3D, established using properties of affine functions in \cref{sec:fv2}.
\end{remark}

\begin{remark}[Fast implementation GBJ-$BJ(K)$]
We say this limiting procedure is of BJ-type, since each subcell evaluation of $p_K$ is bounded with the same cell wise maximum principle. This simplification improves compute speed, since the limiter can be refactored to limit worst violators of the cell-defined BJ maximum principle as follows:
\begin{align}
\alpha_K=\min\left(\alpha\left(m^{BJ}_K,M^{BJ}_K,\bar{u}_K,\max_{x_q\in K}p_K(x_q)\right),\alpha\left(m^{BJ}_K,M^{BJ}_K,\bar{u}_K,\min_{x_q\in K}p_K(x_q)\right)\right), \label{eq:bjfast}
\end{align}
giving two correction factor computations per cell.
\end{remark}
% \begin{remark}
%     For linear representations there exists a decomposition of the cell average onto flux contributing facet midpoints in both 2D and 3D established using properties of affine functions in \cref{sec:fv2}, which allows limiting flux contributing quadrature points to be sufficient for a DLMP.
% \end{remark}

\begin{corollary}
The arbitrary order FV scheme in \cref{eq:dg} with the GBJ-$BJ(K)$ limiter under the conditions in \cref{thm: slope limiters}, does \textbf{\emph{not}} attain 
\begin{align}
    \min_{L\in BJ(K)} \bar{u}_{L}^{n} \leq \bar{u}^{n+1}_K \leq \max_{L\in BJ(K)} \bar{u}_{L}^{n},
\end{align}
but attains the weaker discrete local maximum principle:
\begin{align}
   \min_{L\in N(K)\cup K} \left(\min_{M\in BJ(L)} \bar{u}_{M}^{n} \right) \leq \bar{u}^{n+1}_K \leq \max_{L\in N(K)\cup K} \left(\max_{M\in BJ(L)} \bar{u}_{M}^{n}  \right).\label{eq:BJ_cell_mean_values}
\end{align}
\end{corollary}

\begin{proof}
Generalised-Barth-Jespersen-$BJ(K)$-limiter  correction factors at flux contributing quadrature points are computed via the bounds:
\begin{align}
p_K(\b x_{q})\in [m^{BJ}_{K},M^{BJ}_{K}]= \left[\min_{i\in BJ(K) }\bar{u}^n_{i} , \max_{i\in BJ(K) }\bar{u}^n_{i}\right],\quad \forall q\in \sigma_{KL}^{I},\quad \forall L\in N(K). 
\end{align}
However, $p_{L}(\b x_{q})$ for $q\in \sigma_{KL}^I$ doesn't share these bounds and satisfies
\begin{align}
p_L(\b x_{q})\in [m^{BJ}_{L},M^{BJ}_{L}]= \left[\min_{i\in BJ(L) }\bar{u}^n_{i} , \max_{i\in BJ(L) }\bar{u}^n_{i}\right],\quad \forall q\in \sigma_{KL}^{I},\quad \forall L\in N(K).
\end{align}
This is not sufficient for \cref{thm: slope limiters}, to apply on the $BJ(K)$-region. Instead, this is only sufficient for the extended maximum principle in \cref{eq:BJ_cell_mean_values}, via \cref{thm: slope limiters}.
\end{proof}

\subsection{Factors affecting accuracy: comparing generalised local limiting procedures}\label{sec:Comparing generalised local limiting procedures}
We require the use of the following lemma when comparing the relative size of different correction factors.
\begin{lemma}[Correction factor inequalities]\label{lemma:correction factor inequalities} Given $\bar{u}_K$, $p_K(\bx_q)$, if $[m_1,M_1]
\subseteq [m_2,M_2]$, then the correction factors $\alpha_q$ defined in \cref{BJ_correction_factors}, satisfy 
\begin{align}
\alpha(m_1,M_1,\bar{u}_K,p_K(\bx_q))\leq \alpha(m_2,M_2,\bar{u}_K,p_K(\bx_q)).
\end{align}
\end{lemma}
\begin{proof}
If $p_K(\bx_q) -\bar{u}_{K}>0$, then $\min\left(1,\frac{M_1-\bar{u}_{K}}{p_K(\bx_q) -\bar{u}_K}\right)\leq\min\left(1,\frac{M_2-\bar{u}_{K}}{p_K(\bx_q) -\bar{u}_K}\right)$, since $M_1\leq M_2$. Similarly, if  $p_K(\bx_q) -\bar{u}_K<0$, then $\min\left(1,\frac{m_1-\bar{u}_{K}}{p_K(\bx_q) -\bar{u}_K}\right)\leq\min\left(1,\frac{m_2-\bar{u}_{K}}{p_K(\bx_q) -\bar{u}_K}\right)$, since $m_1\geq m_2$. If $p_K(\bx_q) -\bar{u}_K=0$, then $1 = 1$. 
\end{proof}

\begin{theorem}
Let $GN(K)$ denote a general neighbourhood of cells containing cell $K$. Let conditions in \cref{thm: slope limiters}, hold. 
The FV scheme \cref{eq:dg} with either the Sharp-$GN(K)$-limiter with  correction factor $\alpha^{\sharp}$ or the $GN(K)$-limiter with correction factor $\alpha^{GN}$ (from the same arbitrary initial condition $u^n$) produce different cell mean values satisfying the same discrete local maximum principle,
\begin{align}
   \min_{L\in GN(K)} \bar{u}_{L}^n \leq \bar{u}_{K}^{n+1}\leq \max_{L\in GN(K)} \bar{u}^n_{L}, \quad \forall K\in \mathcal{M}.
\end{align}
However, the correction factors are less limiting in the Sharp--$GN(K)$ limiter, 
\begin{align}
  \alpha_{K}^{GN}  \leq \alpha_{K}^{\sharp}, \quad \forall K\in \mathcal{M}.
\end{align}
\end{theorem}
\begin{proof}
Discrete local maximum principles follow from \cref{thm: slope limiters} and \cref{cor:limiting proceedure}. The correction factors only differ at flux contributing quadrature points. 
Per face $\sigma_{KL}$, in the Sharp--$GN(K)$ limiter we compute  the local face defined maximum principle bounds via
\begin{align}
[m^{\sharp}_{\sigma_{KL}},M^{\sharp}_{\sigma_{KL}}]= \left[\max \left\lbrace \min_{i\in GN(K)}\bar{u}_{i}^{n} , \min_{i\in GN(L)}\bar{u}_{i}^{n}\right\rbrace ,\min\left\lbrace\max_{i\in GN(K)}\bar{u}_{i}^{n} , \max_{i\in GN(L)}\bar{u}_{i}^{n} \right\rbrace \right],\quad \forall L\in N(K).
	\end{align}
Per face $\sigma_{KL}$, in the $GN(K)$ limiter we compute the local face defined maximum principle bounds via
\begin{align}
[m_{\sigma_{KL}},M_{\sigma_{KL}}]= \left[\min_{i\in GN(K)\cap GN(L) }\bar{u}^{n}_{i} , \max_{i\in GN(K)\cap GN(L) }\bar{u}^{n}_{i}\right],\quad \forall L\in N(K).
\end{align}
Therefore, since $GN(K)\cap GN(L) \subseteq GN(K)$ and $GN(K)\cap GN(L) \subseteq GN(L)$, one attains
\begin{align}
[m_{\sigma_{KL}},M_{\sigma_{KL}}]  \subseteq [m^{\sharp}_{\sigma_{KL}},M^{\sharp}_{\sigma_{KL}}] .
\end{align}
This is sufficient to conclude via \cref{lemma:correction factor inequalities} that $\alpha_{q}\leq \alpha^{\sharp}_{q}$, $\forall q\in \sigma_{KL}^{I}$, $\forall L\in N(K)$, $\forall K\in \mathcal{M}$, and since $\alpha_{q}= \alpha^{\sharp}_{q}$, $\forall q\in K^{nfc}$, $\forall K\in \mathcal{M}$ one has
\begin{align}   
    \alpha^{GN}_{K}\leq \alpha^{\sharp}_{K}, \quad \forall K \in \mathcal{M}.
\end{align}

\end{proof}

\begin{theorem}\label{thm: comparing limiters GBJ-SharpGN-GN} Given any Barth-Jespersen cell set $BJ(K)$,
define $GN(K)= \bigcup_{L\in N(K)\cup K} BJ(L)$. Let conditions in \cref{thm: slope limiters}, hold. 
The scheme \cref{eq:dg} with either the Sharp-$GN(K)$-limiter with  correction factor $\alpha^{\sharp}_K$ or the $GN(K)$-limiter with correction factor $\alpha^{GN}_K$ or the Generalised-Barth-Jespersen--$BJ(K)$ limiter with correction factor $\alpha^{GBJ}_K$ (from the same arbitrary initial condition $u^n$)
all satisfy the same discrete local maximum principle,
\begin{align}
   \min_{L\in GN(K)} \bar{u}_{L}^n \leq  \bar{u}_{K}^{n+1}\leq \max_{L\in GN(K)} \bar{u}_{L}^n, \quad \forall K\in \mathcal{M}.
\end{align}
However, the correction factors satisfy
\begin{align}
  \alpha_{K}^{GBJ}  \leq \alpha_{K}^{GN}  \leq \alpha_{K}^{\sharp}, \quad \forall K\in \mathcal{M}.\label{eq:GBJ_GN_SGN}
\end{align} 
\end{theorem}
\begin{proof} Maximum principles are already proven via \cref{thm: slope limiters}.
Generalised--Barth--Jespersen--$BJ(K)$ correction factors at flux contributing quadrature points are computed via the bounds:
\begin{align}
p_K(\b x_{q})\in \left[m^{BJ}_{K},M^{BJ}_{K}\right]= \left[\min_{i\in BJ(K) }\bar{u}_{i} , \max_{i\in BJ(K) }\bar{u}_{i}\right],\quad \forall q\in \sigma_{KL}^{I},\quad \forall L\in N(K),\\
p_L(\b x_{q})\in \left[m^{BJ}_{L},M^{BJ}_{L}\right]= \left[\min_{i\in BJ(L) }\bar{u}_{i} , \max_{i\in BJ(L) }\bar{u}_{i}\right],\quad \forall q\in \sigma_{KL}^{I},\quad \forall L\in N(K),
\end{align}
whereas the $GN(K)$-limiter has correction factors computed via the bounds
\begin{align}
p_{K}(\b x_{q}),p_{L}(\b x_{q})\in [m_{\sigma_{KL}},M_{\sigma_{KL}}]= \left[\min_{i\in GN(K)\cap GN(L) }\bar{u}_{i} , \max_{i\in GN(K)\cap GN(L) }\bar{u}_{i}\right],\quad \forall L\in N(K).
\end{align}
The Sharp--$GN(K)$-limiter has correction factors computed via the bounds
\begin{align}
[m^{\sharp}_{\sigma_{KL}},M^{\sharp}_{\sigma_{KL}}]= \left[\max \left\lbrace \min_{i\in GN(K)}\bar{u}_{i}^{n} , \min_{i\in GN(L)}\bar{u}_{i}^{n}\right\rbrace ,\min\left\lbrace\max_{i\in GN(K)}\bar{u}_{i}^{n} , \max_{i\in GN(L)}\bar{u}_{i}^{n} \right\rbrace \right],\quad \forall L\in N(K).
\end{align}
One can by basic set inclusions and properties of the maxima and minima show that $[m^{BJ}_{K},M^{BJ}_{K}]\subseteq [m_{\sigma_{KL}},M_{\sigma_{KL}}]\subseteq [m^{\sharp}_{\sigma_{KL}},M^{\sharp}_{\sigma_{KL}}]$ sufficient 
via \cref{lemma:correction factor inequalities} to conclude that $\forall q\in \sigma_{KL}^{I}$, $\forall L\in N(K)$, $\forall K\in \mathcal{M}$,
\begin{align}
    \alpha(m^{BJ}_{K},M^{BJ}_{K},\bar{u}_K,p_{K}(\b x_q)) \leq \alpha(m_{\sigma_{KL}},M_{\sigma_{KL}},\bar{u}_K,p_{K}(\b x_q)) \leq \alpha(m^{\sharp}_{\sigma_{KL}},M^{\sharp}_{\sigma_{KL}},\bar{u}_K,p_{K}(\b x_q)).
\end{align}
Similarly via \cref{lemma:correction factor inequalities}, one can conclude for non-flux contributing points in the cell mean decomposition ($\forall q\in K^{nfc} ,\forall K\in \mathcal{M}$) one has
\begin{align}
\alpha(m^{BJ}_{K},M^{BJ}_{K},\bar{u}_K,p_{K}(\b x_q)) \leq \alpha(m_K,M_{K},\bar{u}_K,p_{K}(\b x_q)) = \alpha(m^{\sharp}_{K},M^{\sharp}_{K},\bar{u}_K,p_{K}(\b x_q)).
\end{align}
Taking the minimum over all (flux and non-flux contributing points) points in a Zhang-acceptable decomposition of the cell average gives \cref{eq:GBJ_GN_SGN}.
\end{proof}

In summary, we have proven the following \textbf{\emph{mesh-independent}} results.  Given an arbitrary order FV-scheme with a Zhang-acceptable decomposition of the cell mean, we have three classes of limiters capable of preserving discrete local maximum principles. This includes a natural generalisation of the Barth-Jespersen limiter, a Sharp--$GN(K)$ limiter family defined by the sufficient conditions of \cref{thm: slope limiters}, and a slightly relaxed $GN(K)$ limiter family where more interpretable face defined bounds are used. From here, we demonstrate that given any Generalised-Barth-Jespersen-type limiter (including the original), one can always define a $GN(K)$-limiter with less severe correction factors, preserving the same cell mean maximum principle. Similarly, given any $GN(K)$-limiter one can always define a Sharp--$GN(K)$-limiter with less severe correction factors, preserving the same cell mean maximum principle. Furthermore, the Sharp--$GN(K)$-limiter, and the $GN(K)$-limiter, can preserve a more general class of local maximum principles than Generalised Barth-Jespersen-type limiters, and do so with a polynomial representation closer to the higher order polynomial, and evolve the cell mean with a flux closer to the higher order flux.

Whilst the Sharp-$GN(K)$ limiter variant is provably more accurate (less limiting) on all meshes and preserves the same cell mean maximum principle as the GBJ-$BJ(K)$ limiter or the $GN(K)$ limiter, in practice there are often coding limitations, and operational costs. For large enough stencil size all limiters become equivalent and tend towards the same global boundedness limiter (GlobalMP) introduced in \cite{zhang2010maximum}. Therefore given a large enough stencil size, differences between limiters provably diminish and the improvement in accuracy from less severe correction factors will eventually diminish, and computational cost will eventually be a dominant factor in efficiency. The GBJ-$BJ(K)$ limiter family is the easiest conceptually and the most similar to existing limiting techniques making it more likely to be implemented, the GBJ-$BJ(K)$ limiter also admits a computational refactor which allows for a slightly faster implementation which would offset the decrease in accuracy (provided a large enough stencil size is chosen for the DLMP). The $GN(K)$-limiter family uses slightly more natural face defined bounds, and we will later show similar performance to its corresponding sharp variant (for large enough stencil size). In summary for small stencils, one should use sharp variants where the correction factors play a dominant source of error, for medium stencil sizes all three limiter family's can be competitive and hardware/software specific testing is required. For large enough maximum principles, use a GBJ type (or preferably just a global maximum principle limiter if the notion of local has grown too weak to be useful in controlling oscillations).

\subsection{Instantiated local boundedness limiters}\label{sec:instanciations_of_limiters}
We use the results established in \cref{thm: slope limiters}, to create specific local maximum principle limiters upon local stencil choices of $GN(K)$, and $BJ(K)$. We group limiters capable of preserving the same cell mean maximum principle on cell mean values. These are instantiations of the general limiters: \cref{method:sharp_m_K_M_K}, \cref{method:Dull_m_K_M_K}, \cref{method:GBJ}, with specific stencils used in later numerical experiments.

\begin{method}[Sharp--$N(K)\cup K$ limiter]\label{method:sharp_NKUK}
\noindent
\Cref{method:sharp_m_K_M_K}, with $GN(K)=N(K)\cup K$.
\end{method}
\begin{method}[$N(K)\cup K$ limiter]\label{method:NKUK}
\noindent
\Cref{method:Dull_m_K_M_K}, with $GN(K)=N(K)\cup K$.
\end{method}
\begin{method}[GBJ--$K$ limiter]\label{method:GBJ_K}
\Cref{method:GBJ} with $BJ(K)=K$.
\end{method}

\begin{corollary}
The FV scheme (\cref{eq:dg}) with either the Sharp--$N(K)\cup K$ limiter (\cref{method:sharp_NKUK}), the 
$N(K)\cup K$ limiter (\cref{method:NKUK}), or the GBJ-${K}$ limiter (\cref{method:GBJ_K}), with correction factors $\alpha_{K}^{\sharp}$, $\alpha^{GN}_{K}$, $\alpha_{K}^{GBJ}$ respectively have a resulting cell mean value satisfying the same discrete local maximum principle,
\begin{align}
\min_{L \in N(K) \cup K} \bar{u}^n_{L}\leq \bar{u}_{K}^{n+1} \leq \max_{L \in N(K) \cup K}\bar{u}^n_{L},
\end{align}
irrespective of chosen limiter.
However, the correction factors satisfy
\begin{align}
  \alpha_{K}^{GBJ}  \leq \alpha_{K}^{GN}  \leq \alpha_{K}^{\sharp}, \quad \forall K\in \mathcal{M}.\label{eq:GBJ_GN_SGN_NKUK}
\end{align} 
\end{corollary}
\begin{proof}
Strict subcase of \cref{thm: comparing limiters GBJ-SharpGN-GN}, where $BJ(K)=K$, $GN(K)=N(K)\cup K$.
\end{proof}

\begin{method}
[Sharp--$N^2(K)\cup N(K)$ limiter]\label{method:sharp_N2KUNK}
\noindent
\Cref{method:sharp_m_K_M_K}, with $GN(K)=N^2(K)\cup N(K)$.
\end{method}
\begin{method}
[$N^2(K)\cup N(K)$ limiter]\label{method:dull_N2KUNK}
\noindent
\Cref{method:Dull_m_K_M_K}, with $GN(K)=N^2(K)\cup N(K)$.
\end{method}
\begin{method}
[GBJ-$N(K)\cup K$ limiter]\label{method:GBJ_NKUK}
\noindent
\Cref{method:GBJ}, with $BJ(K)=N(K)\cup K$. 
\end{method}

\begin{corollary}Let \cref{eq:dg} denote the FV scheme and let conditions in \cref{thm: slope limiters} hold. If the FV scheme is limited via either
the Sharp--$N^2(K)\cup N(K)$ limiter (\cref{method:sharp_N2KUNK}), the $N^2(K)\cup N(K)$ limiter (\cref{method:dull_N2KUNK}), or the GBJ-$N(K)\cup K$-limiter \cref{method:GBJ_NKUK}, with correction factors $\alpha_{K}^{\sharp}$, $\alpha^{GN}_{K}$, $\alpha_{K}^{GBJ}$ respectively, then the solution attains the discrete local maximum principle irrespective of chosen limiter
\begin{align}
\min_{L \in N^2(K) \cup N(K)} \bar{u}^n_{L}\leq \bar{u}_{K}^{n+1} \leq \max_{L \in N^2(K) \cup N(K)}\bar{u}^n_{L}.
\end{align}
The correction factors satisfy
\begin{align}
\alpha_{K}^{GBJ}  \leq \alpha_{K}^{GN}  \leq \alpha_{K}^{\sharp}, \quad \forall K\in \mathcal{M}.
\end{align} 
\end{corollary}
\begin{proof}
Strict subcase of \cref{thm: comparing limiters GBJ-SharpGN-GN} where $BJ(K)=N(K)\cup K$, $GN(K)=N^2(K)\cup N(K)$.
\end{proof}
\begin{remark}
    In the strict subcase that the FV scheme is second order with a linear subcell representation in each cell the GBJ-$GN(K)$-limiter with $GN(K)=N(K)\cup K$ is exactly the Barth-Jespersen limiter, and we have proven that the Sharp--$N^2(K)\cup N(K)$ limiter (\cref{method:sharp_N2KUNK}), and the $N^2(K)\cup N(K)$ limiter (\cref{method:dull_N2KUNK}) perform less limiting than the BJ-limiter with the same cell mean maximum principle.
\end{remark}

\begin{method}
[Sharp--$VN(K)$ limiter]\label{method:sharp_VNK}
\noindent
\Cref{method:sharp_m_K_M_K}, with $GN(K)=VN(K)$.
\end{method}

\begin{method}
    [$VN(K)$ limiter]\label{method:dull_VNK}
\noindent
\Cref{method:Dull_m_K_M_K}, with
$GN(K)=VN(K)$.
\end{method}

\begin{corollary}\label{thm:HO_NVK}
The FV scheme (\cref{eq:dg}) with either the Sharp-$VN(K)$-limiter (\cref{method:sharp_VNK}) or the $VN(K)$ limiter (\cref{method:dull_VNK}) with respective correction factors $\alpha_{K}, \alpha^{\sharp}_{K}$ has a resulting cell mean value $\bar{u}_{K}^{n+1}$ satisfying the discrete local maximum principle
	\begin{align}
		\min_{L \in VN(K)} \bar{u}^n_{L}\leq \bar{u}_{K}^{n+1} \leq \max_{L \in VN(K)}\bar{u}^n_{L}. \label{eq:statement2}
	\end{align}
The correction factors satisfy
\begin{align}
  \alpha_{K}  \leq \alpha^{\sharp}_{K}.
\end{align}
\end{corollary}
\begin{proof}
Strict subcase of \cref{thm: comparing limiters GBJ-SharpGN-GN} where $GN(K)=VN(K)$.
\end{proof}

\begin{method}[Global boundedness limiter \cite{zhang2010maximum}]
\cref{method:Dull_m_K_M_K}, \cref{method:sharp_m_K_M_K}, \cref{method:GBJ} with $BJ(K)$, and $GN(K)$ covering the entire mesh $\mathcal{M}$.
\end{method}

\begin{theorem}
For a second order finite volume method with linear subcell representation, the FV2 scheme (\cref{eq:dgFV2}) with either the Kuzmin limiter \cite{kuzmin2010vertex} (equivalently MLP-\cite{park2010multi} limiter), the Sharp-$VN(K)$-limiter (\cref{method:sharp_VNK}) or the $VN(K)$ limiter (\cref{method:dull_VNK}) with respective correction factors $\alpha^{kuz}_{K}, \alpha_{K}, \alpha^{\sharp}_{K}$ have a resulting cell mean value $\bar{u}_{K}^{n+1}$ satisfying the same vertex neighbour discrete local maximum principle,
	\begin{align}
		\min_{L \in VN(K)} \bar{u}^n_{L}\leq \bar{u}_{K}^{n+1} \leq \max_{L \in VN(K)}\bar{u}^n_{L}. 
	\end{align}
The correction factors satisfy
\begin{align}
  \alpha^{Kuz}_{K}  \leq \alpha_{K}  \leq \alpha^{\sharp}_{K}.
\end{align}

\end{theorem}
\begin{proof}
This is true for all 2D and 3D admissible meshes, but the proof is not entirely a consequence of previous theory. Instead, a constructive and dimension dependent proof is required using specific properties about affine functions on convex polygons and convex polyhedra and so is relegated to a section in the appendix, with the required notational setup in \cref{sec:fv2}.
\end{proof}

\begin{remark}
We do not compare limiters which enforce different maximum principles, as results become mesh dependent and different maximum principles are more natural on different meshes. For instance, on a uniform square mesh by set inclusion properties one can prove that the correction factors are less severe for the $N^2(K)\cup N(K)$-limiter than the $VN(K)$-limiter. However, on a uniform equilateral  triangle mesh one can prove the exact opposite result. 
\end{remark}

\section{Application 1: Second order finite volume (FV2)}\label{sec: Application: second order method: chapter 2}
We consider conditions for a second order finite volume scheme (introduced for arbitrary meshes in 2D and 3D in \cref{sec:fv2}) to have a local maximum principle on a uniform square mesh, of cell width $\Delta x$ and height $\Delta y$ respectively. We do so in the context of theory and limiters introduced in \cref{sec: Theory: theorems: chapter 2,sec:instanciations_of_limiters} which simplifies significantly since there is a Zhang-acceptable decomposition of the cell mean onto solely flux contributing facet midpoints in both 2D and 3D (\cref{sec:fv2}). For a second order finite volume scheme, the interpolating polynomial is a linear subcell representation
\begin{align}
p_{i,j}(x,y) &= \bar{u}_{i,j} + \alpha (u_x)_{{i,j}}(x-x_i) + \alpha(u_y)_{i,j}(y-y_j),\\
			    	(u_x)_{{i,j}} & = \frac{\bar{u}_{i+1} - \bar{u}_{i-1} }{2\Delta x}, \quad (u_y)_{{i,j}} = \frac{\bar{u}_{j+1} - \bar{u}_{j-1} }{2\Delta y},
\end{align}
	where $\alpha$ arises from the slope limiter. This subcell representation satisfies the conservation property 
 \begin{align}
 \frac{1}{\Delta x \Delta y}\int_{x_{i-1/2}}^{x_{i+1/2}}\int_{y_{j-1/2}}^{y_{j+1/2}} p_{i,j}(x,y) dx dy = \bar{u}_{i,j}.
 \end{align} The flux contributing quadrature points are at the midpoint of each face, and the quadrature point evaluations for cell $(i,j)$ are the right left up and down values defined as
	\begin{align}
	u^R_{i,j} &= p_{i,j}(x_{i+1/2},y_j),\quad
	u^L_{i,j} = p_{i,j}(x_{i-1/2},y_j),\quad
	u^U_{i,j} = p_{i,j}(x_{i},y_{j+1/2}),\quad 
	u^D_{i,j} = p_{i,j}(x_{i},y_{j-1/2}).
	\end{align}
	A Zhang-acceptable decomposition of the cell average can be conveniently found in terms of solely flux contributing quadrature points as $\bar{u}_{i,j} = \frac{1}{4}(u^R_{i,j}+u^{L}_{i,j}+u^{U}_{i,j}+u^{D}_{i,j})$. 
	The right-hand flux through face $(i+1/2,j)$ is computed by the midpoint, $\frac{| \Delta y|}{| \Delta x \Delta y|}F_{i,i+1}(u^R_{i,j},u^{L}_{i+1,j},\b v \cdot n_{i,i+1})$.
So that the method can be written as the sum of 4 local HHLK schemes solving Riemann problems
 \begin{align}
\operatorname{HHLK}_{i+1/2,j}(u^{R}_{i,j},u^L_{i+1,j},\b v \cdot n_{i,i+1})= \frac{1}{4} \left[ u^{R}_{i,j} - \frac{4\Delta t}{\Delta x} F(u^{R}_{i,j},u^L_{i+1,j}, \b v \cdot n_{i,i+1}) \right],
\end{align}
there are only face-defined flux-contributing quadrature points located at the midpoint of each face. 
Therefore the cell mean evolution equation for cell $(i,j)$,
			    \begin{align}
			    	\bar{u}_{i,j}^{n+1} = \bar{u}_{i,j}^n - \frac{\Delta t  }{ \Delta x } F_{i,i+1} - \frac{\Delta t  }{ \Delta x }F_{i,i-1}- \frac{\Delta t  }{ \Delta y } F_{j,j+1} - \frac{\Delta t  }{ \Delta y }F_{j,j-1},
			    \end{align}
			    is a monotonic function of the edge defined quadrature points $u^R_{i,j}$, $u^{L}_{i+1,j}$,$ u^L_{i,j}$, $u^{R}_{i-1,j}$, $u^U_{i,j}$, $u^{D}_{i,j+1}$, $u^D_{i,j}$, $u^{U}_{i,j-1}$, when the following local Courant number conditions holds $\frac{\Delta t}{\Delta x} (\b v \cdot n_{i,i+1})^+,  \frac{\Delta t}{\Delta x} (\b v \cdot n_{i,i-1})^+, \frac{\Delta t}{\Delta y} (\b v \cdot n_{j,j+1})^+, \frac{\Delta t}{\Delta y} (\b v \cdot n_{j,j-1})^+ \leq \frac{1}{4}.$
		    
		    	The Courant number is now a concept to be interpreted on edges
		    	\begin{align}
		    	\frac{\Delta t|\sigma_{KL}| (\b v \cdot n_{KL})^+}{|K|} \leq 1/4, \quad \forall K\in \mathcal{M}, \quad \forall L\in N(K).
		    	\end{align}
			    Though one can pessimistically write this in terms of a more convenient (less sharp) cell defined Courant number as 
			    \begin{align}
			    	C_{K} = \sum_{L\in N(K)} \frac{\Delta t|\sigma_{KL}| (\b v \cdot n_{KL})^+}{|K|} \leq 1/4, \quad \forall K \in \mathcal{M},
			    \end{align}
			    where if one assumes incompressibility $1/4$ becomes $1/2$.
We demonstrated the FV2 method can be written as a monotone function of quadrature points. We could have equivalently used \cref{thm: slope limiters} with $w^{\sigma_{KL}}_q = 1$, $\gamma^{K}_{q,L} = 1/4$ and identified that; there are no non-flux contributing quadrature points, and only one flux contributing quadrature point per face. This generalises into 2D and 3D on arbitrary meshes (\cref{sec:fv2}).

\subsection{Application 1a: Advection}\label{sec:fv2_advection}
In this application, we study solid body rotation of several profiles of differing regularity, consisting of a cone, a cosine hump and a slotted cylinder. This is considered a challenging test of monotonicity and shape preservation \cite{Z_1979}, \cite{leveque_2002}. 
We define the numerical domain as $\Omega = [0,1]^2$, the flow as solid body rotation defined by the stream function $\Psi = (y-x)$ and the LeVeque initial condition \cite{Leveque_test_ic} is defined as
\begin{align}
	u_{0}  = 
	\begin{cases}
			1 &\text{for} \sqrt{ (x-0.5)^2+(y-0.75)^2} \leq 0.15, \quad\text{and}\quad (x\leq 0.475), \\
			1 & \text{for} \sqrt{(x-0.5)^2+(y-0.75)^2} \leq 0.15, \quad\text{and}\quad (x>0.525), \\
			1 & \text{for}\sqrt{(x-0.5)^2+(y-0.75)^2} \leq 0.15, \quad\text{and}\quad (y\geq 0.85),\quad\text{and} \quad(0.475<x\leq 0.525), \\
			(1-\frac{r_{cone}}{0.15}) &\text{for}\quad (r_{cone} = \sqrt{(x-0.5)^2+(y-0.25)^2}\leq 0.15), \\
			\frac{1}{4}(1+\cos(\pi  \frac{r_{cos}}{0.15}) & \text{for}\quad (r_{cos} = \sqrt{(x-0.25)^2+(y-0.5)^2}\leq 0.15),\\
			0&\mbox{otherwise.} 
	\end{cases}\label{test:LeVeque} 
\end{align}

\begin{figure}[H]
\centering

% --- Row 1 (3 panels) ---
\begin{subfigure}[t]{0.32\textwidth}
  \centering
  \includegraphics[width=\linewidth]{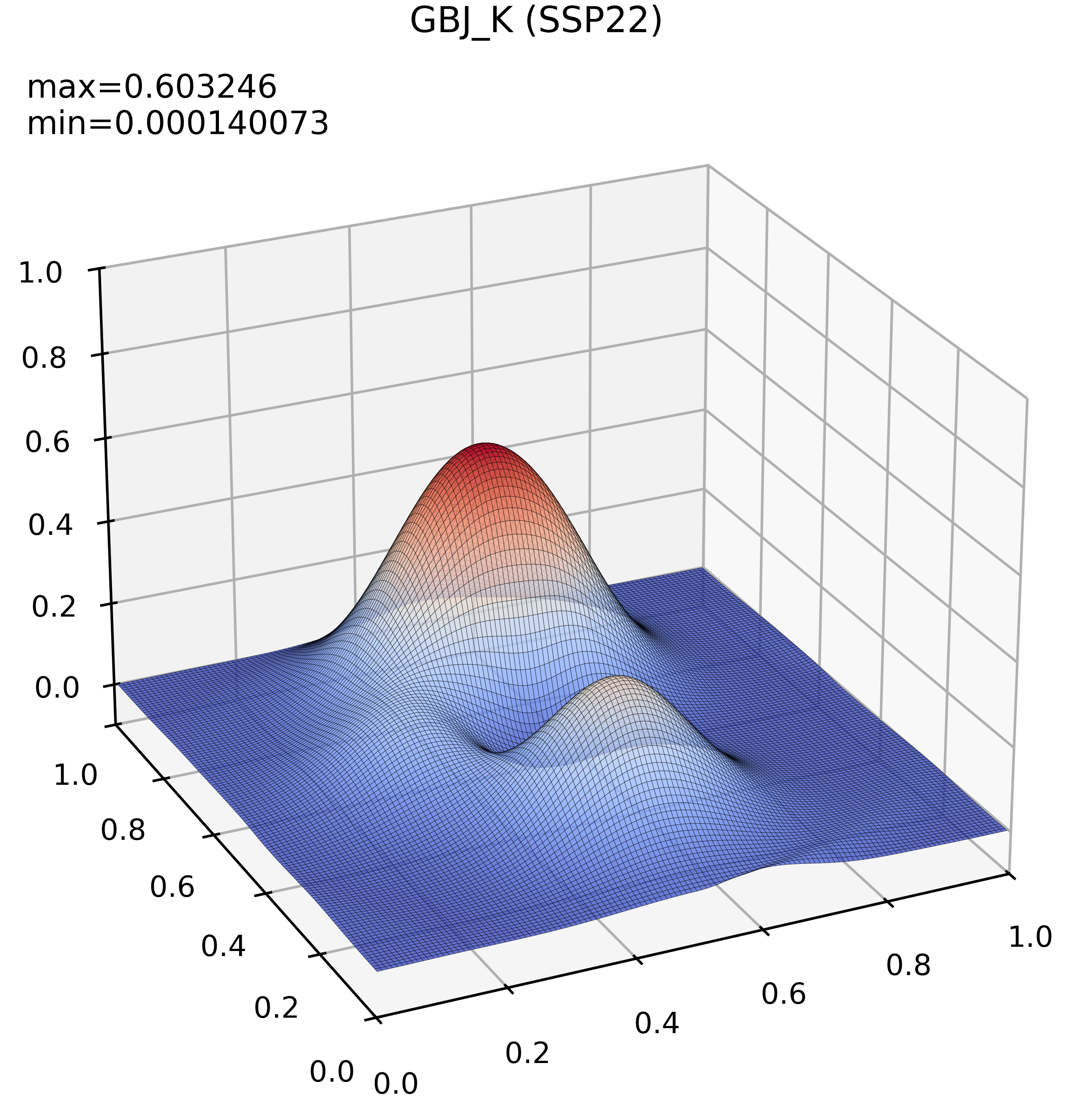}
  \caption{GBJ--$K$.}
\end{subfigure}\hfill
\begin{subfigure}[t]{0.32\textwidth}
  \centering
  \includegraphics[width=\linewidth]{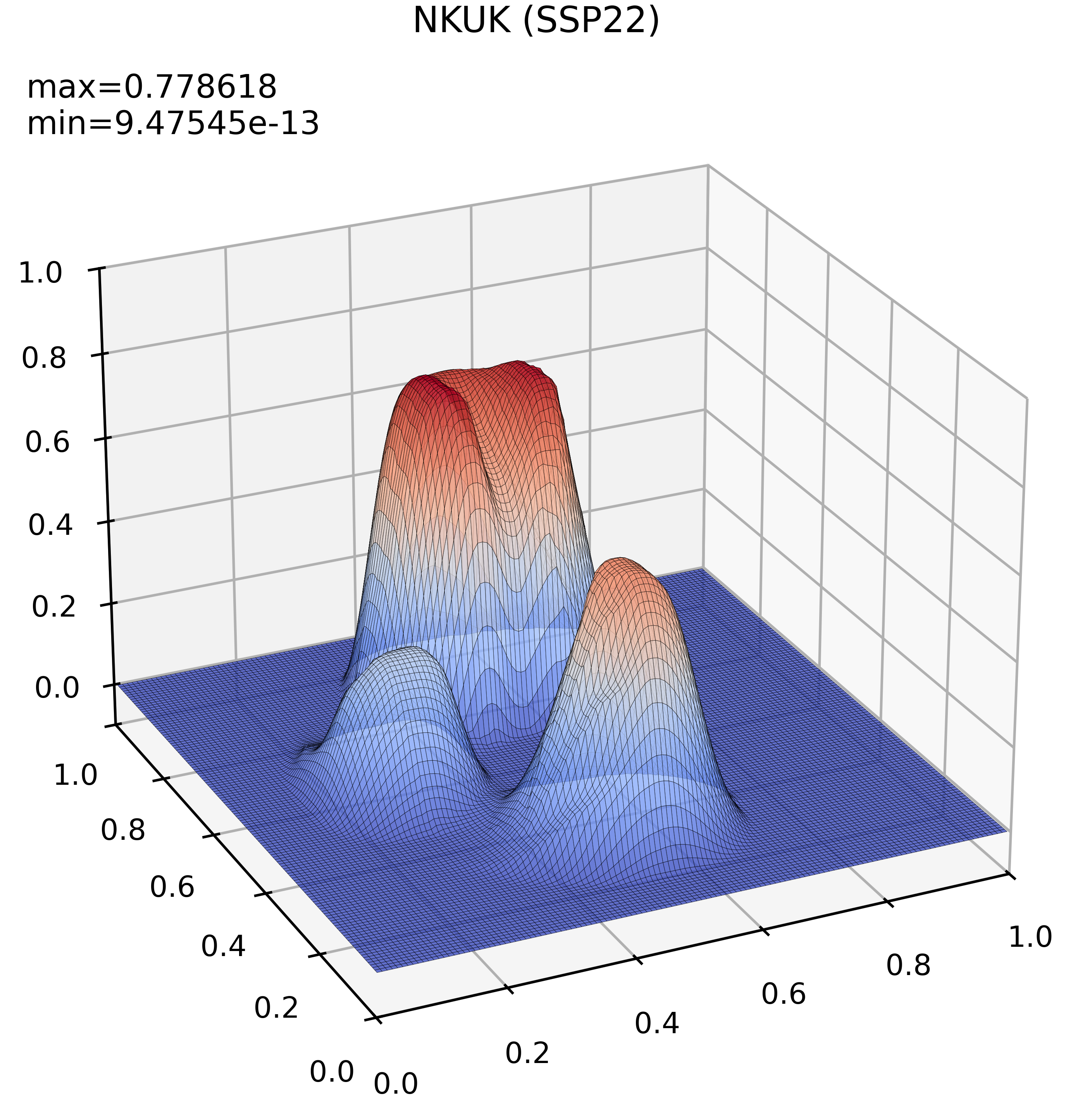}
  \caption{NKUK.}
\end{subfigure}\hfill
\begin{subfigure}[t]{0.32\textwidth}
  \centering
  \includegraphics[width=\linewidth]{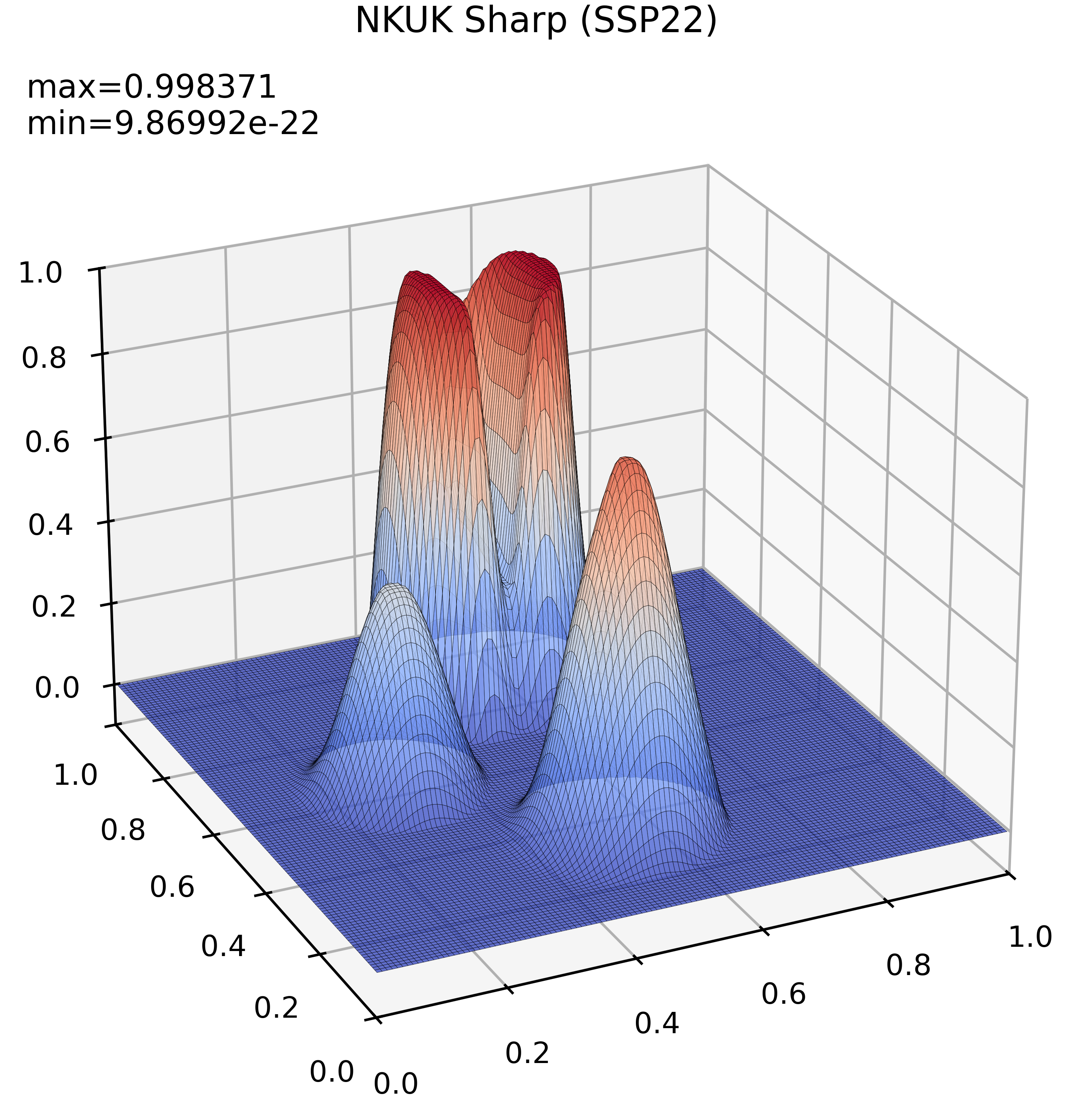}
  \caption{Sharp--NKUK.}
\end{subfigure}

\vspace{0.6em}

% --- Row 2 (3 panels) ---
\begin{subfigure}[t]{0.32\textwidth}
  \centering
  \includegraphics[width=\linewidth]{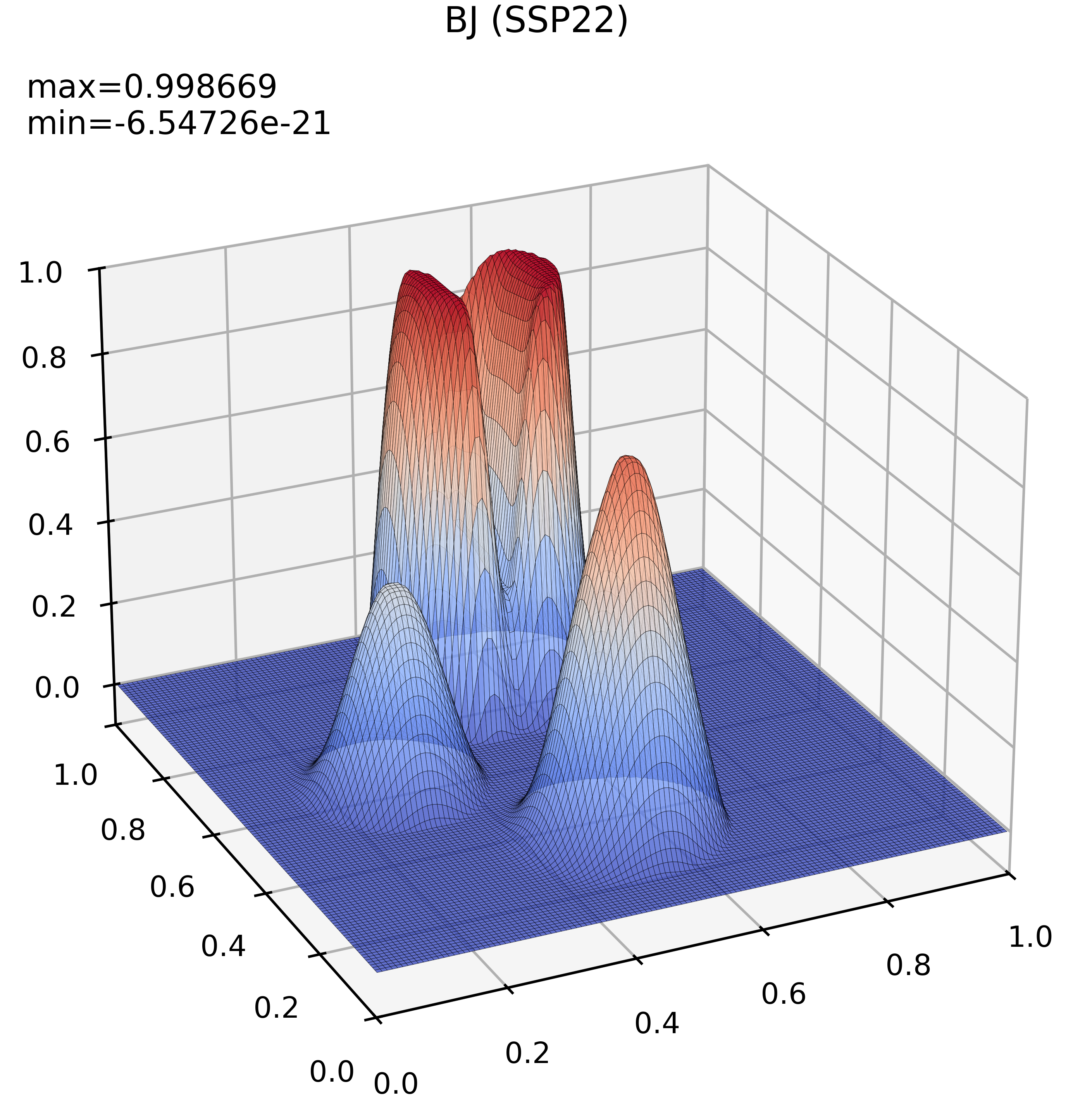}
  \caption{BJ.}
\end{subfigure}\hfill
\begin{subfigure}[t]{0.32\textwidth}
  \centering
  \includegraphics[width=\linewidth]{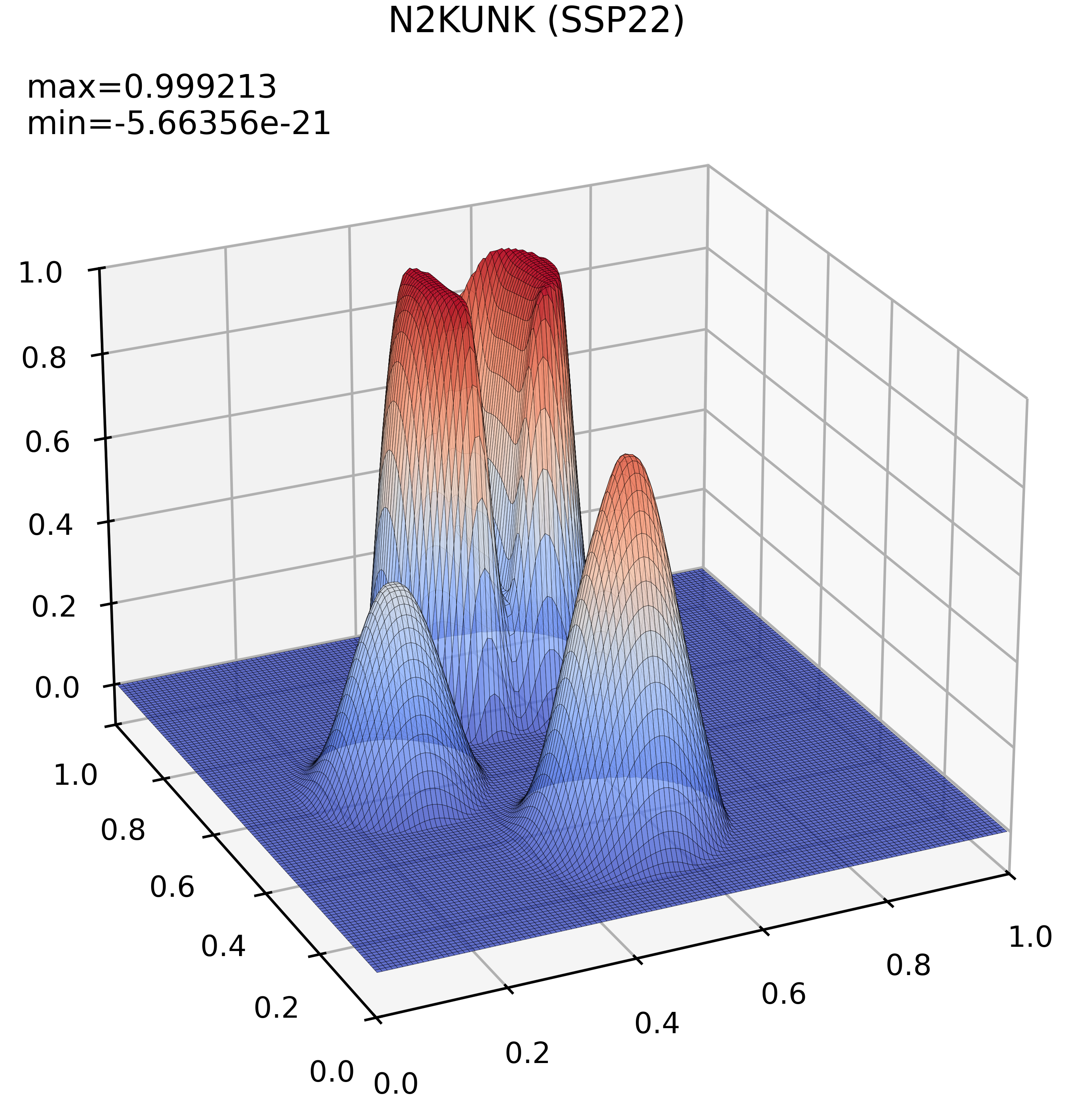}
  \caption{N2KUNK.}
\end{subfigure}\hfill
\begin{subfigure}[t]{0.32\textwidth}
  \centering
  \includegraphics[width=\linewidth]{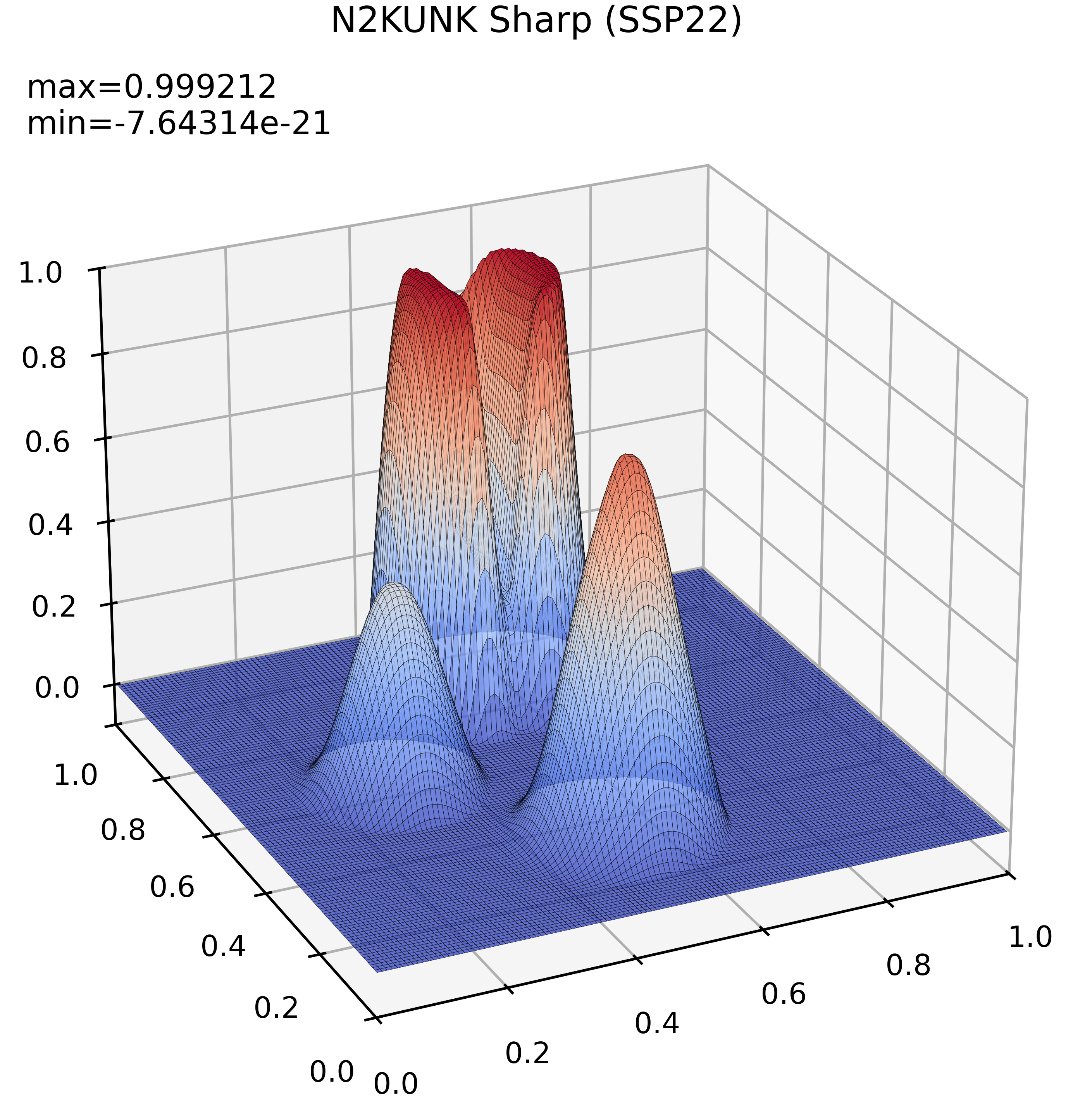}
  \caption{Sharp--N2KUNK.}
\end{subfigure}

\vspace{0.6em}

% --- Row 3 (3 panels) ---
\begin{subfigure}[t]{0.32\textwidth}
  \centering
  \includegraphics[width=\linewidth]{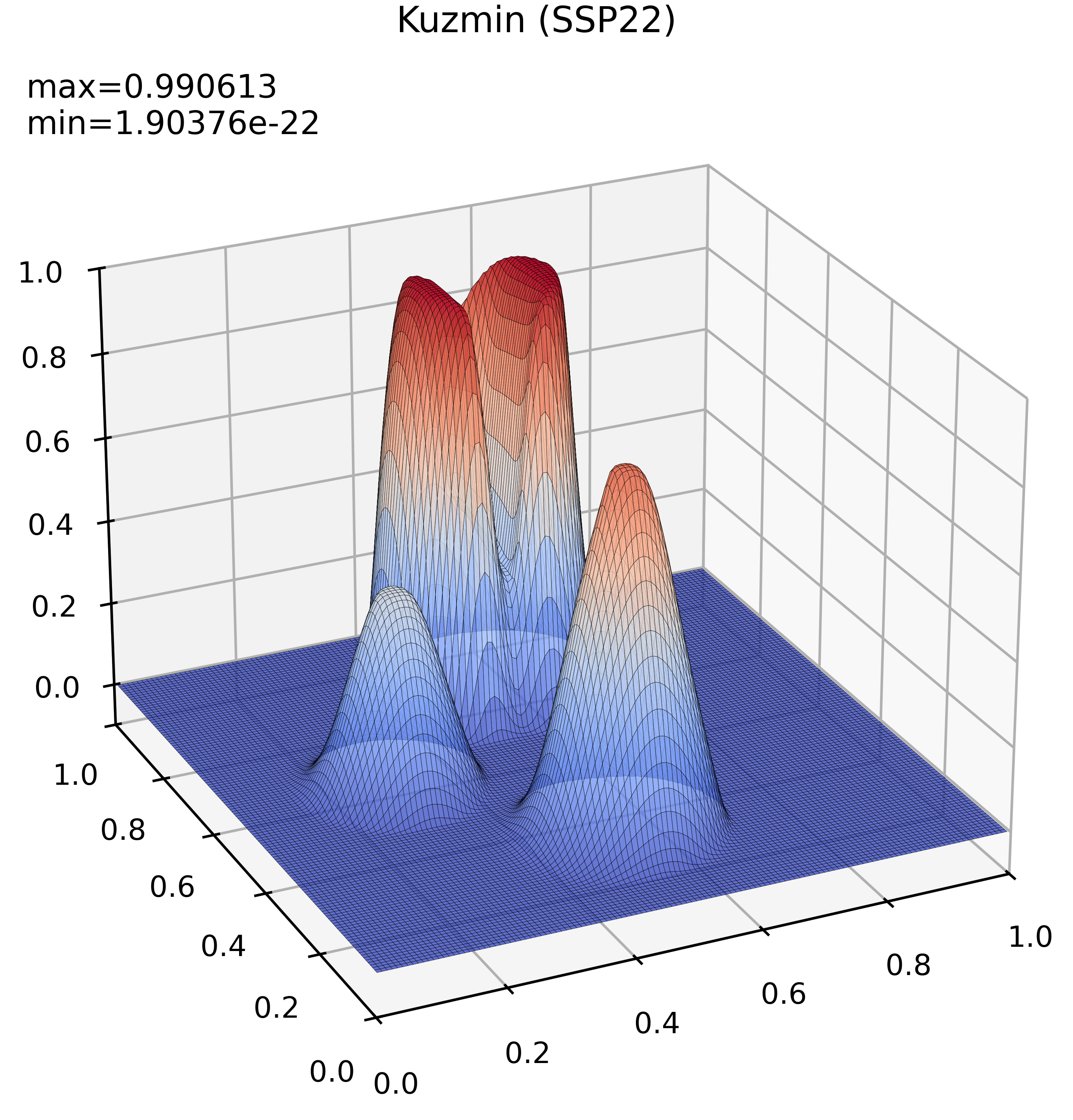}
  \caption{Kuzmin.}
\end{subfigure}\hfill
\begin{subfigure}[t]{0.32\textwidth}
  \centering
  \includegraphics[width=\linewidth]{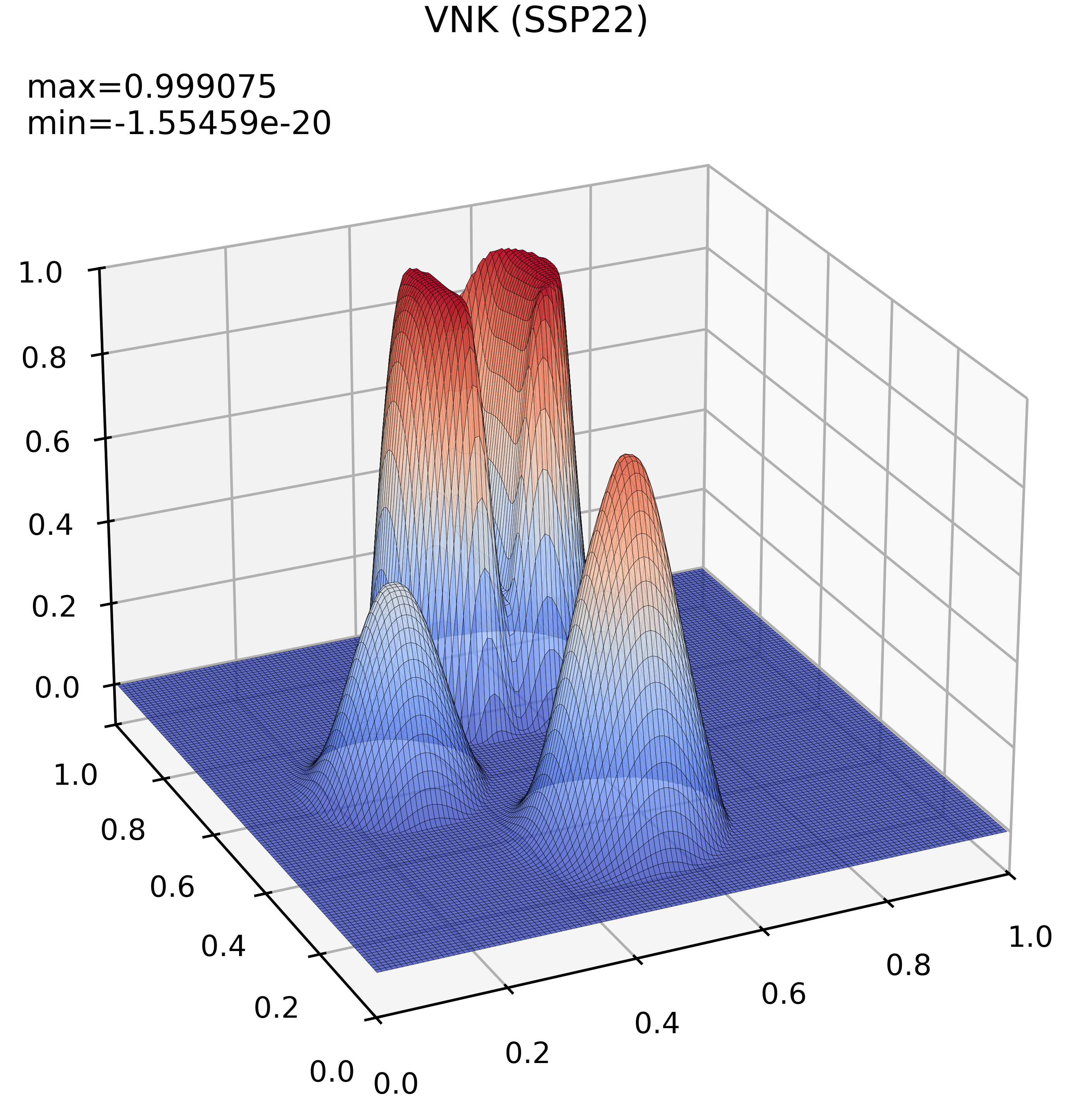}
  \caption{VNK.}
\end{subfigure}\hfill
\begin{subfigure}[t]{0.32\textwidth}
  \centering
  \includegraphics[width=\linewidth]{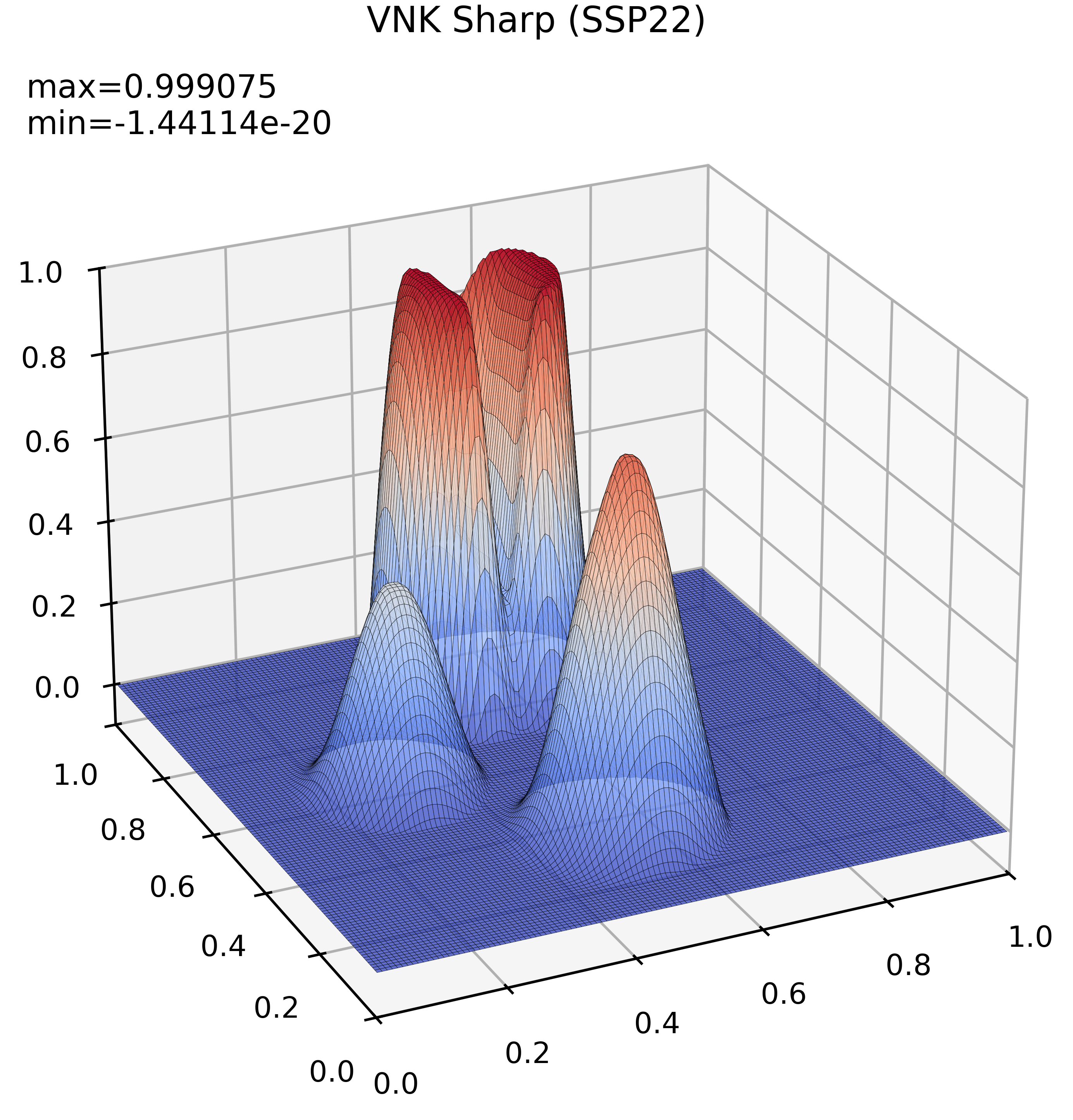}
  \caption{Sharp--VNK.}
\end{subfigure}

\caption{LeVeque solid-body rotation at $n_x=128$. 3D surface plots at final time for three DLMP families: condensed-$N(K)\cup K$-stencil (top), BJ-type $N^2(K)\cup N(K)$ (middle), and Vertex-Neighbour-$VN(K)$-based (bottom).}
\label{fig:leveque-3d-grid-128}
\end{figure}

The results are summarised as follows. 
\Cref{fig:leveque-3d-grid-128} contains 3D plots of the final timestep for three separate classes of limiters, each enforcing the same DLMP in each row. Maxima and minima are displayed in the legend of \cref{fig:leveque-3d-grid-128}, showing no global bounds are violated for any limiter beyond machine precision error. The Sharp--$N(K)\cup K$ limiter is significantly more accurate than other limiters in its corresponding $N(K)\cup K$-DLMP class. The Sharp-$VN(K)$ and $VN(K)$ limiters are slightly better at preserving peaks and marginally more accurate to those in their corresponding $VN(K)$-DLMP class. Similarly, for the larger stencil $N^2(K)\cup N(K)$-DLMP, the sharp--$N^2(K)\cup N(K)$ limiter and the $N^2(K)\cup N(K)$ limiter are marginally more accurate and peak preserving than the Barth-Jespersen limiter.

\begin{figure}[H]
\centering

\begin{subfigure}[t]{\linewidth}
  \centering
  \includegraphics[width=0.80\linewidth]{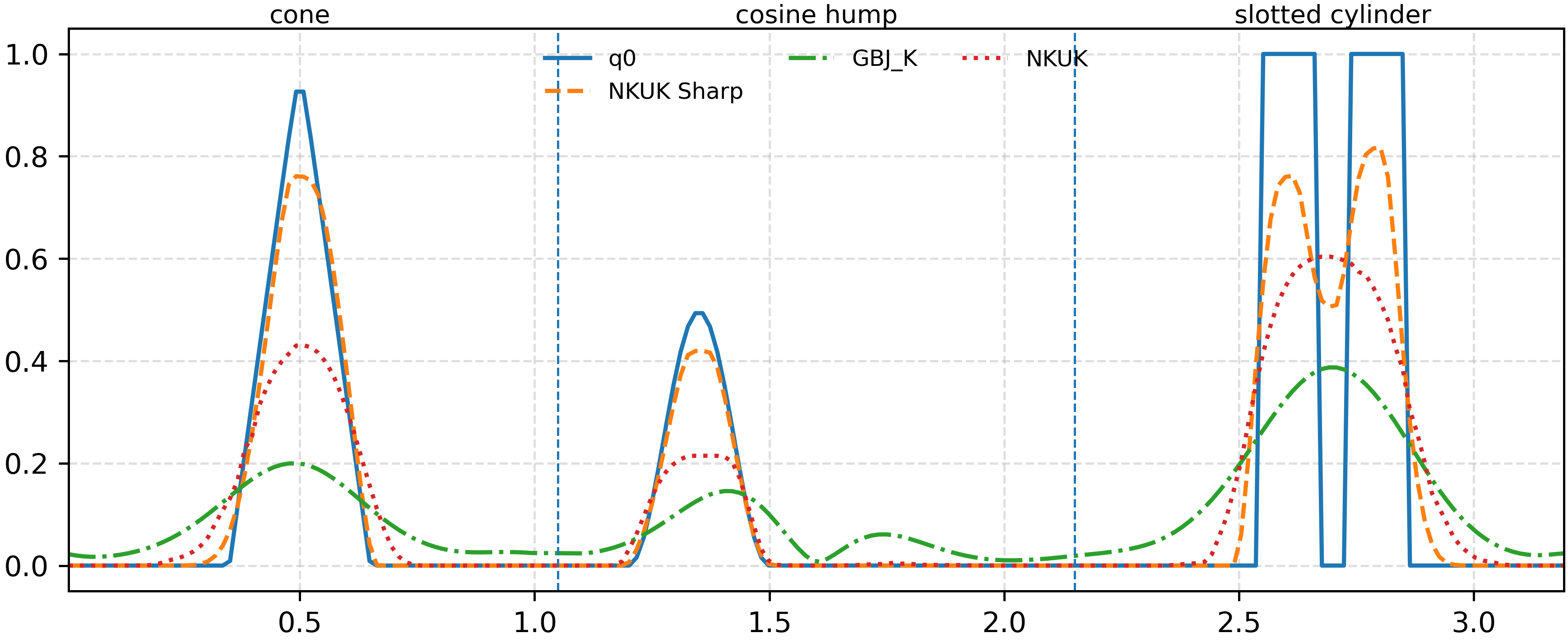}
  \caption{$N(K)\cup K$ DLMP family.}
  \label{fig:slice64_NKUK}
\end{subfigure}

\vspace{0.6em}

\begin{subfigure}[t]{\linewidth}
  \centering
  \includegraphics[width=0.80\linewidth]{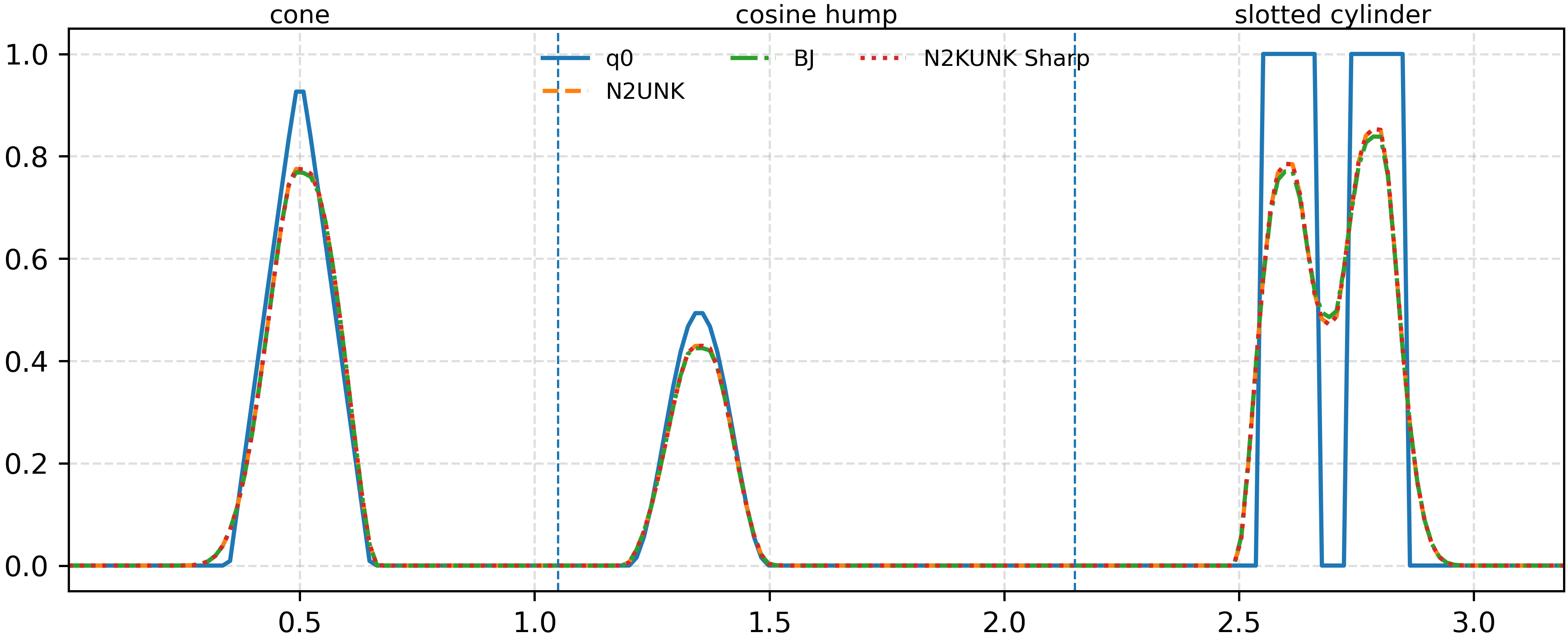}
  \caption{$N^2(K)\cup N(K)$ DLMP family (N2KUNK variants compared with BJ).}
  \label{fig:slice64_N2KUNK}
\end{subfigure}

\vspace{0.6em}

\begin{subfigure}[t]{\linewidth}
  \centering
  \includegraphics[width=0.80\linewidth]{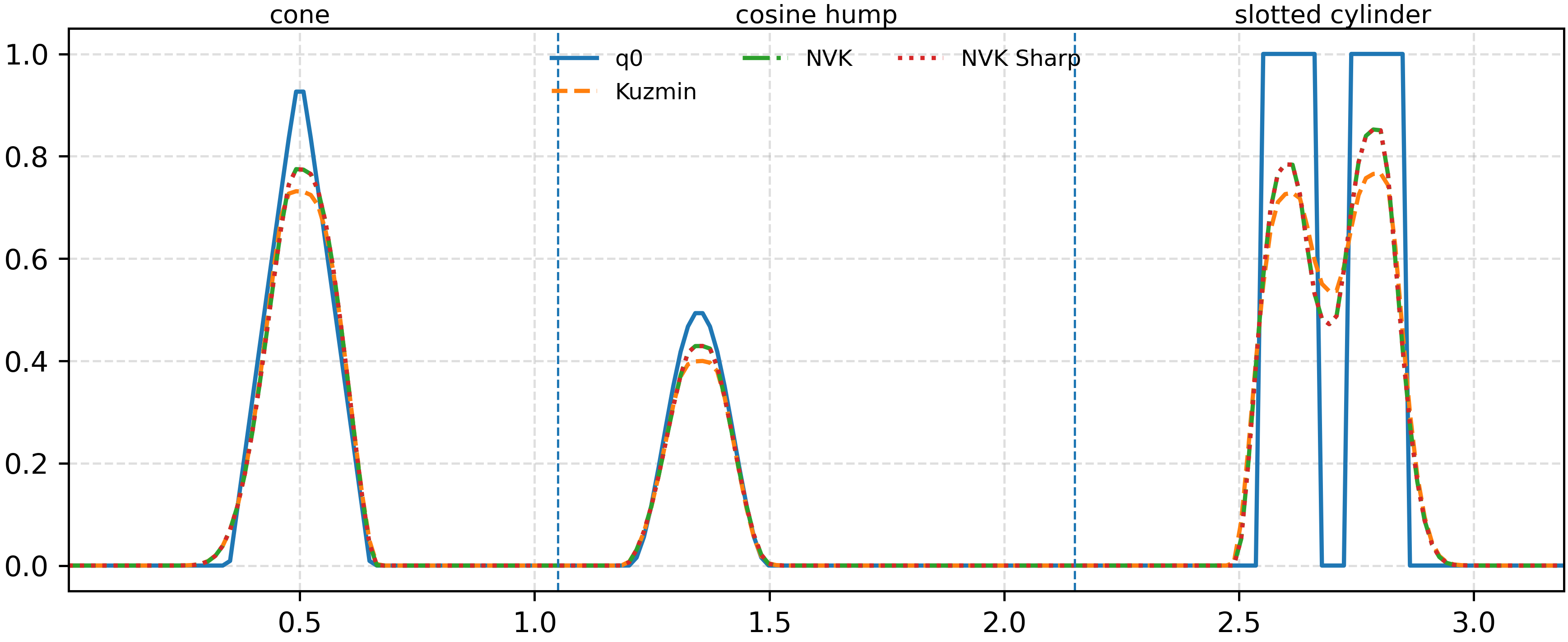}
  \caption{$VN(K)$ DLMP family (VNK variants compared with Kuzmin).}
  \label{fig:slice64_NVK}
\end{subfigure}

\caption{LeVeque test ($n_x=64$). Concatenated slices through the three profiles (cone, cosine hump, slotted cylinder) for three DLMP families.}
\label{fig:slices64_leveque_joined}
\end{figure}

\Cref{fig:slices64_leveque_joined} contains concatenated horizontal slices through shape (cone cosine slotted cylinder) midpoints at resolution $64^2$. Improvements are seen using the Sharp or $GN(K)$ class of limiter for all DLMP classes in each region (cone, cosine, slotted cylinder), and the effect is more pronounced for the more restrictive maximum principles. 

\begin{figure}[H]
\centering
\captionsetup[subfigure]{font=small}
\captionsetup{font=small}

\begin{subfigure}[t]{0.32\linewidth}
  \centering
  \includegraphics[width=\linewidth]{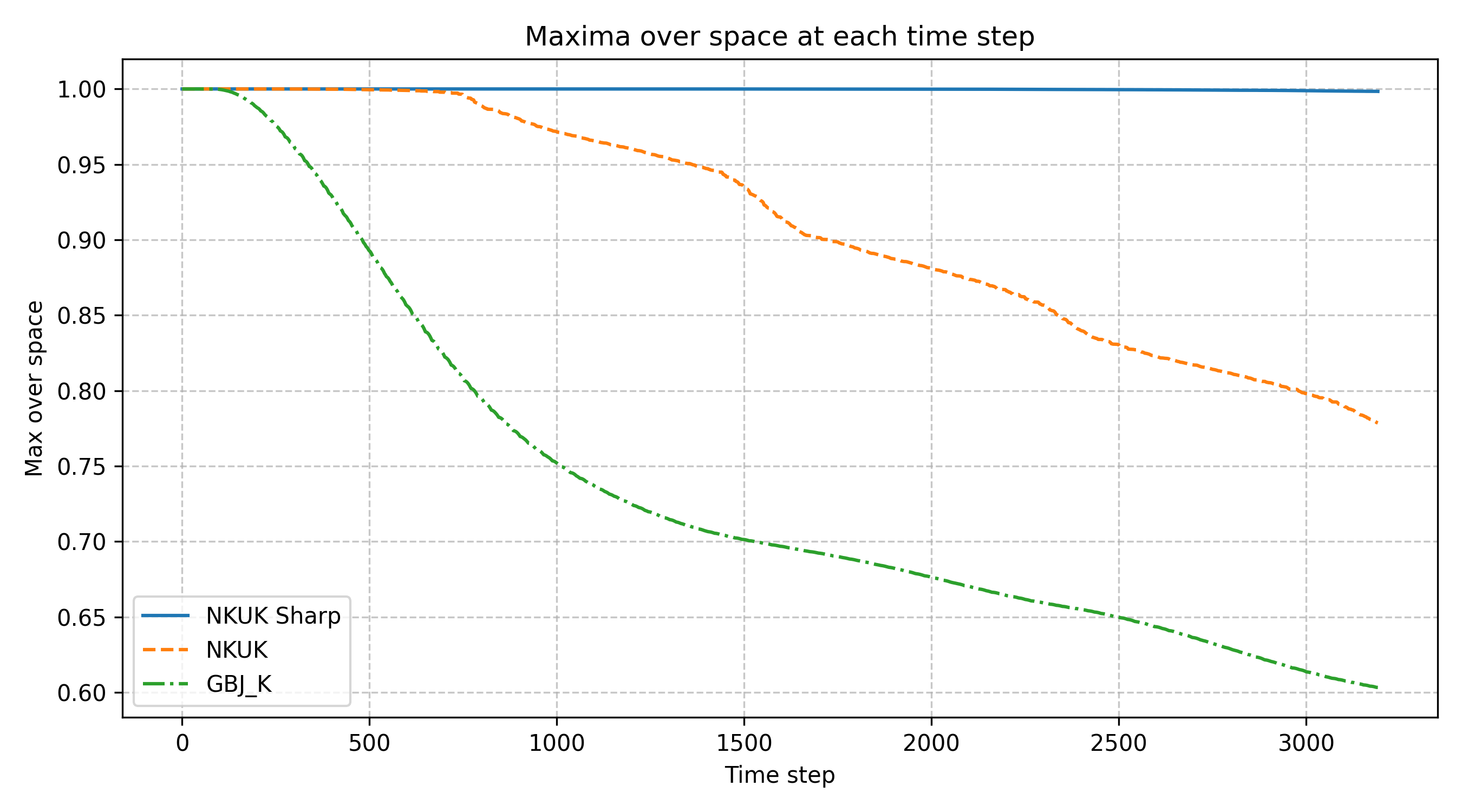}
  \caption{$N(K)\cup K$ DLMP: GBJ--$K$, NKUK, Sharp--NKUK.}
  \label{fig:leveque-maxima-timeseries-nkuk-128}
\end{subfigure}\hfill
\begin{subfigure}[t]{0.32\linewidth}
  \centering
  \includegraphics[width=\linewidth]{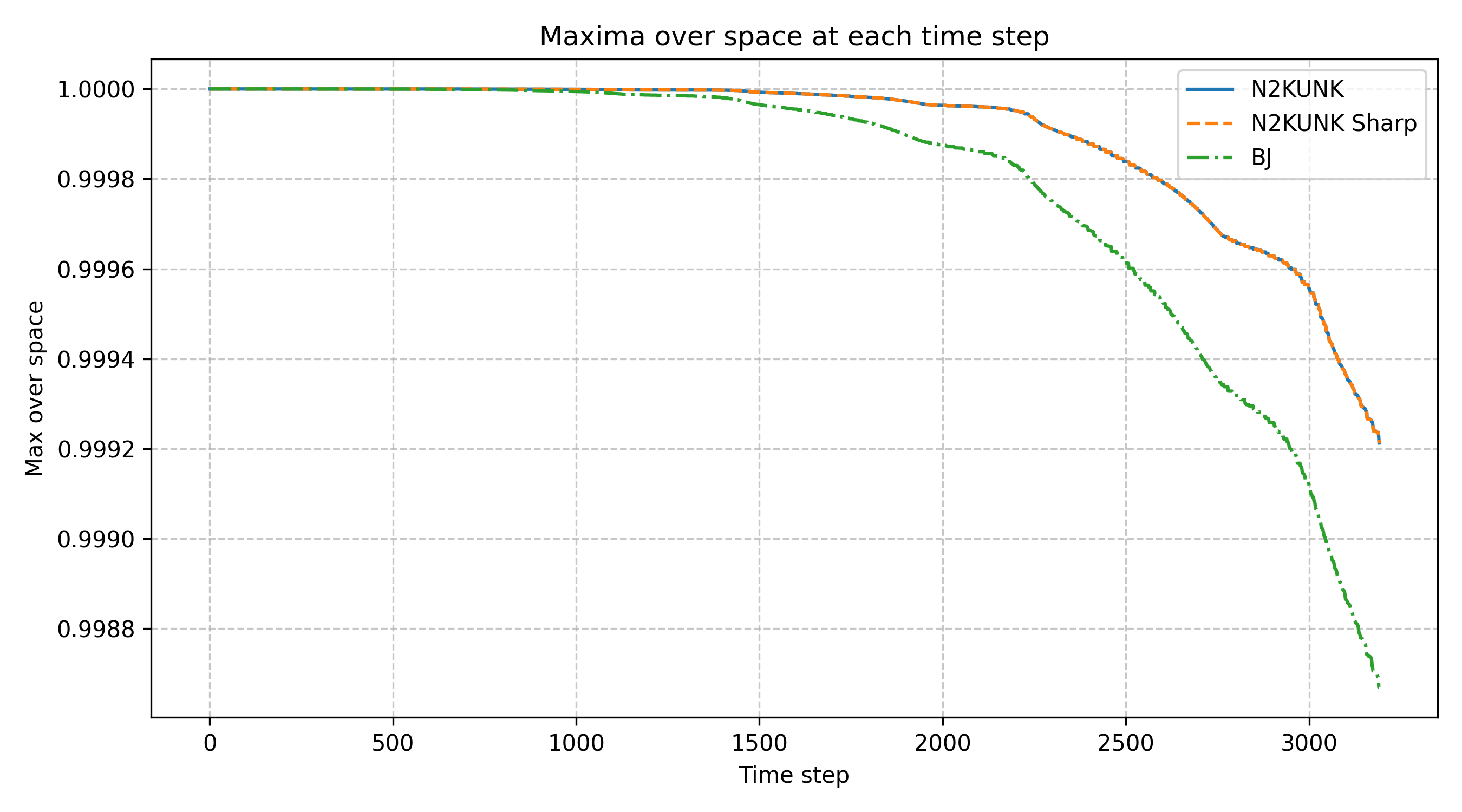}
  \caption{$N^2(K)\cup N(K)$ DLMP: BJ, N2KUNK, Sharp--N2KUNK.}
  \label{fig:leveque-maxima-timeseries-n2kunk-128}
\end{subfigure}\hfill
\begin{subfigure}[t]{0.32\linewidth}
  \centering
  \includegraphics[width=\linewidth]{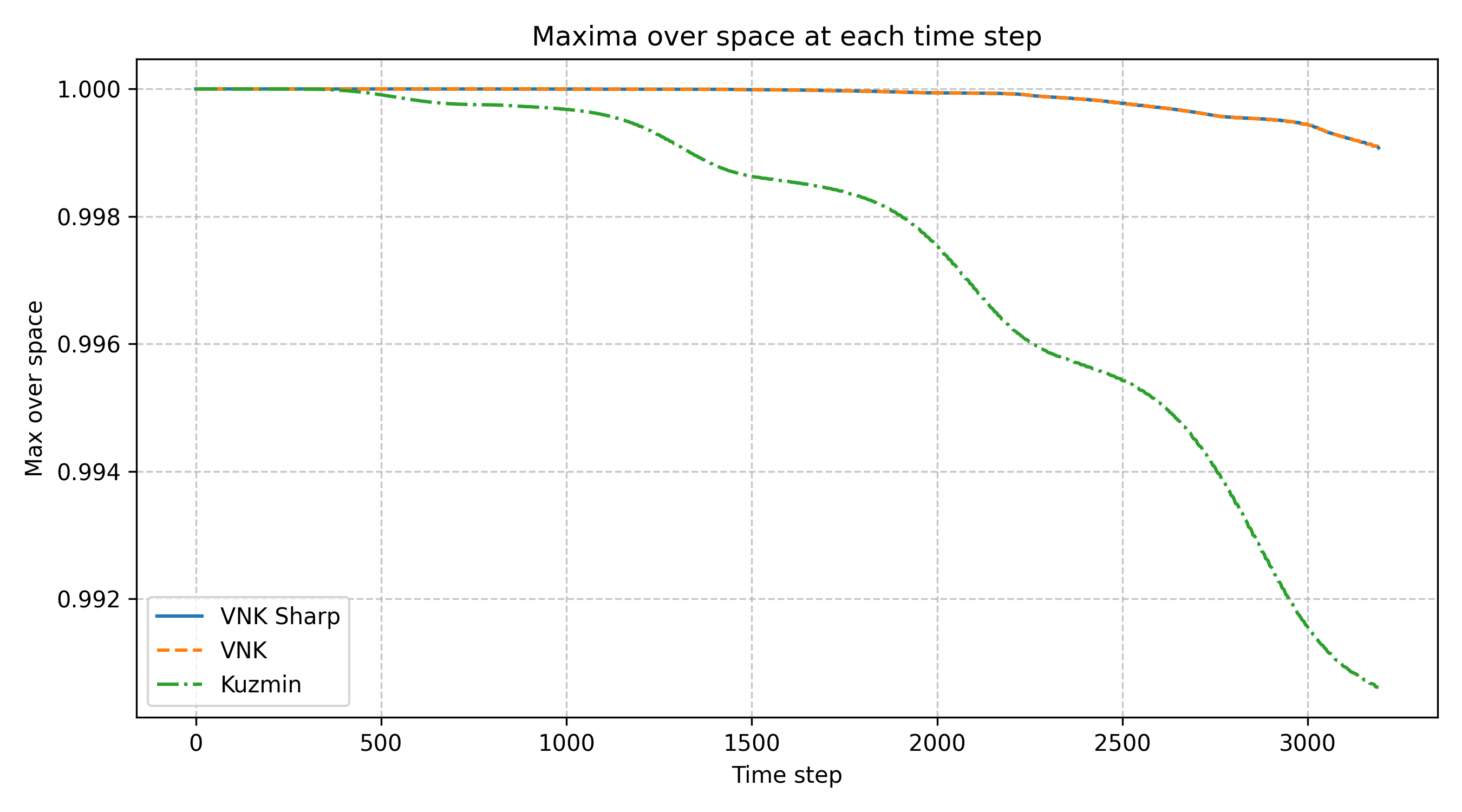}
  \caption{$VN(K)$ DLMP: Kuzmin, VNK, Sharp--VNK.}
  \label{fig:leveque-maxima-timeseries-vnk-128}
\end{subfigure}

\caption{LeVeque solid-body rotation ($n_x=128$). Time series of the spatial maximum $\max_{\Omega} q(\cdot,t)$ for the three DLMP families. We see better maxima preservation in the $GN(K)$-limiter and the Sharp $GN(K)$ limiter more so. In the second figure the BJ-limiter is more diffusive. In the third figure the Kuzmin limiter is more diffusive.}
\label{fig:leveque-maxima-timeseries-128}
\end{figure}

Similarly, we find improvements using the Sharp or $GN(K)$ class in each region and in each DLMP class in each region. These are reported via relative L2 metrics about each shape in 
\cref{tab:leveque-region-relL2-dlmp-groups}. Improvements decrease as the DLMP incorporates a larger stencil, as expected. Finally, plotting spatial maxima as a function of time to assess peak retention (\cref{fig:leveque-maxima-timeseries-vnk-128}), we see improvements using the Sharp or $GN(K)$ limiters in each DLMP class due to reduced correction factors. Improvements decrease as the DLMP incorporates a larger stencil, but do not completely diminish.

\begin{table*}[t]
\centering
\small
\setlength{\tabcolsep}{5pt}
\renewcommand{\arraystretch}{1.12}
\caption{Region-wise relative $L^2$ errors at $n_x\in\{64,128\}$, grouped by DLMP families. Rates are $p=\log_2(e_{64}/e_{128})$ (reported on the $128$ row). Within each DLMP group, the smallest error at each resolution is in bold.}
\label{tab:leveque-region-relL2-dlmp-groups}
\begin{tabular}{ll|cccc|cccc}
\hline
Method & $n_x$ & Cone & Cosine & Slotted & Whole & $p_{\mathrm{cone}}$ & $p_{\mathrm{cosine}}$ & $p_{\mathrm{slotted}}$ & $p_{\mathrm{whole}}$ \\
\hline\hline
\multicolumn{10}{l}{\textbf{FV2}}\\
\hline
No limiter
& 64  & 1.1901e-01 & 2.1059e-01 & 4.6712e-01 & 4.8851e-01 &  &  &  &  \\
& 128 & 5.0650e-02 & 8.4883e-02 & 3.3905e-01 & 3.5482e-01 & 1.232 & 1.311 & 0.462 & 0.461 \\
\hline\hline

\multicolumn{10}{l}{\textbf{DLMP group: $N(K)\cup K$}}\\
\hline
GBJ-K
& 64  & 6.7669e-01 & 7.0721e-01 & 7.2826e-01 & 7.6921e-01 &  &  &  &  \\
& 128 & 5.2782e-01 & 6.4076e-01 & 6.1734e-01 & 6.6193e-01 & 0.358 & 0.142 & 0.238 & 0.217 \\
\hline
NKUK
& 64  & 3.7659e-01 & 4.9085e-01 & 5.7783e-01 & 6.0329e-01 &  &  &  &  \\
& 128 & 1.9986e-01 & 2.5998e-01 & 4.4853e-01 & 4.6985e-01 & 0.914 & 0.917 & 0.365 & 0.361 \\
\hline
SharpNKUK
& 64  & \textbf{1.0853e-01} & \textbf{1.4960e-01} & \textbf{4.6618e-01} & \textbf{4.6846e-01} &  &  &  &  \\
& 128 & \textbf{4.6316e-02} & \textbf{4.7353e-02} & \textbf{3.1998e-01} & \textbf{3.2934e-01} & 1.228 & 1.660 & 0.543 & 0.508 \\
\hline\hline

\multicolumn{10}{l}{\textbf{DLMP group: $N^2(K)\cup N(K)$}}\\
\hline
BJ
& 64  & 1.0785e-01 & 1.5145e-01 & 4.6813e-01 & 4.6944e-01 &  &  &  &  \\
& 128 & 4.6113e-02 & 4.7237e-02 & 3.2125e-01 & 3.2987e-01 & 1.226 & 1.681 & 0.543 & 0.509 \\
\hline
N2KUNK
& 64  & \textbf{1.0637e-01} & \textbf{1.4930e-01} & \textbf{4.6527e-01} & \textbf{4.6681e-01} &  &  &  &  \\
& 128 & \textbf{4.5645e-02} & \textbf{4.7147e-02} & \textbf{3.2007e-01} & \textbf{3.2885e-01} & 1.221 & 1.663 & 0.540 & 0.505 \\
\hline
SharpN2KUNK
& 64  & \textbf{1.0637e-01} & 1.4963e-01 & 4.6532e-01 & 4.6684e-01 &  &  &  &  \\
& 128 & \textbf{4.5645e-02} & 4.7148e-02 & 3.2019e-01 & 3.2892e-01 & 1.221 & 1.666 & 0.539 & 0.505 \\
\hline\hline

\multicolumn{10}{l}{\textbf{DLMP group: $VN(K)$}}\\
\hline
Kuzmin
& 64  & 1.1331e-01 & 1.5499e-01 & 4.7375e-01 & 4.7639e-01 &  &  &  &  \\
& 128 & 4.8305e-02 & 4.8108e-02 & 3.2363e-01 & 3.3272e-01 & 1.230 & 1.688 & 0.550 & 0.518 \\
\hline
VNK
& 64  & \textbf{1.0660e-01} & \textbf{1.4969e-01} & \textbf{4.6565e-01} & 4.6722e-01 &  &  &  &  \\
& 128 & \textbf{4.5720e-02} & \textbf{4.7169e-02} & \textbf{3.2035e-01} & \textbf{3.2913e-01} & 1.221 & 1.666 & 0.540 & 0.505 \\
\hline
SharpVNK
& 64  & \textbf{1.0660e-01} & 1.4974e-01 & 4.6567e-01 & \textbf{4.6720e-01} &  &  &  &  \\
& 128 & \textbf{4.5720e-02} & 4.7170e-02 & 3.2037e-01 & \textbf{3.2913e-01} & 1.221 & 1.667 & 0.540 & 0.505 \\
\hline
\end{tabular}
\end{table*}
\newpage

\subsection{Application 1b: FV2 Burgers equation}
\label{sec:fv2_burgers}
In these experiments we take $\b v =(1,1)$ and $f(u)=u^2/2$ and solve the 2D Burgers equation
\begin{align}
u_t + (u^2/2)_x + (u^2/2)_y = 0,\quad u(x,0)=u_0,
\end{align}
with the following mesh discretisation $n_x,n_t,T_{max}=128,1028,0.5$, subject to the discontinuous initial condition 
\begin{align}
u_0 =
\begin{cases}
1.0 & \text{if } 0.1 < x < 0.6 \text{ and } 0.1 < y < 0.6, \\
0 & \text{otherwise}.
\end{cases}
\end{align}
In this specific instance of a Burgers flux, one can solve the Riemann problem exactly at the boundary of every cell (Godunov's method with a Godunov Flux), or one can approximately solve the Riemann problem with a Local-Lax-Friedrich (Rusanov) flux. The LLF and the Godunov fluxes are both monotone and defined below.

\begin{definition}[LLF-Burgers-equation-Flux]
The Local Lax Friedrich flux 
is 
\begin{align}
f_{i+1/2} &= \frac{1}{2}\left(f(u^L_{i+1/2}) + f(u^R_{i+1/2})\right) - \frac{1}{2}\alpha_{i+1/2}(u^R_{i+1/2} - u^L_{i+1/2}),
\end{align}
where
\begin{align}
 f(u) = \frac{1}{2} u^2,\quad
    \alpha_{i+1/2} = \max_{u\in \lbrace u^L_{i+1/2},u^R_{i+1/2}\rbrace} \left| \frac{\partial}{\partial u}f(u)\right| = \max \left\lbrace \left|u^L_{i+1/2}\right|, \left|u^R_{i+1/2} \right| \right\rbrace, \label{eq:alpha wavespeed}
\end{align}
\end{definition}

\begin{definition}[Godunov-Burgers-Flux]
Godunov's method solves the discontinuous initial value problems at the boundaries of cells. Let $u^L_{i+1/2},u^R_{i+1/2}$ denote the values to the left and right of the discontinuity ($u^R_{i},u^L_{i+1}$ respectively in the previous cell notation).  Burgers equation has a flux function $f(u) = u^2/2$, and numerical flux \begin{align}
f_{i+1/2}(u^L,u^R) =
\begin{cases}
f(u^L), & u^L > u^R \ \text{and}\ \dfrac{u^L+u^R}{2} \ge 0,\\[6pt]
f(u^R), & u^L > u^R \ \text{and}\ \dfrac{u^L+u^R}{2} \le 0,\\[6pt]
f(u^L), & u^L \le u^R \ \text{and}\ u^L \ge 0,\\[3pt]
f(u^R), & u^L \le u^R \ \text{and}\ u^R \le 0,\\[3pt]
0,      & u^L < 0 < u^R.
\end{cases}
\end{align}
\end{definition}

\begin{figure}[H]
\centering
\begin{subfigure}[t]{0.3\textwidth}
\includegraphics[width=.95\textwidth]{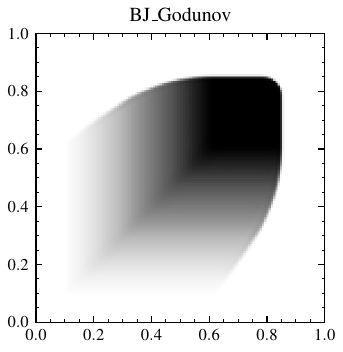}
\end{subfigure}
\begin{subfigure}[t]{0.3\textwidth}
\includegraphics[width=.95\textwidth]{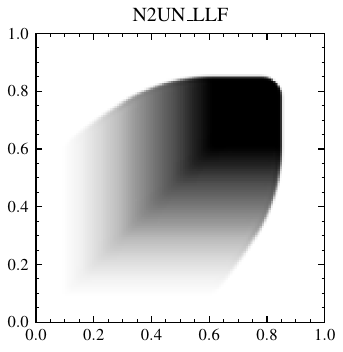}
\end{subfigure}
\begin{subfigure}[t]{0.3\textwidth}
\includegraphics[width=.95\textwidth]{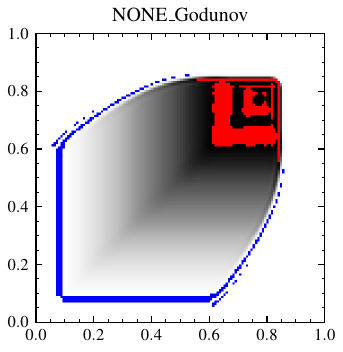}
\end{subfigure}\\
\begin{subfigure}[t]{0.3\textwidth}
\includegraphics[width=.95\textwidth]{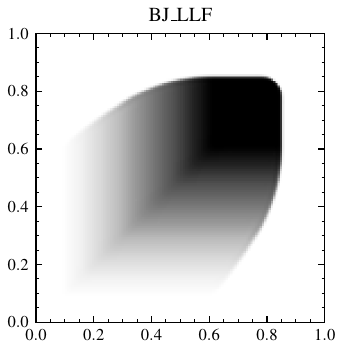}
\end{subfigure}
\begin{subfigure}[t]{0.3\textwidth}
\includegraphics[width=.95\textwidth]{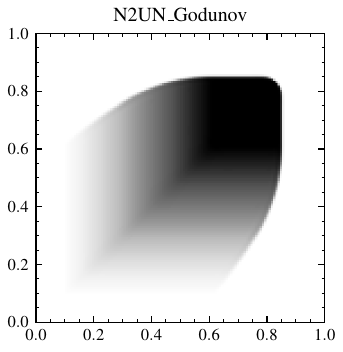}
\end{subfigure}
\begin{subfigure}[t]{0.3\textwidth}
\includegraphics[width=.95\textwidth]{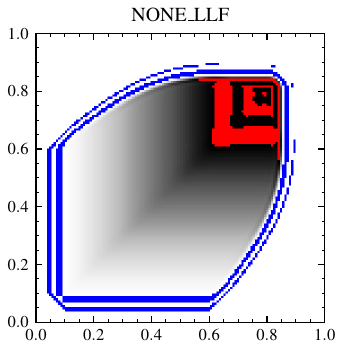}
\end{subfigure}
\caption{Here we plot the final timestep solution of Burgers equation, undershoots $u<2e-15$ are masked blue, overshoots $u>1+2e-15$ are masked red on a perceptually uniform grey colour scale between [0,1].  Solved with Godunov (row 1) and Local Lax Friedrich flux (row 2), and using Barth-Jespersen's slope limiter (column 1), $N^2(K)\cup N(K)$-limiter (column 2), and with no limiting in (row 3).}\label{fig: Burgers Final timestep Solution}
\end{figure}

We plot the solution at the final timestep in \cref{fig: Burgers Final timestep Solution} for the SSP22 time-stepping with the numerical scheme FV2, with three different limiters and both Godunov/LLF Riemann solvers. We observe that both Barth-Jespersen's slope limiter (column 1) and $N^2(K)\cup N(K)$-limiter (column 2) remain bounded, but the unlimited scheme produces both undershoots and overshoots in (column 3). 

In \cref{table:L2 relative error Burgers equation} we compare the solutions to a high-resolution reference set in the relative $l^2$ error. From here we observe that using no limiting is more accurate in the relative L2 norm than when using the BJ-limiter or $N^2(K)\cup N(K)$-limiter at the specific resolutions $32,128,256$, despite the monotonicity violations in \cref{fig: Burgers Final timestep Solution}. In \cref{table:L2 relative error Burgers equation} the unlimited scheme at resolution $64$ is less accurate than the slope limited schemes. Additionally, the Local Lax Friedrich (Russanov) flux is less accurate than the Godunov solver.

\begin{table}[H]
\centering
\begin{tabular}{||c|c|c|c|c|c||}
\hline
Limiter & Riemann Solver & 32&64&128&256 \\
\hline
NONE & Godunov & 0.342797 & 0.248308 & 0.0729384 & 0.044121 \\
BJ & Godunov & 0.343336 & 0.246393 & 0.0859516 & 0.0578801 \\
N2UN & Godunov & 0.343292 & 0.246425 & 0.085589 & 0.0575774 \\
\hline
NONE & LLF & 0.341244 & 0.246801 & 0.0745044 & 0.0458193 \\
BJ & LLF & 0.343864 & 0.246236 & 0.0914636 & 0.061922 \\
N2UN & LLF & 0.343747 & 0.246279 & 0.0910518 &0.0615696\\ 
\hline
\end{tabular}
\caption{L2 relative error as compared with a high-resolution ($2048\times 2048$) simulation.}
\label{table:L2 relative error Burgers equation}
\end{table}

We also observe that the $N^2(K)\cup N(K)$ limiter is marginally less accurate than the BJ limiter only at the $64\times 64$ resolution, where at this specific resolution the higher order flux without limiting (NONE) leads to a less accurate solution and the $N^2(K)\cup N(K)$ limiter uses more of this (less accurate) flux as compared with the Barth-Jespersen Limiter. 

Conversely, in \cref{table:L2 relative error Burgers equation} the $N^2(K)\cup N(K)$ limiter is more accurate than the Barth-Jespersen limiter when using either a Godunov or Local-Lax-Friedrichs (Russanov) solver at the resolutions $32^2,128^2,256^2$. This is because the higher order flux leads to a more accurate solution, and the $N^2(K)\cup N(K)$ limiter uses more of the higher order flux than the BJ limiter. In summary, less restrictive correction factors tend to lead to more accurate solutions when the higher order unlimited flux leads to a more accurate solution.

We then repeat the same experiment (initial conditions and equation) with additional limiters, and extend to twice the runtime $T_{max}=1$, fix the CFL near 0.125, and compare to a $1024^2$ $N^2(K)\cup N(K)$ high resolution reference. Results are found in \cref{tab:dlmp-summary-stacked-64-128_burgers}, which reports more comprehensive accuracy and boundedness metrics at $n_x=64$ and $n_x=128$, for more limiters grouped by limiter class for the Godunov flux.

From this extended experiment, in \cref{tab:dlmp-summary-stacked-64-128_burgers} we observe that, from within the $N(K)\cup K$ class, NKUK-Sharp is least diffusive (smallest errors, largest maxima), whereas NKUK and GBJ-$K$ exhibit stronger peak damping and larger errors. In the $N^{2}(K)\cup N(K)$ class, N2KUNK (and its sharp variant) and BJ are closely matched but ordered consistent with the fact that the N2KUNK is provably less restrictive than BJ (i.e. applies weaker limiting for the same cell-mean bounds), however the present Burgers test case does not strongly separate their relative L2 error metrics or peak preservation ability. The vertex-based group behaves similarly: VNK and Sharp-VNK are essentially indistinguishable and slightly outperform Kuzmin in both error and peak retention, in line with the analytical result that VNK has less severe correction factors than Kuzmin. The unlimited scheme attains the smallest $L^2$ errors but exhibits clear undershoots and overshoots.

\begin{table}[t]
\centering
\caption{Summary statistics, grouped by limiter class. For each method, the $n_x=64$ row is followed by the $n_x=128$ row. Within each limiter class, the smallest \emph{rL2} (per resolution) is in bold, and the largest \emph{max} (per resolution) is in bold (ties are both in bold).}
\label{tab:dlmp-summary-stacked-64-128_burgers}
\small
\setlength{\tabcolsep}{5pt}
\begin{tabular}{llcccc}
\toprule
Class & Method & $n_x$ & rL2 & min & max \\
\midrule

\multirow{6}{*}{$N(K)\cup K$ DLMP}
& NKUK Sharp & 64  & \textbf{1.8007e-01} & -3.3747e-36 & \textbf{9.2501e-01} \\
&           & 128 & \textbf{1.0087e-01} & -3.2439e-36 & \textbf{9.5207e-01} \\
& NKUK       & 64  & 2.7114e-01 & -5.6385e-36 & 8.6068e-01 \\
&           & 128 & 1.6593e-01 & -6.5467e-36 & 9.0164e-01 \\
& GBJ\_K     & 64  & 3.1893e-01 &  0.0000e+00 & 8.1386e-01 \\
&           & 128 & 2.0434e-01 &  0.0000e+00 & 8.6951e-01 \\
\midrule

\multirow{6}{*}{$N^{2}(K)\cup N(K)$ DLMP}
& N2KUNK Sharp & 64  & \textbf{1.7752e-01} & -2.4497e-36 & \textbf{9.2769e-01} \\
&             & 128 & \textbf{9.9662e-02} & -2.5182e-36 & \textbf{9.5283e-01} \\
& N2KUNK       & 64  & 1.7800e-01 & -2.7967e-36 & 9.2768e-01 \\
&             & 128 & 9.9898e-02 & -2.8544e-36 & \textbf{9.5283e-01} \\
& BJ           & 64  & 1.7918e-01 & -3.3446e-36 & 9.2723e-01 \\
&             & 128 & 1.0064e-01 & -3.8740e-36 & 9.5258e-01 \\
\midrule

\multirow{6}{*}{$VN(K)$ DLMP}
& VNK Sharp & 64  & \textbf{1.7817e-01} & -2.5959e-36 & \textbf{9.2758e-01} \\
&          & 128 & \textbf{1.0007e-01} & -2.7479e-36 & \textbf{9.5277e-01} \\
& VNK      & 64  & \textbf{1.7817e-01} & -2.5959e-36 & \textbf{9.2758e-01} \\
&          & 128 & \textbf{1.0007e-01} & -2.7479e-36 & \textbf{9.5277e-01} \\
& Kuzmin   & 64  & 1.8152e-01 &  0.0000e+00 & 9.2421e-01 \\
&          & 128 & 1.0136e-01 &  0.0000e+00 & 9.5213e-01 \\
\midrule

\multirow{2}{*}{No DLMP (unlimited)}
& No Limiter & 64  & \textbf{1.3012e-01} & -2.1448e-02 & \textbf{1.1058e+00} \\
&           & 128 & \textbf{4.3011e-02} & -1.9110e-02 & \textbf{1.1107e+00} \\
\bottomrule
\end{tabular}
\end{table}

\subsection{Application 1c: FV2 Incompressible Euler}\label{sec:fv2_euler}
\begin{figure}[H]
    \centering
\includegraphics[width=1.0\linewidth]{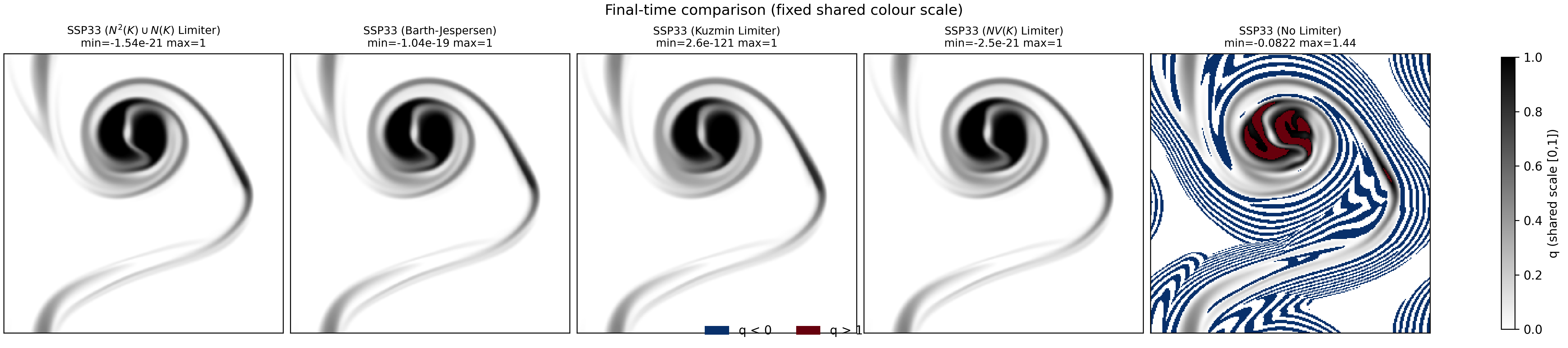}
    \caption{Final time solution of Euler's equation (global overshoots red, global undershoots blue). Limiters remain bounded and are visually similar.}
    \label{fig:Euler}
\end{figure}

For an additional demonstration, we solve 2D incompressible Euler's equation 
\begin{align}
q_t + \operatorname{div}(\b u q)=0,\quad \b u = -\nabla^{\perp}\psi, \quad \psi = -\Delta^{-1}q. 
\end{align}
An FFT-based spectral solver solves for the stream function, and a standard C-grid formulation averages cell centred stream function values to vertex values, and second order difference convert stream function into edge defined velocities. FV2 with SSP33 \cite{shu1988efficient} is the solver. Initial conditions are LeVeque in the vorticity variable. Final time vorticity solutions are plotted in \cref{fig:Euler} at $t_{max}=40$ with no overshoots when limiting, indicating that at second order, new limiters are competitive with existing ones. The next section shows new limiting techniques can be applied to higher order schemes as well.

\section{Application 2: Higher order limiting}\label{sec: Application2: chapter 2}

\subsection{FV4: Fourth order finite volume}\label{sec: Application2: fv4 chapter 2}
We define a fourth order finite volume method, which is directly applicable for a 2D orthogonal grid and bears some similarity to the MCORE \cite{ullrich2012mcore} finite volume dynamical core, but does not use a convolution and deconvolution strategy for the fluxes. Instead, the scheme uses direct evaluations at Gauss points from the high order subcell representation. 

It can be defined by a sequence of compositions,
\begin{align}
	\bar{u}^{n+1} &= (\mathcal{E} \circ\mathcal{R}\circ \mathcal{Q} \circ \mathcal{G} \circ \mathcal{P} )\circ \bar{u},
\end{align}
in pseudocode format, as follows.
\begin{enumerate}
	\item We use the following fourth order projection map $\mathcal{P}_4: \bar{u}_{i,j} \mapsto u_{i,j} + O(\Delta x^4 + \Delta y^4)$ to approximate point values from cell mean values. It is consistent with respect to constants.
	\begin{align}
		u_{i,j} &= \bar{u}_{i,j}  - \frac{1}{24}[ \bar{u}_{i+1,j} - 2 \bar{u}_{i,j} + \bar{u}_{i-1,j}]
		- \frac{1}{24}[ \bar{u}_{i,j+1} - 2 \bar{u}_{i,j} + \bar{u}_{i,j-1}], \quad \forall (i,j)
	\end{align}
	
	\item We use the gradient map 
	\begin{align}
		\mathcal{G}_3: u \mapsto u, u_x, u_{y}, u_{xx}, u_{xy}, u_{yy},  u_{xxx},u_{xxy},u_{xyy}, u_{yyy},
	\end{align}
	defined by the fourth order centred finite difference weights
	\begin{align}
		w1 &= 1/12( [-1,8,0,-8,1])\\
		w2 &= 1/12 ([-1,16,-30,16,-1])\\
		w3 &= 1/8 ([-1,8,-13,0,13,-8,1]),
	\end{align}
	to construct $u_{x},u_{xx},u_{xxx}$ and $u_{y},u_{yy},u_{yyy}$ from the newly computed point values. We use these newly computed values, and the finite difference stencil $1/12( [-1,8,0,-8,1])$, to compute all the missing cross term derivatives $u_{xy},u_{xxy},u_{xyy}$ within each subcell representation.
	\item 
    %, 
    We evaluate a quadrature map $\mathcal{Q}_4: (x_q,y_q) \mapsto p_{K}(x_q,y_q)$, $\forall (x_q,y_q) \in K$ to compute a set of fourth order accurate quadrature point evaluations at $(x_q,y_q)$ for all cells. We do so by evaluating the following formula of the subcell representation,
%	\begin{align}
%		P_{i,j}(x,y) &= u_{i,j} + (x-x_i)u_x + (y-y_j ) u_y + \\
%		&\frac{1}{2}[  (x-x_i)^2 u_{xx} + 2 (x-x_i)(y-y_i) u_{xy} + (y-y_i)^2 u_{yy}  ] \\
%		&\frac{1}{3!}[  (x-x_i)^3 u_{xxx} + 3 (x-x_i)^2(y-y_i) u_{xxy} \\
%		&+ 3 (x-x_i)(y-y_i)^2 u_{xyy} + (y-y_i)^3 u_{yyy}  ]. 
%	\end{align}
\begin{align}
p_{i,j}(x,y) &= \bar{u}_{i,j} + (x-x_i)u_x + (y-y_j ) u_y + \\
&\frac{1}{2} \left[  \left( (x-x_i)^2 -\frac{\Delta x^2}{12}\right) u_{xx} + 2 (x-x_i)(y-y_j) u_{xy} + \left( (y-y_j)^2  -\frac{\Delta y^2}{12}\right) u_{yy}  \right] \\
&+\frac{1}{3!} \left[  (x-x_i)^3 u_{xxx} + 3 (x-x_i)^2(y-y_j) u_{xxy} + 3 (x-x_i)(y-y_j)^2 u_{xyy} + (y-y_j)^3 u_{yyy}  \right]. 
\end{align}
	\item We define the map $\mathcal{R}_{4}$ which resolve Riemann problems and reconstructs a flux as follows. We resolve the local quadrature defined Riemann problems 
	\begin{align}
		F(x_{i+1/2},y_q,v(x_q,y_q)) = p_{K}(x_{i+1/2},y_q)v(x_{i+1/2},y_q)^+ + p_{L}(x_{i+1/2},y_q)v(x_{i+1/2},y_q)^- .
	\end{align}
	The flux is computed by a fourth order Gauss quadrature, for example the right edge is computed using
	\begin{align}
		F_{i,i+1} &= \sum_{q_k \in \sigma_{i,i+1} } w_{q_{k}} \left[F\left(x_{i+1/2},y_{q_{k}},v(x_{i+1/2},y_{q_{k}})\right)\right], \\
		\text{where} \quad \lbrace w_{q_1}, w_{q_2}\rbrace &= \lbrace 1/2,1/2\rbrace, \quad \lbrace y_{q_1},y_{q_2}\rbrace = \left\lbrace y_{j+1/2}-\frac{\Delta y}{2\sqrt{3}}, y_{j+1/2}+\frac{\Delta y}{2\sqrt{3}}\right\rbrace.
	\end{align}
	\item  The final stage, $\mathcal{E}$, involves the normal cell mean evolution procedure, where the fluxes on each face are used to update the solution.
	\begin{align}
\bar{u}^{n+1}_{i,j} = \bar{u}^{n}_{i,j} - \frac{\Delta t }{|K|}\sum_{L\in N(K)} |\sigma_{KL}| F_{K,L}.
	\end{align}
\end{enumerate}
We employ the upwind/downwind variant of the standard RK4 scheme in \cite{shu1988efficient}, with radius of monotonicity $2/3$. Technically, this changes the scheme into an additive Runge-Kutta method. We employ limiting in each forward Euler stage of the SSP(A)RK method. Note that SSP schemes only preserve convex functionals, discrete local maximum principles are not preserved in the same way, one ensures a DLMP holds for the substages of the RK scheme owing to the convex combination \cite{woodfield2024new}.

\subsection{$N^{2}(K)\cup N(K)$ limiter for FV4}
Based on \cref{thm: slope limiters}, we wish to employ the $N^2(K)\cup N(K)$ limiter (and other similar limiters) introduced in \cref{sec:instanciations_of_limiters} to the fourth order finite volume method FV4. We first remark on some non-trivial facts about this specific finite volume construction and how the limiter interacts with the scheme.
\begin{enumerate}
\item The eight flux contributing quadrature points $(\b x_q)$ for cell $(i,j)$ are located at the positions
\begin{align}
\left[x_{i} \pm \frac{\Delta x}{2}, y_{j}\pm \frac{\Delta y}{2\sqrt{3}}\right], \quad \left[x_{i} \pm \frac{\Delta x}{2\sqrt{3}}, y_{j}\pm \frac{\Delta y}{2}\right].
\end{align}
These points are limited by face defined maximum principles,
two points per face share the same maximum principle.
\item The subcell representation is cell mean preserving:
\begin{align}
\frac{1}{h^2}\int_{-h/2}^{h/2} \int_{-h/2}^{h/2} p_{i,j}(x,y)dx dy =  \bar{u}_{i,j}.
\end{align}
\item There exists a convex Zhang-acceptable decomposition of the cell average onto flux contributing quadrature points and the cell midpoint,
\begin{align}
\bar{u}_{i,j} = \frac{1}{2}p_{ij}(x_i,y_j) + \frac{1}{16}\sum_{q\in K^{fc}} p_{ij}(\b x_q). \label{eq:zhangFV4decomposition1}
\end{align}
The cell midpoint $(x_i,y_j)$ is not flux contributing and must be limited separately to the flux contributing quadrature points.
\item 
The local Riemann problem at the upper quadrature point on the right face takes the form
 \begin{align}
	Rei_{i+1/2,j+\frac{\Delta y}{2\sqrt{3}}}= \frac{1}{16}[ u^{R1}_{i,j} - \frac{8\Delta t}{\Delta x} F(u^{R1}_{i,j},u^{L1}_{i+1,j}, \b v_{i+1/2,j+\frac{\Delta y}{2\sqrt{3}}} \cdot n_{i,i+1})].
\end{align}
\item The Courant number limit is 1/8 for compressible flow 
 \begin{align}
	C_{K} = \sum_{L\in N(K)} \frac{\Delta t|\sigma_{KL}| (\b v \cdot n_{KL})^+}{|K|} \leq 1/8, \quad \forall K \in \mathcal{M},
\end{align}
 and $1/4$ for incompressible flow. This can be identified by making the associations $w^{\sigma_{KL}}_{q} = 1/2$ and $\gamma^{K}_{q,L} = 1/16$.
\end{enumerate}

Given these properties of the scheme, we now provide a demonstration on how to use the $N^2(K)\cup N(K)$ limiter (and other similar variants) in \cref{def41}.
%We have stated enough about the scheme to use the $N^2(K)\cup N(K)$ limiter (and other variants similar to it). 

\begin{figure}
\hspace{0cm} 
 %\centering
\begin{subfigure}[b]{0.35\textwidth}
	\includegraphics[scale=0.28]{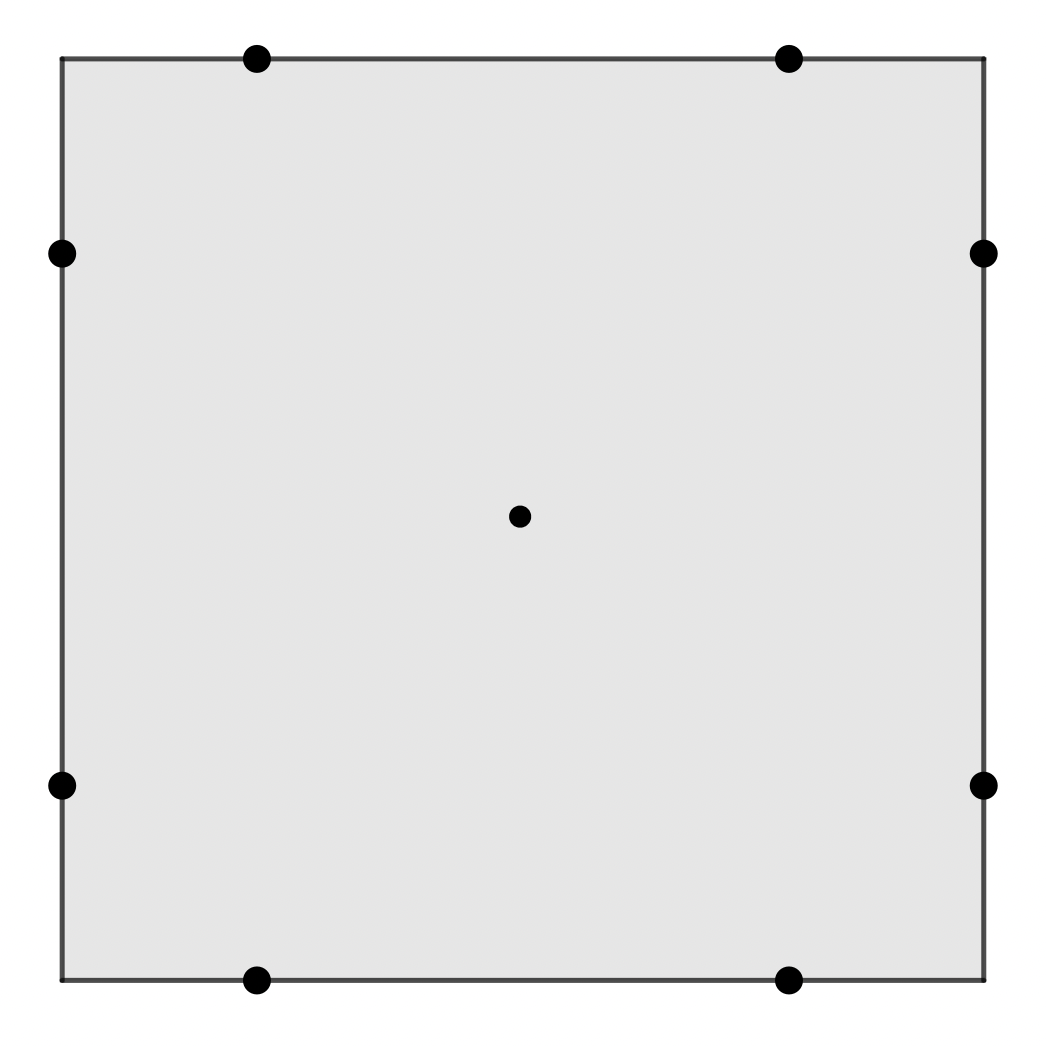}%{/Users/jmw/Desktop/Figures_diagrams/Points.png}
	\subcaption{Points used in the Zhang-acceptable cell mean decomposition of FV4 as in \Cref{eq:zhangFV4decomposition1}. There are two flux contributing quadrature points per face at Gauss nodes and one cell cell midpoint evaluation.}
	\end{subfigure}
\hspace{1cm} 
	\begin{subfigure}[b]{0.5\textwidth}
	\includegraphics[scale=0.33]{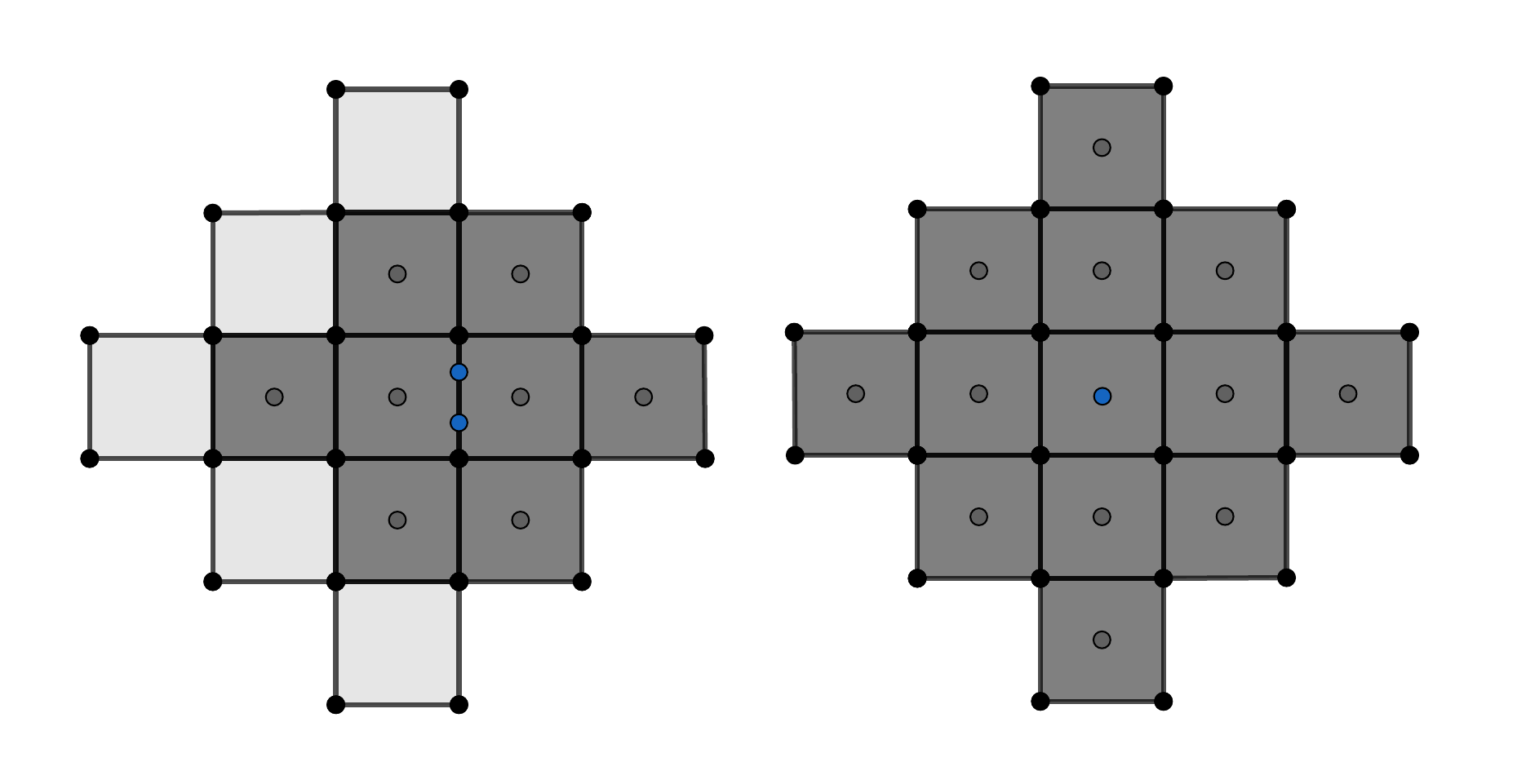}%{/Users/jmw/Desktop/Figures_diagrams/atty1.pdf}
	\subcaption{Flux contributing quadrature points at the edge $\sigma_{KL}$ are limited based on being bounded by the cell mean values $\bar{u}$ in the $N(K) \cup N(L)$ region (darker grey left diagram).  Non-flux contributing quadrature point evaluation of the midpoint $p(x_i,y_j)$ is limited to be locally bounded by the cell mean values $\bar{u}$ in the $N^2(K) \cup N(K)$ region (darker grey right diagram).}
	\end{subfigure}
\caption{Points from the FV4 cell mean decomposition, and interaction with the $N^2(K)\cup N(K)$ limiter. }
\label{fig:FV4:N2UNK limiter principle.}
\end{figure}

\begin{definition}[simplification of $N^2(K)\cup N(K)$-limiter]\label{def41}
	We point to \cref{fig:FV4:N2UNK limiter principle.} and captions in \cref{fig:FV4:N2UNK limiter principle.} for a description of this limiter with diagrams. 
	\begin{enumerate}
		\item Per face $\sigma_{KL}$, we compute and associate the local face defined maximum principle bounds,
		\begin{align}
[m_{\sigma_{KL}},M_{\sigma_{KL}}]= \left[\min_{M \in N(L) \cup N(K)} \bar{u}^n_{M},\max_{M \in N(L) \cup N(K)} \bar{u}^n_{M}\right].
		\end{align}
		This principle is associated to both quadrature points $\b x_{q} \in \sigma_{KL}$ at the face. 
		\item  Per cell $K$ we compute the desired maximum principle
		\begin{align}
			[m_{K^{nfc}} , M_{K^{nfc}}] = \left[\min_{L \in N^2(K)\cup N(K)} \bar{u}^n_{L},\max_{L \in N^2(K)\cup N(K)} \bar{u}^n_{L}\right].
		\end{align}
		This is associated to the one non-flux contributing quadrature point $\b x_{q} \in K^{nfc}$ located at the cell midpoint.
		\item 
		We then per cell compute all the Barth-Jespersen quadrature corrections factors $\alpha_{q}$, to ensure $\tilde{p}_{K}(x) = \alpha (p_{K}(x) - \bar{u}_{K}) + \bar{u}_{K}$, satisfies the conditions in \cref{thm: slope limiters}
		\begin{align}
			&\tilde{p}_{K}(\b x_{q}) \in [m_{K^{nfc}},M_{K^{nfc}}],\quad \forall q \in K^{nfc},\\
			&\tilde{p}_{K}(\b x_{q}) \in [m_{\sigma_{KL}},M_{\sigma_{KL}}],\quad \forall q \in \sigma_{KL} \cap K^{fc}, \quad \forall L\in N(K),
		\end{align}
		by choosing the smallest value
		\begin{align}
			\alpha = \min_{\forall q \in K} \alpha_{q}.
		\end{align}
	This ensures the limited internal subcell representation $\tilde{p}_{K}(x) = \alpha (p_{K}(x) - \bar{u}_{K}) + \bar{u}_{K}$ satisfies the required edge sharing quadrature maximum principles for both flux contributing quadrature points, and the cell midpoint satisfies a non-flux contributing quadrature point maximum principle (\cref{fig:FV4:N2UNK limiter principle.}).
	\end{enumerate}
\end{definition}

We employ the upwind bias SSP variant of RK4 in \cite{shu1988efficient} in the implementation. 

\begin{remark}\label{remark:zhang nonuniqueness}
There exists other Zhang-acceptable decompositions of the cell mean such as 
\begin{align}
	\bar{u}_{i,j} &=\frac{\theta}{2}p_{ij}(x_i,y_{j})+ \frac{1-\theta}{8}[p_{ij}(x_i,y_{j+1/2})+p_{ij}(x_i,y_{j-1/2})+p_{ij}(x_{i+1/2},y_{j}) +p_{ij}(x_{i-1/2},y_{j})] \\
	&+ \frac{1}{16}\sum_{q\in K^{fc}} p_{ij}(\b x_q), \quad \theta \in [0,1],
\end{align}
 such that the free parameter $\theta$ could be locally varied to minimise the Barth-Jespersen correction factors arising from the non-flux contributing quadrature principle. This could be used for increased accuracy. We take $\theta=1$. 
\end{remark}

% \section{FV4: Experiments}
\subsection{FV4: Experimental Convergence}\label{sec:FV4: Experimental Convergence}
We consider the time-reversing deformational flow generated by
\begin{align}
\psi(x,y,t) &= 0.125\sin(2\pi x)\sin(2\pi y)\cos(\pi t/T_{\max}),\
q_0(x,y) = 0.5+0.5\sin(2\pi x)\cos(2\pi y),\quad T_{\max}=1.
\end{align}
Using the analytic solution at $t=T_{\max}$ as reference, we report $L^p$ errors and refinement rates for FV4 with various limiters. Results in \cref{tab:convergence-errors-orders-dmpgroup} are grouped by enforced DMP-region, with the smallest error per norm and resolution highlighted. We verify that the new finite volume scheme achieves fourth-order convergence in $(L^1,L^2,L^\infty)$, when unlimited. The global bound-preserving limiter in \cite{zhang2010maximum} removes global maxima and minima violations, but reduces the observed orders to about $(3.2,2.4,1.5)$ in $(L^1,L^2,L^\infty)$ for this flow. The same order behaviour is observed for many of the local DMP variants considered here (GBJ-$N(K)\cup K$ \cref{method:GBJ_NKUK}, $N^2(K)\cup N(K)$ \cref{method:dull_N2KUNK}, Sharp-$N^2(K)\cup N(K)$ \cref{method:sharp_N2KUNK}, Sharp-$VN(K)$ \cref{method:sharp_VNK}, $VN(K)$ \cref{method:dull_VNK}) whilst preserving stronger local boundedness principles. However, the global bound-preserving limiter is slightly more accurate.

The results in \cref{tab:convergence-errors-orders-dmpgroup} show that within each DLMP-region, Sharp--$GN(K)$ variants give the smallest errors, while GBJ-type variants are most diffusive, consistent with correction factor inequalities. Under the $N(K)\cup K$ DLMP constraint, GBJ-K and NKUK exhibit strong order reduction, whereas Sharp-$N(K)\cup K$ retains orders near $(2.9,2.1,1.4)$. GBJ-NKUK is essentially indistinguishable from the limiter in WZL-\cite{WANG2004716} (which doesn't limit interior points arising from a decomposition of a cell mean), suggesting that limiting non-flux contributing quadrature points arising in a decomposition of a cell mean is necessary for a DLMP and positivity with negligible additional accuracy penalty.

Repeating the experiment but simplifying the flow to a to a unidirectional diagonal periodic flow, we observe 
similar convergence behaviour across limiters and observe true order (4,4,4) for the unlimited scheme, as seen in \cref{tab:convergence-errors-orders-dmpgroup_constant_flow}. We can therefore conclude, as above, that within each DLMP-region, Sharp--$GN(K)$ variants give the smallest errors consistent with correction factor inequalities. 

\subsection{FV4: Accuracy and monotonicity}\label{sec:FV4:Accuracy and monotonicity}
\begin{figure}[H]
\centering

% --- Row 1 (3 panels) ---
\begin{subfigure}[t]{0.32\textwidth}
  \centering
  \includegraphics[width=\linewidth]{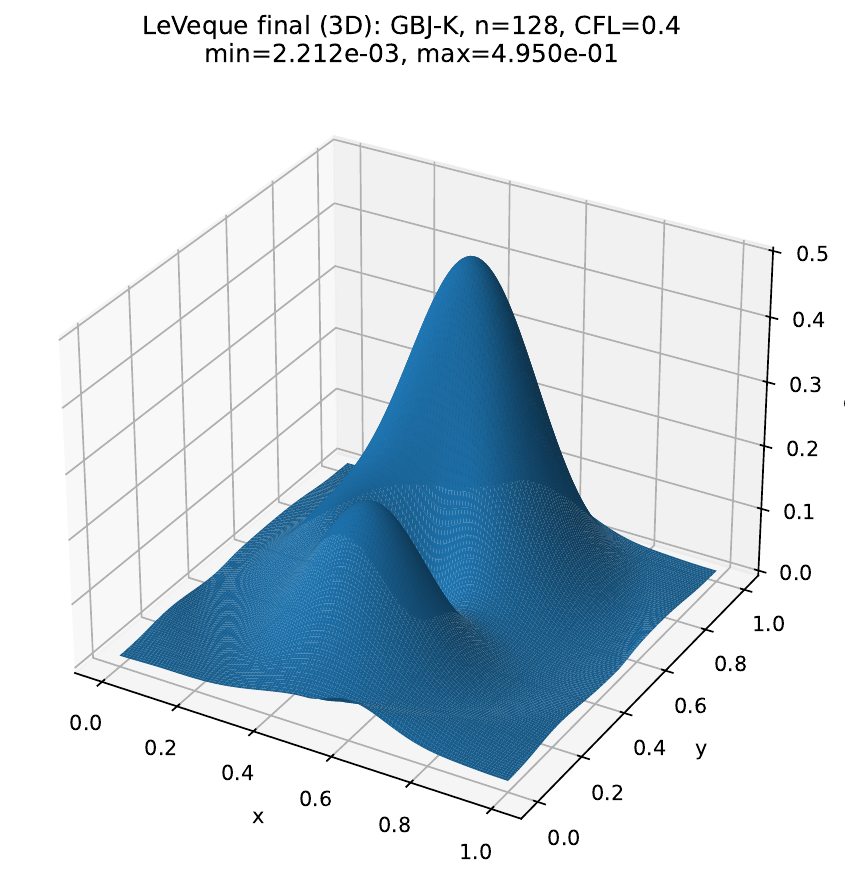}
  \caption{GBJ--$K$.}
\end{subfigure}\hfill
\begin{subfigure}[t]{0.32\textwidth}
  \centering
  \includegraphics[width=\linewidth]{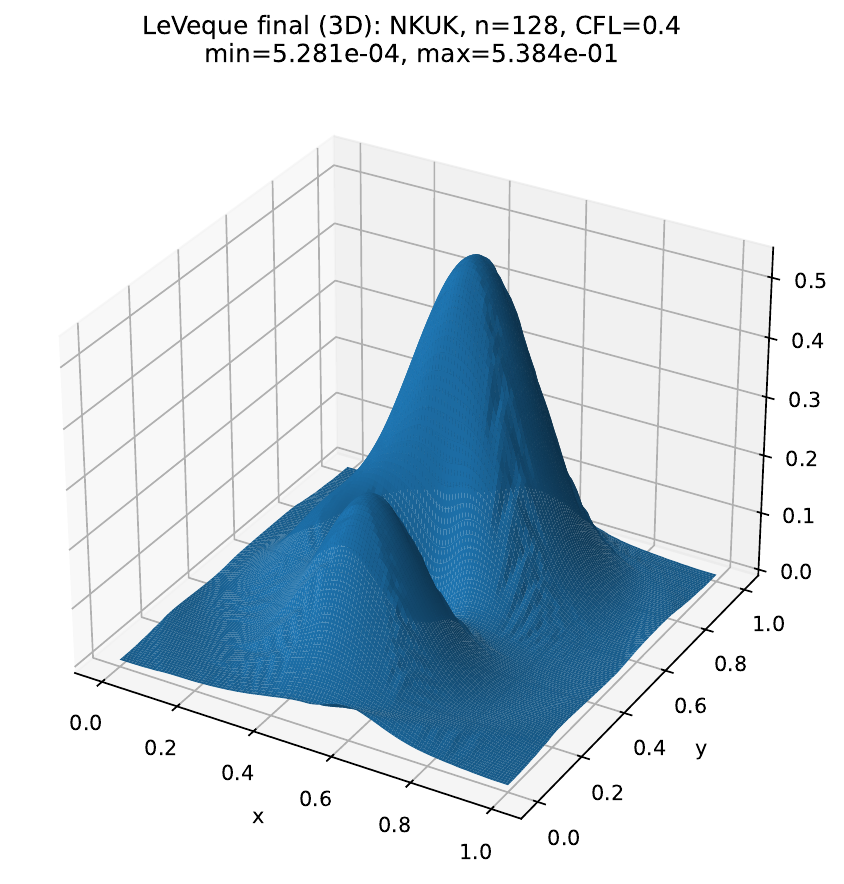}
  \caption{NKUK.}
\end{subfigure}\hfill
\begin{subfigure}[t]{0.32\textwidth}
  \centering
  \includegraphics[width=\linewidth]{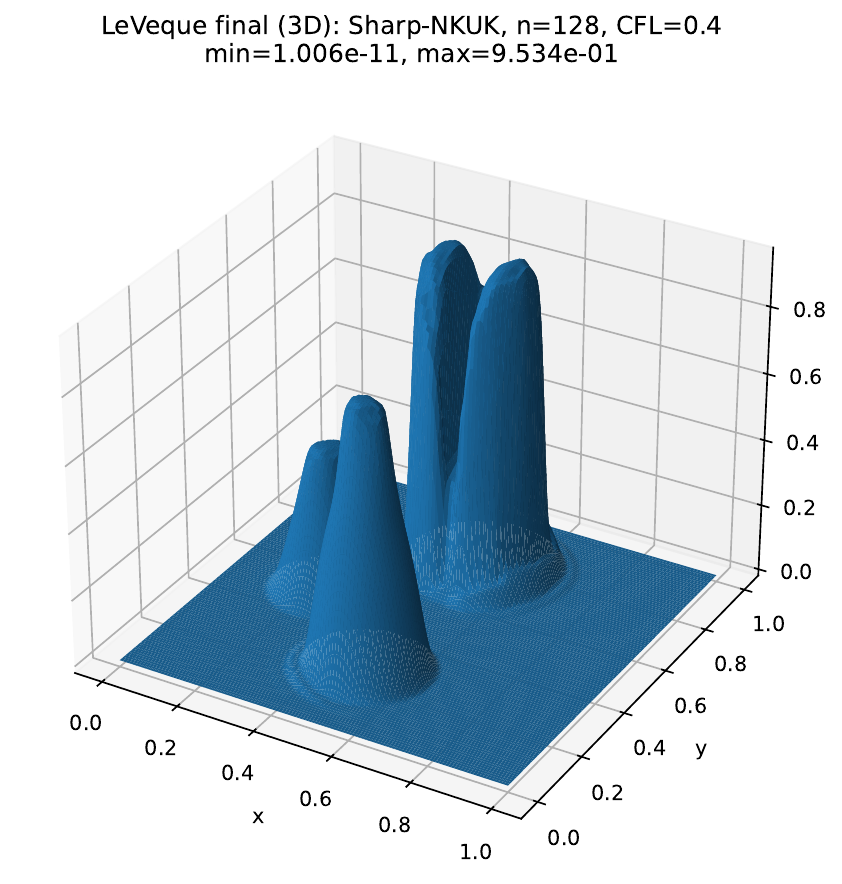}
  \caption{Sharp--NKUK.}
\end{subfigure}

\vspace{0.6em}

% --- Row 2 (3 panels) ---
\begin{subfigure}[t]{0.32\textwidth}
  \centering
  \includegraphics[width=\linewidth]{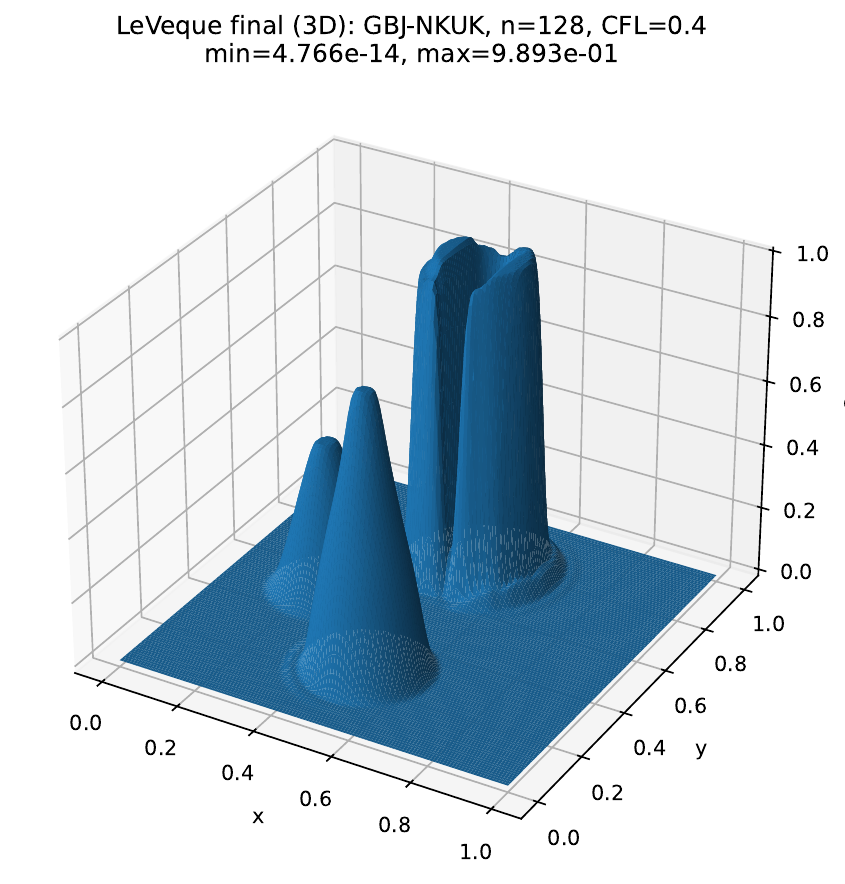}
  \caption{GBJ-NKUK}
\end{subfigure}\hfill
\begin{subfigure}[t]{0.32\textwidth}
  \centering
  \includegraphics[width=\linewidth]{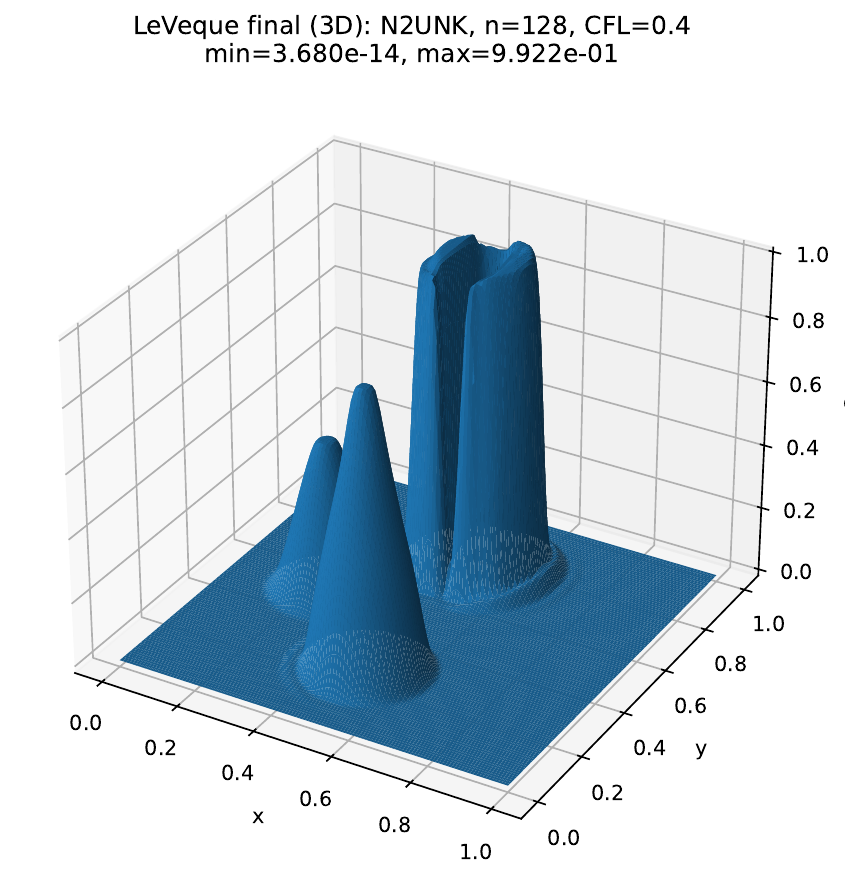}
  \caption{N2KUNK.}
\end{subfigure}\hfill
\begin{subfigure}[t]{0.32\textwidth}
  \centering
  \includegraphics[width=\linewidth]{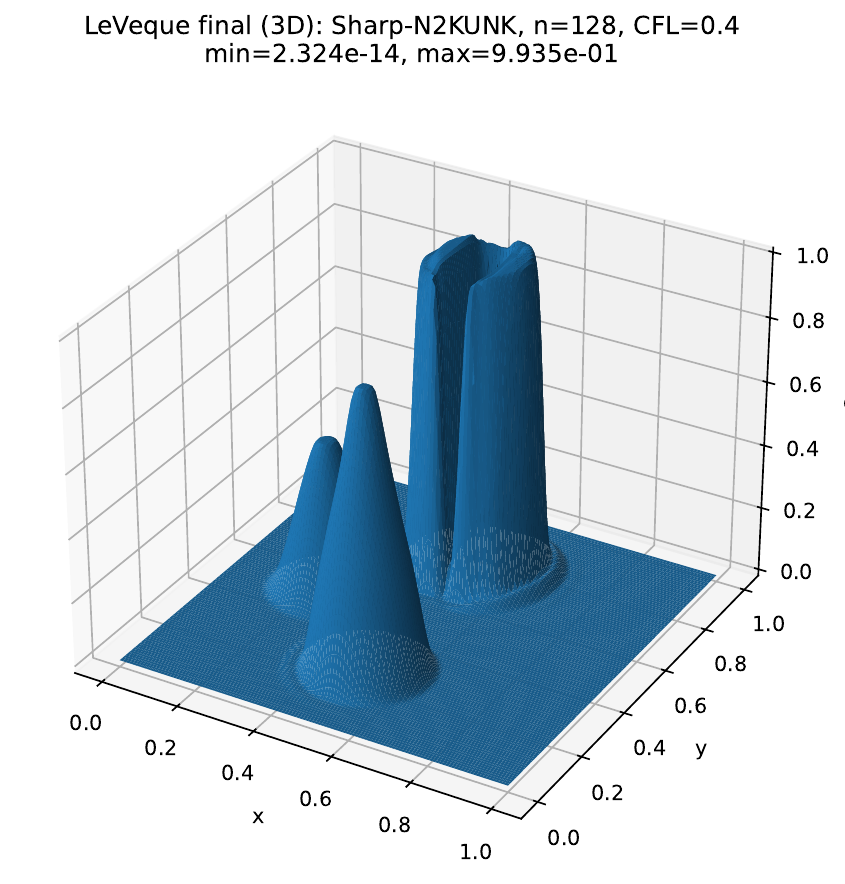}
  \caption{Sharp--N2KUNK.}
  \label{fig:leveque-3d-n2kunk-sharp-128}
\end{subfigure}

\vspace{0.6em}

% --- Row 3 (3 panels) ---
\begin{subfigure}[t]{0.32\textwidth}
  \centering
  \includegraphics[width=\linewidth]{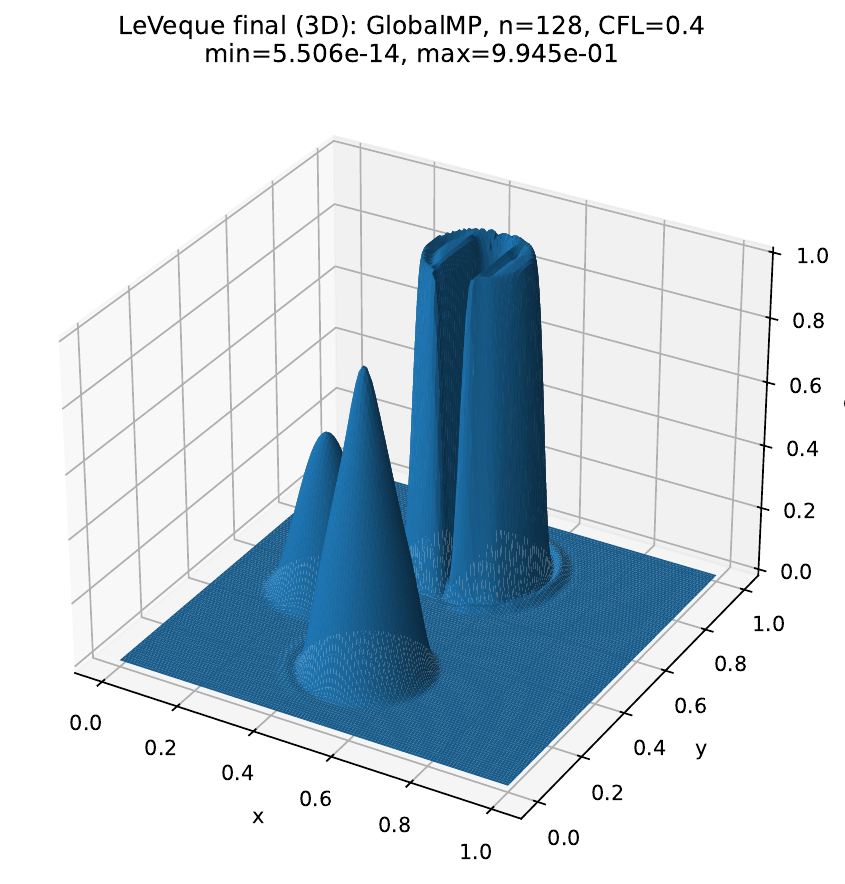}
  \caption{Global Boundedness limiter}
  \label{fig:leveque-3d-kuzmin-128}
\end{subfigure}\hfill
\begin{subfigure}[t]{0.32\textwidth}
  \centering
  \includegraphics[width=\linewidth]{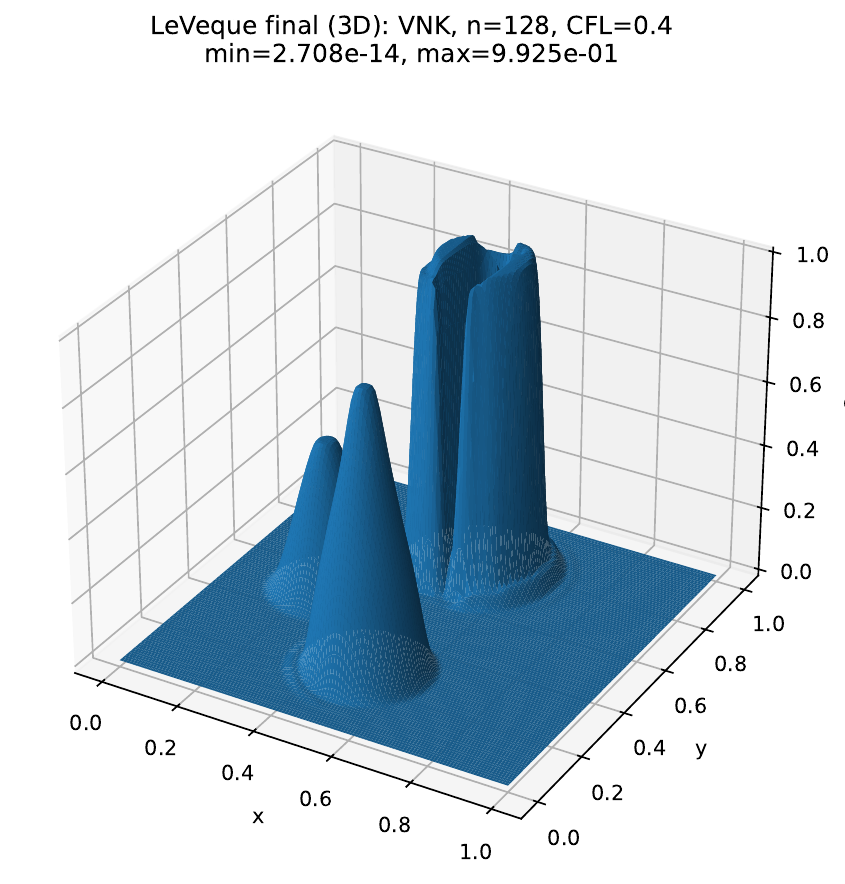}
  \caption{VNK.}
  \label{fig:leveque-3d-vnk-128}
\end{subfigure}\hfill
\begin{subfigure}[t]{0.32\textwidth}
  \centering
  \includegraphics[width=\linewidth]{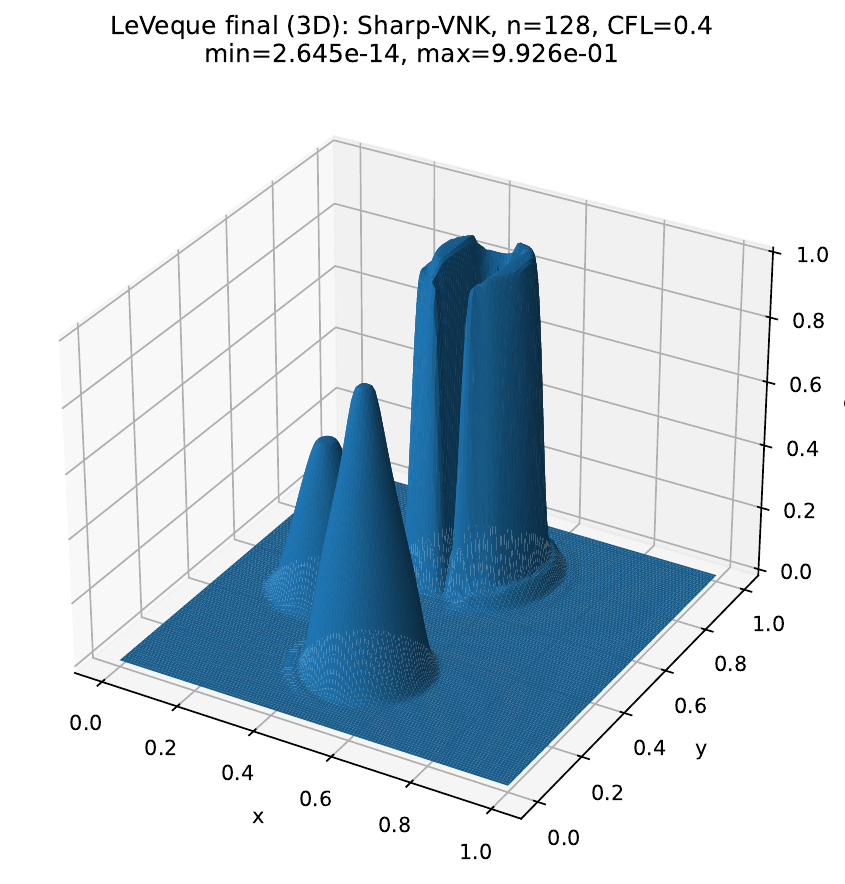}
  \caption{Sharp--VNK.}
  \label{fig:leveque-3d-vnk-sharp-128}
\end{subfigure}

\caption{LeVeque solid-body rotation at $n_x=128$. 3D surface plots at final time. Scheme: FV4-SSP(A)RK44.}
\label{fig:leveque-FV4}
\end{figure}

\begin{figure}[H]
\centering

\begin{subfigure}[t]{0.48\linewidth}
  \centering
  \includegraphics[width=\linewidth]{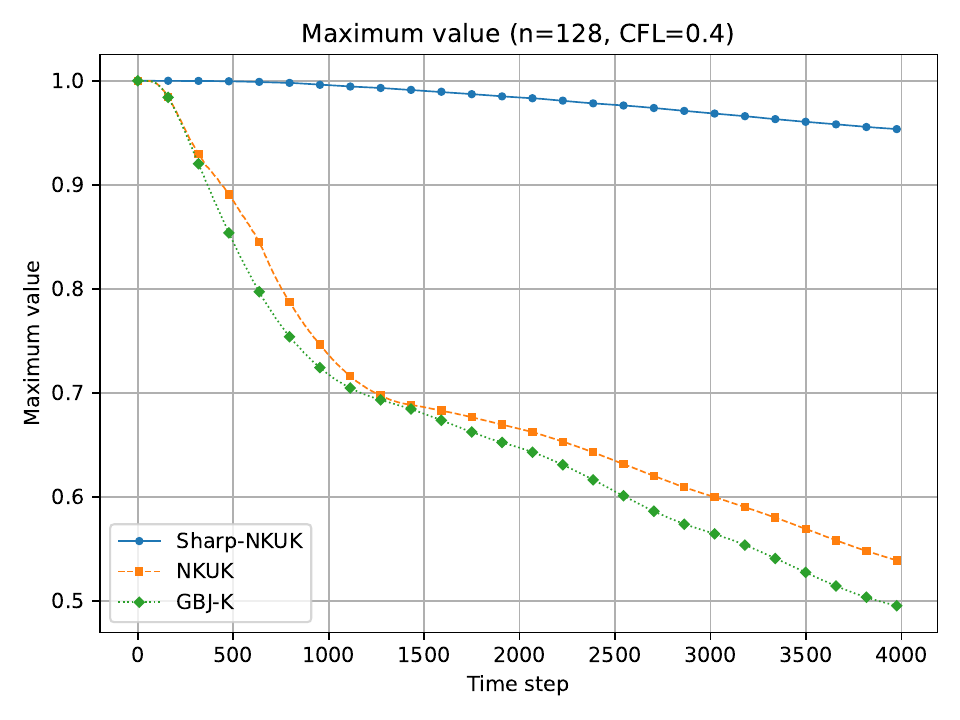}
\end{subfigure}\hfill
\begin{subfigure}[t]{0.48\linewidth}
  \centering
  \includegraphics[width=\linewidth]{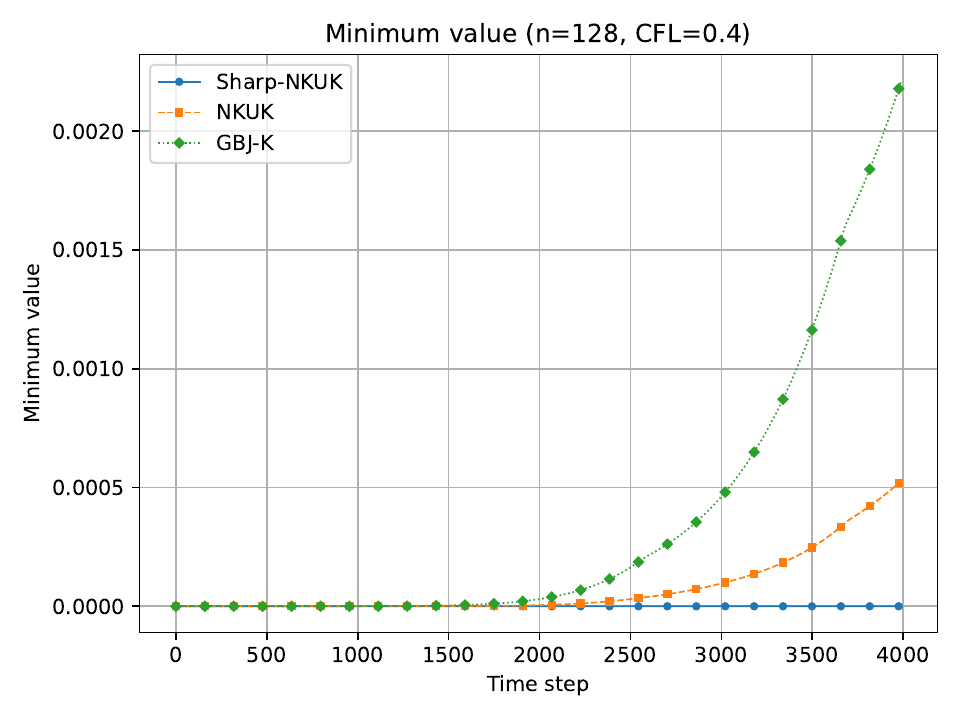}
\end{subfigure}

\vspace{0.4em}

\begin{subfigure}[t]{0.48\linewidth}
  \centering
  \includegraphics[width=\linewidth]{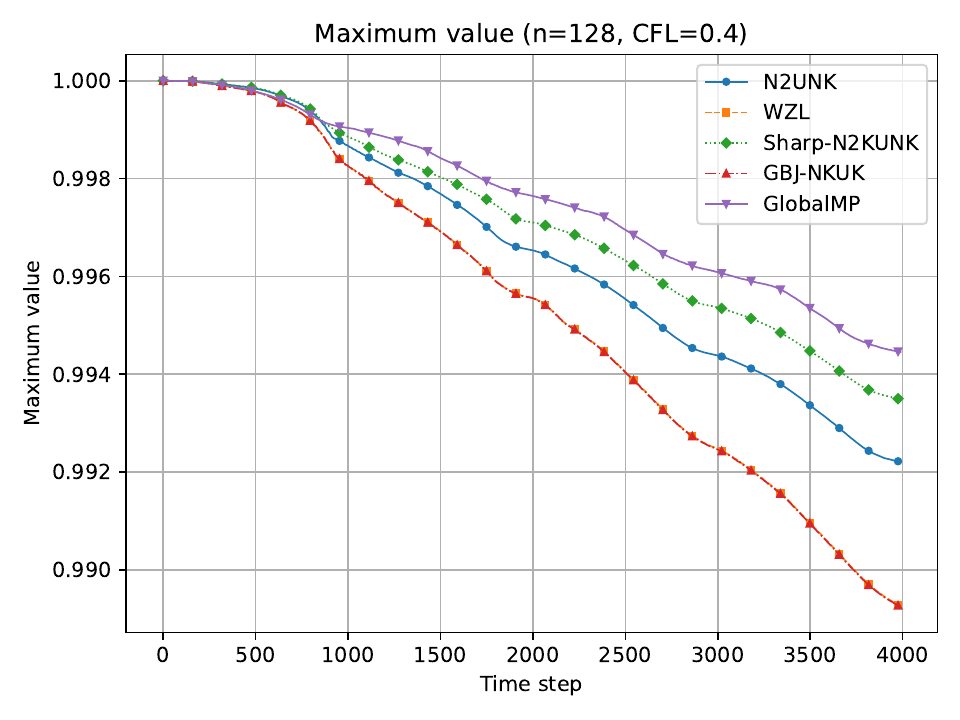}
\end{subfigure}\hfill
\begin{subfigure}[t]{0.48\linewidth}
  \centering
  \includegraphics[width=\linewidth]{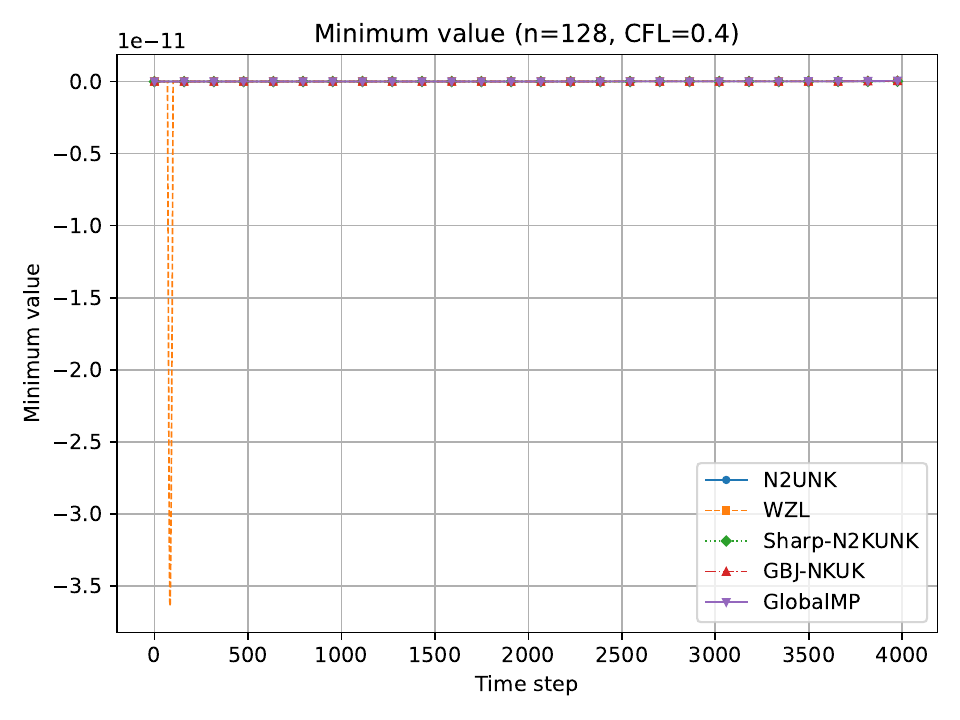}
\end{subfigure}

\vspace{0.4em}

\begin{subfigure}[t]{0.48\linewidth}
  \centering
  \includegraphics[width=\linewidth]{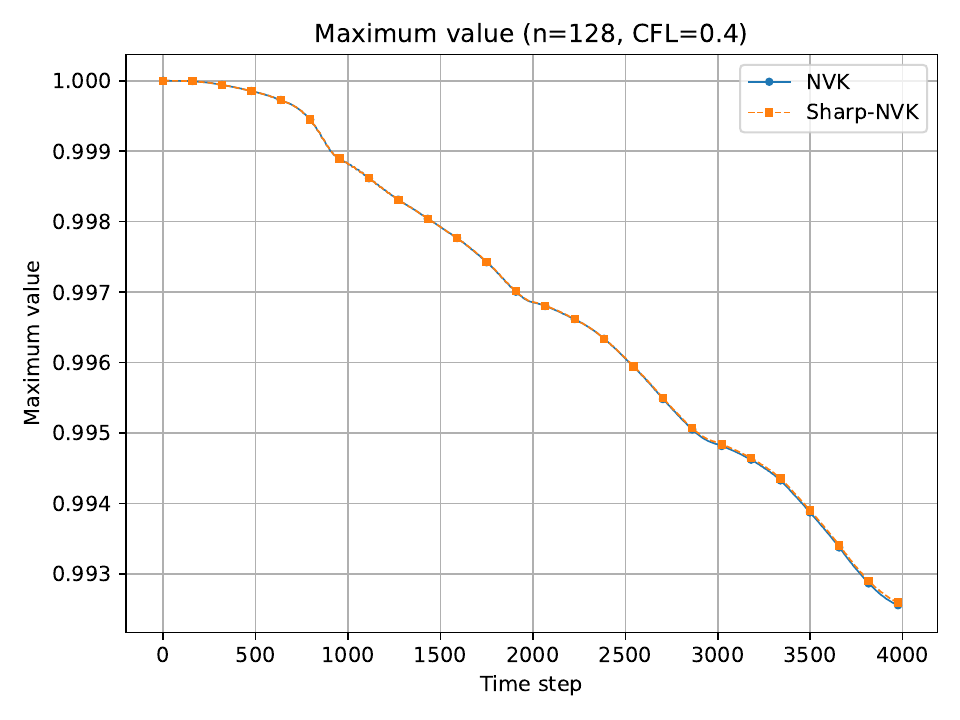}
\end{subfigure}\hfill
\begin{subfigure}[t]{0.48\linewidth}
  \centering
  \includegraphics[width=\linewidth]{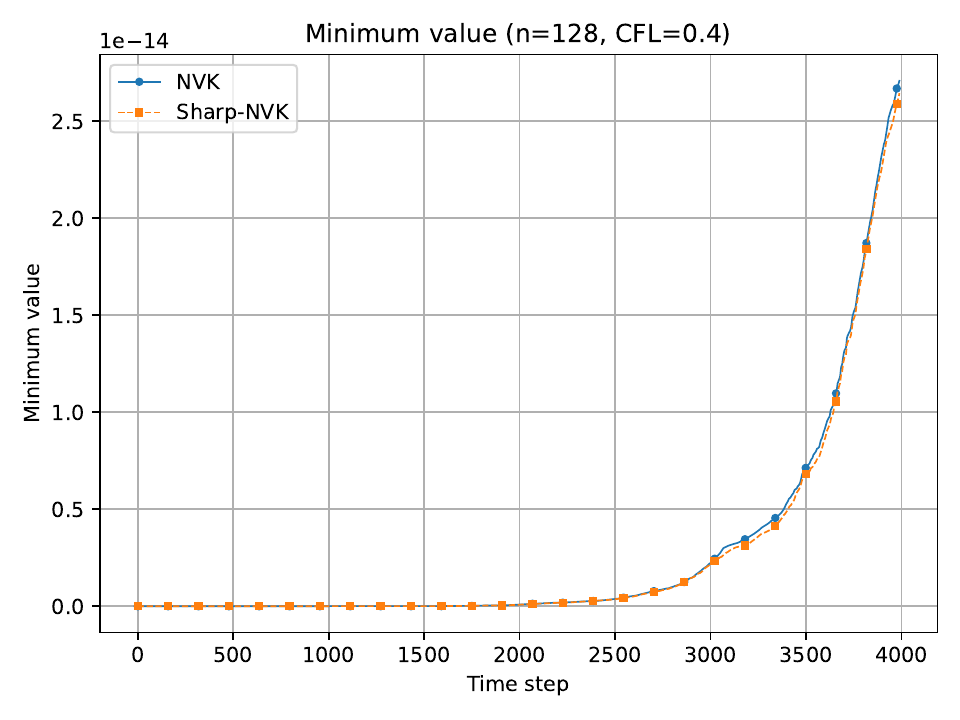}
\end{subfigure}
\caption{Here we plot spatial maxima(left column) and minima(right column) of the Zalesak initial condition, as a function of time under solid body rotation for the FV4 scheme. Sharp-limiters are better at maxima preservation indicating less diffusion. The WZL limiter is identical to GBJ-NKUK but does not limit the midpoint evaluation, and we observe a positivity violation.}
\label{fig:max_min_fv4_limited}
\end{figure}

To test the accuracy of the limiters under solid body rotation, we consider the LeVeque initial condition (\cref{test:LeVeque}) at $128^2$ and $64^2$ resolutions. We compute relative errors and convergence behaviour for the FV4 scheme under different limiters, and group results by discrete local maximum principles, seen in \cref{tab:fv4-dlmp-relL2-rates}. Final time snapshots are shown in \cref{fig:leveque-FV4}. Horizontal slices through final time shapes are concatenated and plotted together in \cref{fig:slice128_FV4_all}. Spatial maxima and minima are plotted as a function of time in \cref{fig:max_min_fv4_limited}, to assess maxima retention and boundedness.

Overall, we observe expected behaviours. Sharp limiters are slightly more accurate at all resolutions in all regions (cone, cosine and slotted cylinder) for any DLMP enforced, as is clear from \cref{tab:fv4-dlmp-relL2-rates} and in \cref{fig:slice128_FV4_all}. In \cref{fig:leveque-FV4}, we see final time solutions have expected clipping of local extrema, and results are bounded. In \cref{tab:fv4-dlmp-relL2-rates} the global maximum principle limiter \cite{zhang2010maximum} is slightly more accurate than local maximum principle limiters (one can prove correction factor inequalities consistent with this), this is expected as global maximum principle limiters do not control oscillations locally. Convergence rates of the global boundedness limiter are similar to the local boundedness limiter variants in \cref{tab:fv4-dlmp-relL2-rates}, but allow the creation of new local extrema in \cref{fig:leveque-FV4}, \cref{fig:slice128_FV4_all}. We emphasise that all schemes remain bounded to machine precision, with exception of the unlimited and WZL-\cite{WANG2004716}-limiting, under solid body rotation.

\Cref{fig:max_min_fv4_limited} contains a plot of spatial maxima and minima, of LeVeque initial conditions the FV4 scheme under solid body rotation. 
In the first row, the sharp limiter is by far the least diffusive in its corresponding DLMP limiter class.
In the second row we plot the GBJ-$N(K)\cup K$ limiter, the $N^2(K)\cup N(K)$-limiter and the Sharp-$N^2(K)\cup N(K)$-limiter, each becoming better at peak preservation, consistent with theoretical results proving less severe correction factors. We also plot in this graph (\Cref{fig:max_min_fv4_limited} row 2), the global maximum principle limiter \cite{zhang2010maximum}, it does not have any local boundedness principles and serves as a reference, this scheme is equivalent to any local neighbourhood limiter with an infinitely large stencil maximum principle. 

In the second row of \cref{fig:max_min_fv4_limited}, we compare the WLZ-\cite{WANG2004716} limiter, to the GBJ-NKUK limiter, and see that they both are identical in peak preservation, but the WLZ-\cite{WANG2004716} limiter creates a significant minimum value violation of -4e-11. This is entirely due to the fact that the  WLZ-\cite{WANG2004716} limiter doesn't limit non-flux contributing points arising in the cell mean decomposition (since this is the only difference from the GBJ-$N(K)\cup K$-limiter). 
In the third and final row, we demonstrate that the Sharp-$VN(K)$ and $VN(K)$ variants perform similarly and remain bounded.

\begin{table}[H]
\centering
\caption{Convergence results (errors and observed orders), grouped by DLMP-region. Within each DLMP class and each $n$, the smallest error in each norm is in bold.}
\label{tab:convergence-errors-orders-dmpgroup}
\small
\setlength{\tabcolsep}{2.2pt}
\renewcommand{\arraystretch}{0.88}
\setlength{\dashlinedash}{1.2pt}  % dot/dash length
\setlength{\dashlinegap}{1.2pt}   % gap
\setlength{\arrayrulewidth}{0.35pt}

\begin{tabular}{cll|rr|rrr|rrr}
\toprule
 &  \multicolumn{2}{c}{FV4}  & \multicolumn{2}{c}{Parameters} & \multicolumn{3}{c}{Errors} & \multicolumn{3}{c}{Orders} \\
\cmidrule(lr){4-5}\cmidrule(lr){6-8}\cmidrule(lr){9-11}
 & DLMP-region & Limiter & $n$ & $\Delta t$ & $L^1$ & $L^2$ & $L^\infty$ & $p(L^1)$ & $p(L^2)$ & $p(L^\infty)$ \\
\midrule
& \multirow{15}{*}{$N(K)\cup K$}
& GBJ-K (\cref{method:GBJ_K})      & 16  & 5.102e-03 & 9.498e-02 & 1.193e-01 & 2.346e-01 &        &        &       \\
&  & GBJ-K (\cref{method:GBJ_K})       & 32  & 2.500e-03 & 6.399e-02 & 8.292e-02 & 1.933e-01 & 0.570  & 0.525  & 0.279 \\
&  & GBJ-K  (\cref{method:GBJ_K})     & 64  & 1.245e-03 & 3.939e-02 & 5.292e-02 & 1.818e-01 & 0.700  & 0.648  & 0.088 \\
&  & GBJ-K  (\cref{method:GBJ_K})     & 128 & 6.219e-04 & 2.244e-02 & 3.106e-02 & 1.241e-01 & 0.812  & 0.769  & 0.551 \\
&  & GBJ-K  (\cref{method:GBJ_K})     & 256 & 3.108e-04 & 1.210e-02 & 1.709e-02 & 7.000e-02 & 0.891  & 0.862  & 0.827 \\
\cdashline{3-11}
&  & NKUK  (\cref{method:NKUK})      & 16  & 5.102e-03 & 8.372e-02 & 1.059e-01 & 2.210e-01 &        &        &       \\
&  & NKUK    (\cref{method:NKUK})    & 32  & 2.500e-03 & 5.472e-02 & 7.115e-02 & 1.944e-01 & 0.613  & 0.574  & 0.185 \\
&  & NKUK   (\cref{method:NKUK})     & 64  & 1.245e-03 & 3.485e-02 & 4.677e-02 & 1.818e-01 & 0.651  & 0.605  & 0.096 \\
&  & NKUK   (\cref{method:NKUK})     & 128 & 6.219e-04 & 2.098e-02 & 2.852e-02 & 1.241e-01 & 0.733  & 0.714  & 0.551 \\
&  & NKUK   (\cref{method:NKUK})     & 256 & 3.108e-04 & 1.159e-02 & 1.593e-02 & 7.000e-02 & 0.856  & 0.840  & 0.827 \\
\cdashline{3-11}
&  & Sharp-NKUK (\cref{method:sharp_NKUK})  & 16  & 5.102e-03 & \textbf{3.220e-02} & \textbf{3.995e-02} & \textbf{9.299e-02} &        &        &       \\
&  & Sharp-NKUK (\cref{method:sharp_NKUK}) & 32  & 2.500e-03 & \textbf{6.761e-03} & \textbf{1.176e-02} & \textbf{4.688e-02} & 2.252  & 1.764  & 0.988 \\
&  & Sharp-NKUK (\cref{method:sharp_NKUK}) & 64  & 1.245e-03 & \textbf{1.264e-03} & \textbf{3.289e-03} & \textbf{1.897e-02} & 2.419  & 1.838  & 1.305 \\
&  & Sharp-NKUK (\cref{method:sharp_NKUK}) & 128 & 6.219e-04 & \textbf{1.660e-04} & \textbf{7.354e-04} & \textbf{7.170e-03} & 2.929  & 2.161  & 1.404 \\
&  & Sharp-NKUK (\cref{method:sharp_NKUK}) & 256 & 3.108e-04 & \textbf{2.249e-05} & \textbf{1.723e-04} & \textbf{2.709e-03} & 2.884  & 2.093  & 1.404 \\
\midrule

& \multirow{15}{*}{$N^2(K)\cup N(K)$}
& GBJ-NKUK  (\cref{method:GBJ_NKUK})    & 16  & 5.102e-03 & 1.000e-02 & 1.691e-02 & 5.672e-02 &        &        &       \\
&  & GBJ-NKUK    (\cref{method:GBJ_NKUK})    & 32  & 2.500e-03 & 1.727e-03 & 4.378e-03 & 2.560e-02 & 2.534  & 1.949  & 1.148 \\
&  & GBJ-NKUK    (\cref{method:GBJ_NKUK})    & 64  & 1.245e-03 & 2.545e-04 & 9.904e-04 & 9.747e-03 & 2.762  & 2.144  & 1.393 \\
&  & GBJ-NKUK    (\cref{method:GBJ_NKUK})    & 128 & 6.219e-04 & 3.182e-05 & 2.066e-04 & 3.519e-03 & 3.000  & 2.261  & 1.470 \\
&  & GBJ-NKUK    (\cref{method:GBJ_NKUK})    & 256 & 3.108e-04 & 3.553e-06 & 4.017e-05 & 1.218e-03 & 3.163  & 2.363  & 1.531 \\
\cdashline{3-11}

&  & N2UNK   (\cref{method:dull_N2KUNK})      & 16  & 5.102e-03 & 8.430e-03 & 1.587e-02 & 5.393e-02 &        &        &       \\
&  & N2UNK     (\cref{method:dull_N2KUNK})    & 32  & 2.500e-03 & 1.539e-03 & 3.954e-03 & 2.402e-02 & 2.454  & 2.005  & 1.167 \\
&  & N2UNK     (\cref{method:dull_N2KUNK})    & 64  & 1.245e-03 & 2.003e-04 & 8.748e-04 & 9.179e-03 & 2.941  & 2.176  & 1.388 \\
&  & N2UNK    (\cref{method:dull_N2KUNK})     & 128 & 6.219e-04 & 2.325e-05 & 1.822e-04 & 3.345e-03 & 3.107  & 2.264  & 1.456 \\
&  & N2UNK    (\cref{method:dull_N2KUNK})     & 256 & 3.108e-04 & \textbf{2.607e-06} & 3.537e-05 & 1.145e-03 & 3.156  & 2.365  & 1.547 \\
\cdashline{3-11}

&  & Sharp-N2KUNK (\cref{method:sharp_N2KUNK}) & 16  & 5.102e-03 & \textbf{8.424e-03} & \textbf{1.586e-02} & \textbf{5.386e-02} &        &        &       \\
&  & Sharp-N2KUNK (\cref{method:sharp_N2KUNK}) & 32  & 2.500e-03 & \textbf{1.533e-03} & \textbf{3.937e-03} & \textbf{2.391e-02} & 2.458  & 2.010  & 1.172 \\
&  & Sharp-N2KUNK (\cref{method:sharp_N2KUNK}) & 64  & 1.245e-03 & \textbf{2.002e-04} & \textbf{8.727e-04} & \textbf{9.158e-03} & 2.936  & 2.174  & 1.385 \\
&  & Sharp-N2KUNK(\cref{method:sharp_N2KUNK})  & 128 & 6.219e-04 & \textbf{2.324e-05} & \textbf{1.817e-04} & \textbf{3.341e-03} & 3.107  & 2.264  & 1.455 \\
&  & Sharp-N2KUNK (\cref{method:sharp_N2KUNK}) & 256 & 3.108e-04 & 2.610e-06 & \textbf{3.524e-05} & \textbf{1.141e-03} & 3.154  & 2.366  & 1.550 \\
\midrule

& \multirow{10}{*}{$VN(K)$}
& VNK  (\cref{method:dull_VNK})     & 16  & 5.102e-03 & \textbf{9.137e-03} & \textbf{1.704e-02} & \textbf{5.788e-02} &        &        &       \\
&  & VNK    (\cref{method:dull_VNK})    & 32  & 2.500e-03 & \textbf{1.681e-03} & \textbf{4.253e-03} & \textbf{2.550e-02} & 2.442  & 2.002  & 1.182 \\
&  & VNK    (\cref{method:dull_VNK})    & 64  & 1.245e-03 & \textbf{2.170e-04} & \textbf{9.260e-04} & \textbf{9.560e-03} & 2.954  & 2.199  & 1.415 \\
&  & VNK    (\cref{method:dull_VNK})    & 128 & 6.219e-04 & \textbf{2.517e-05} & \textbf{1.899e-04} & \textbf{3.427e-03} & 3.108  & 2.286  & 1.480 \\
&  & VNK    (\cref{method:dull_VNK})    & 256 & 3.108e-04 & \textbf{2.762e-06} & \textbf{3.623e-05} & \textbf{1.159e-03} & 3.188  & 2.390  & 1.564 \\
\cdashline{3-11}

&  & Sharp-VNK (\cref{method:sharp_VNK})  & 16  & 5.102e-03 & \textbf{9.137e-03} & \textbf{1.704e-02} & \textbf{5.788e-02} &        &        &       \\
&  & Sharp-VNK (\cref{method:sharp_VNK}) & 32  & 2.500e-03 & \textbf{1.681e-03} & \textbf{4.253e-03} & \textbf{2.550e-02} & 2.442  & 2.002  & 1.182 \\
&  & Sharp-VNK (\cref{method:sharp_VNK}) & 64  & 1.245e-03 & \textbf{2.170e-04} & \textbf{9.260e-04} & \textbf{9.560e-03} & 2.954  & 2.199  & 1.415 \\
&  & Sharp-VNK (\cref{method:sharp_VNK}) & 128 & 6.219e-04 & \textbf{2.517e-05} & \textbf{1.899e-04} & \textbf{3.427e-03} & 3.108  & 2.286  & 1.480 \\
&  & Sharp-VNK (\cref{method:sharp_VNK}) & 256 & 3.108e-04 & \textbf{2.762e-06} & \textbf{3.623e-05} & \textbf{1.159e-03} & 3.188  & 2.390  & 1.564 \\
\midrule
& \multirow{5}{*}{Global}
& GlobalMP-\cite{zhang2010maximum} & 16  & 5.102e-03 & \textbf{8.431e-03} & \textbf{1.586e-02} & \textbf{5.385e-02} &        &        &       \\
&  & GlobalMP-\cite{zhang2010maximum} & 32  & 2.500e-03 & \textbf{1.529e-03} & \textbf{3.921e-03} & \textbf{2.381e-02} & 2.463  & 2.016  & 1.178 \\
&  & GlobalMP-\cite{zhang2010maximum} & 64  & 1.245e-03 & \textbf{1.998e-04} & \textbf{8.681e-04} & \textbf{9.121e-03} & 2.936  & 2.175  & 1.384 \\
&  & GlobalMP-\cite{zhang2010maximum} & 128 & 6.219e-04 & \textbf{2.292e-05} & \textbf{1.794e-04} & \textbf{3.312e-03} & 3.124  & 2.274  & 1.462 \\
&  & GlobalMP-\cite{zhang2010maximum} & 256 & 3.108e-04 & \textbf{2.555e-06} & \textbf{3.455e-05} & \textbf{1.126e-03} & 3.165  & 2.377  & 1.556 \\
\midrule

& \multirow{10}{*}{NONE}
& WZL-\cite{WANG2004716} & 16  & 5.102e-03 & 9.941e-03 & 1.683e-02 & 5.652e-02 &        &        &       \\
&  & WZL-\cite{WANG2004716} & 32  & 2.500e-03 & 1.726e-03 & 4.376e-03 & 2.559e-02 & 2.526  & 1.943  & 1.143 \\
&  & WZL-\cite{WANG2004716} & 64  & 1.245e-03 & 2.545e-04 & 9.903e-04 & 9.746e-03 & 2.762  & 2.144  & 1.393 \\
&  & WZL-\cite{WANG2004716} & 128 & 6.219e-04 & 3.182e-05 & 2.066e-04 & 3.519e-03 & 3.000  & 2.261  & 1.470 \\
&  & WZL-\cite{WANG2004716} & 256 & 3.108e-04 & 3.553e-06 & 4.017e-05 & 1.218e-03 & 3.163  & 2.363  & 1.531 \\
\cdashline{3-11}

&  & off                  & 16  & 5.102e-03 & \textbf{4.581e-03} & \textbf{7.492e-03} & \textbf{2.451e-02} &        &        &       \\
&  & off                  & 32  & 2.500e-03 & \textbf{5.825e-04} & \textbf{1.256e-03} & \textbf{5.833e-03} & 2.975  & 2.577  & 2.071 \\
&  & off                  & 64  & 1.245e-03 & \textbf{3.554e-05} & \textbf{9.447e-05} & \textbf{6.035e-04} & 4.035  & 3.733  & 3.273 \\
&  & off                  & 128 & 6.219e-04 & \textbf{1.390e-06} & \textbf{4.007e-06} & \textbf{2.964e-05} & 4.677  & 4.559  & 4.348 \\
&  & off                  & 256 & 3.108e-04 & \textbf{4.542e-08} & \textbf{1.339e-07} & \textbf{1.012e-06} & 4.935  & 4.904  & 4.871 \\
\bottomrule
\end{tabular}

\normalsize
\end{table}

\begin{table}[H]
\centering
\caption{Convergence results unidirectional flow (errors and observed orders), grouped by DLMP-region. Within each region and each resolution $n$, the smallest error in each norm is in bold.}
\label{tab:convergence-errors-orders-dmpgroup_constant_flow}
\small
\setlength{\tabcolsep}{2.2pt}
\renewcommand{\arraystretch}{0.88}
\setlength{\dashlinedash}{1.2pt}  % dot/dash length
\setlength{\dashlinegap}{1.2pt}   % gap
\setlength{\arrayrulewidth}{0.35pt}

\begin{tabular}{cll|rr|rrr|rrr}
\toprule
 &  \multicolumn{2}{c}{FV4}  & \multicolumn{2}{c}{Parameters} & \multicolumn{3}{c}{Errors} & \multicolumn{3}{c}{Orders} \\
\cmidrule(lr){4-5}\cmidrule(lr){6-8}\cmidrule(lr){9-11}
 & DLMP-region & Limiter & $n$ & $\Delta t$ & $L^1$ & $L^2$ & $L^\infty$ & $p(L^1)$ & $p(L^2)$ & $p(L^\infty)$ \\
\midrule

 & \multirow{15}{*}{$N(K)\cup K$} & GBJ-K (\cref{method:GBJ_K}) & 16 & 1.953e-03 & 2.010e-01 & 2.447e-01 & 4.715e-01 &  &  &  \\
 &  & GBJ-K (\cref{method:GBJ_K}) & 32 & 9.766e-04 & 1.740e-01 & 2.140e-01 & 4.244e-01 & 0.208 & 0.194 & 0.152 \\
 &  & GBJ-K (\cref{method:GBJ_K}) & 64 & 4.883e-04 & 1.260e-01 & 1.553e-01 & 3.100e-01 & 0.466 & 0.463 & 0.453 \\
 &  & GBJ-K (\cref{method:GBJ_K}) & 128 & 2.441e-04 & 7.794e-02 & 9.614e-02 & 1.922e-01 & 0.692 & 0.692 & 0.689 \\
 &  & GBJ-K (\cref{method:GBJ_K}) & 256 & 1.221e-04 & 4.368e-02 & 5.388e-02 & 1.077e-01 & 0.836 & 0.835 & 0.835 \\
\cdashline{3-11}
 &  & NKUK (\cref{method:NKUK}) & 16 & 1.953e-03 & 1.993e-01 & 2.427e-01 & 4.672e-01 &  &  &  \\
 &  & NKUK (\cref{method:NKUK}) & 32 & 9.766e-04 & 1.679e-01 & 2.060e-01 & 4.079e-01 & 0.248 & 0.236 & 0.196 \\
 &  & NKUK (\cref{method:NKUK}) & 64 & 4.883e-04 & 1.197e-01 & 1.464e-01 & 2.902e-01 & 0.488 & 0.493 & 0.492 \\
 &  & NKUK (\cref{method:NKUK}) & 128 & 2.441e-04 & 7.331e-02 & 8.877e-02 & 1.758e-01 & 0.707 & 0.722 & 0.723 \\
 &  & NKUK (\cref{method:NKUK}) & 256 & 1.221e-04 & 4.143e-02 & 4.911e-02 & 9.650e-02 & 0.823 & 0.854 & 0.865 \\
\cdashline{3-11}
 &  & Sharp-NKUK (\cref{method:sharp_NKUK}) & 16 & 1.953e-03 & \textbf{4.136e-02} & \textbf{6.170e-02} & \textbf{1.754e-01} &  &  &  \\
 &  & Sharp-NKUK (\cref{method:sharp_NKUK}) & 32 & 9.766e-04 & \textbf{1.212e-02} & \textbf{1.935e-02} & \textbf{7.740e-02} & 1.771 & 1.673 & 1.180 \\
 &  & Sharp-NKUK (\cref{method:sharp_NKUK}) & 64 & 4.883e-04 & \textbf{1.868e-03} & \textbf{4.723e-03} & \textbf{2.990e-02} & 2.698 & 2.035 & 1.372 \\
 &  & Sharp-NKUK (\cref{method:sharp_NKUK}) & 128 & 2.441e-04 & \textbf{2.831e-04} & \textbf{1.170e-03} & \textbf{1.134e-02} & 2.722 & 2.013 & 1.398 \\
 &  & Sharp-NKUK (\cref{method:sharp_NKUK}) & 256 & 1.221e-04 & \textbf{4.355e-05} & \textbf{2.937e-04} & \textbf{4.441e-03} & 2.700 & 1.994 & 1.353 \\
\midrule

 & \multirow{15}{*}{$N^2(K)\cup N(K)$} & GBJ-NKUK (\cref{method:GBJ_NKUK}) & 16 & 1.953e-03 & 2.161e-02 & 3.104e-02 & 9.588e-02 &  &  &  \\
 &  & GBJ-NKUK (\cref{method:GBJ_NKUK}) & 32 & 9.766e-04 & 3.157e-03 & 6.476e-03 & 3.530e-02 & 2.775 & 2.261 & 1.441 \\
 &  & GBJ-NKUK (\cref{method:GBJ_NKUK}) & 64 & 4.883e-04 & 3.768e-04 & 1.337e-03 & 1.317e-02 & 3.067 & 2.276 & 1.422 \\
 &  & GBJ-NKUK (\cref{method:GBJ_NKUK}) & 128 & 2.441e-04 & 5.169e-05 & 3.096e-04 & 5.197e-03 & 2.866 & 2.111 & 1.342 \\
 &  & GBJ-NKUK (\cref{method:GBJ_NKUK}) & 256 & 1.221e-04 & 6.791e-06 & 6.983e-05 & 1.949e-03 & 2.928 & 2.149 & 1.415 \\
\cdashline{3-11}
 &  & N2UNK (\cref{method:dull_N2KUNK}) & 16 & 1.953e-03 & \textbf{1.466e-02} & \textbf{2.243e-02} & \textbf{7.338e-02} &  &  &  \\
 &  & N2UNK (\cref{method:dull_N2KUNK}) & 32 & 9.766e-04 & \textbf{2.016e-03} & \textbf{4.415e-03} & \textbf{2.663e-02} & 2.862 & 2.345 & 1.463 \\
 &  & N2UNK (\cref{method:dull_N2KUNK}) & 64 & 4.883e-04 & \textbf{2.338e-04} & \textbf{9.020e-04} & \textbf{9.820e-03} & 3.108 & 2.291 & 1.439 \\
 &  & N2UNK (\cref{method:dull_N2KUNK}) & 128 & 2.441e-04 & 3.100e-05 & 2.082e-04 & 3.889e-03 & 2.915 & 2.115 & 1.336 \\
 &  & N2UNK (\cref{method:dull_N2KUNK}) & 256 & 1.221e-04 & \textbf{4.449e-06} & \textbf{4.960e-05} & \textbf{1.532e-03} & 2.801 & 2.070 & 1.344 \\
\cdashline{3-11}
 &  & Sharp-N2KUNK (\cref{method:sharp_N2KUNK}) & 16 & 1.953e-03 & \textbf{1.466e-02} & \textbf{2.243e-02} & \textbf{7.338e-02} &  &  &  \\
 &  & Sharp-N2KUNK (\cref{method:sharp_N2KUNK}) & 32 & 9.766e-04 & \textbf{2.016e-03} & \textbf{4.415e-03} & \textbf{2.663e-02} & 2.862 & 2.345 & 1.463 \\
 &  & Sharp-N2KUNK (\cref{method:sharp_N2KUNK}) & 64 & 4.883e-04 & \textbf{2.338e-04} & 9.023e-04 & 9.826e-03 & 3.108 & 2.291 & 1.438 \\
 &  & Sharp-N2KUNK (\cref{method:sharp_N2KUNK}) & 128 & 2.441e-04 & \textbf{3.084e-05} & \textbf{2.072e-04} & \textbf{3.876e-03} & 2.923 & 2.122 & 1.342 \\
 &  & Sharp-N2KUNK (\cref{method:sharp_N2KUNK}) & 256 & 1.221e-04 & 4.477e-06 & 5.003e-05 & 1.548e-03 & 2.784 & 2.050 & 1.324 \\
\midrule

 & \multirow{10}{*}{$VN(K)$} & VNK (\cref{method:dull_VNK}) & 16 & 1.953e-03 & \textbf{1.532e-02} & \textbf{2.329e-02} & \textbf{7.599e-02} &  &  &  \\
 &  & VNK (\cref{method:dull_VNK}) & 32 & 9.766e-04 & \textbf{2.138e-03} & \textbf{4.643e-03} & \textbf{2.776e-02} & 2.842 & 2.327 & 1.453 \\
 &  & VNK (\cref{method:dull_VNK}) & 64 & 4.883e-04 & \textbf{2.346e-04} & \textbf{8.957e-04} & \textbf{9.725e-03} & 3.188 & 2.374 & 1.513 \\
 &  & VNK (\cref{method:dull_VNK}) & 128 & 2.441e-04 & \textbf{2.979e-05} & \textbf{2.009e-04} & \textbf{3.777e-03} & 2.977 & 2.157 & 1.364 \\
 &  & VNK (\cref{method:dull_VNK}) & 256 & 1.221e-04 & \textbf{4.315e-06} & \textbf{4.857e-05} & \textbf{1.509e-03} & 2.787 & 2.048 & 1.323 \\
\cdashline{3-11}
 &  & Sharp-VNK (\cref{method:sharp_VNK}) & 16 & 1.953e-03 & \textbf{1.532e-02} & \textbf{2.329e-02} & \textbf{7.599e-02} &  &  &  \\
 &  & Sharp-VNK (\cref{method:sharp_VNK}) & 32 & 9.766e-04 & \textbf{2.138e-03} & \textbf{4.643e-03} & \textbf{2.776e-02} & 2.842 & 2.327 & 1.453 \\
 &  & Sharp-VNK (\cref{method:sharp_VNK}) & 64 & 4.883e-04 & \textbf{2.346e-04} & \textbf{8.957e-04} & \textbf{9.725e-03} & 3.188 & 2.374 & 1.513 \\
 &  & Sharp-VNK (\cref{method:sharp_VNK}) & 128 & 2.441e-04 & \textbf{2.979e-05} & \textbf{2.009e-04} & \textbf{3.777e-03} & 2.977 & 2.157 & 1.364 \\
 &  & Sharp-VNK (\cref{method:sharp_VNK}) & 256 & 1.221e-04 & \textbf{4.315e-06} & \textbf{4.857e-05} & \textbf{1.509e-03} & 2.787 & 2.048 & 1.323 \\
\midrule

 & \multirow{5}{*}{Global} & GlobalMP-\cite{zhang2010maximum} & 16 & 1.953e-03 & \textbf{1.429e-02} & \textbf{2.192e-02} & \textbf{7.185e-02} &  &  &  \\
 &  & GlobalMP-\cite{zhang2010maximum} & 32 & 9.766e-04 & \textbf{1.936e-03} & \textbf{4.254e-03} & \textbf{2.579e-02} & 2.884 & 2.365 & 1.478 \\
 &  & GlobalMP-\cite{zhang2010maximum} & 64 & 4.883e-04 & \textbf{2.219e-04} & \textbf{8.574e-04} & \textbf{9.430e-03} & 3.125 & 2.311 & 1.451 \\
 &  & GlobalMP-\cite{zhang2010maximum} & 128 & 2.441e-04 & \textbf{2.790e-05} & \textbf{1.889e-04} & \textbf{3.609e-03} & 2.991 & 2.182 & 1.386 \\
 &  & GlobalMP-\cite{zhang2010maximum} & 256 & 1.221e-04 & \textbf{3.796e-06} & \textbf{4.304e-05} & \textbf{1.391e-03} & 2.878 & 2.134 & 1.375 \\
\midrule

 & \multirow{10}{*}{NONE} & WZL-\cite{WANG2004716} & 16 & 1.953e-03 & 2.161e-02 & 3.104e-02 & 9.588e-02 &  &  &  \\
 &  & WZL-\cite{WANG2004716} & 32 & 9.766e-04 & 3.157e-03 & 6.476e-03 & 3.530e-02 & 2.775 & 2.261 & 1.441 \\
 &  & WZL-\cite{WANG2004716} & 64 & 4.883e-04 & 3.768e-04 & 1.337e-03 & 1.317e-02 & 3.067 & 2.276 & 1.422 \\
 &  & WZL-\cite{WANG2004716} & 128 & 2.441e-04 & 5.170e-05 & 3.098e-04 & 5.206e-03 & 2.866 & 2.110 & 1.339 \\
 &  & WZL-\cite{WANG2004716} & 256 & 1.221e-04 & 6.877e-06 & 7.070e-05 & 1.973e-03 & 2.910 & 2.131 & 1.400 \\
\cdashline{3-11}
 &  & off & 16 & 1.953e-03 & \textbf{1.068e-03} & \textbf{1.228e-03} & \textbf{2.234e-03} &  &  &  \\
 &  & off & 32 & 9.766e-04 & \textbf{6.036e-05} & \textbf{6.781e-05} & \textbf{1.136e-04} & 4.145 & 4.178 & 4.298 \\
 &  & off & 64 & 4.883e-04 & \textbf{3.654e-06} & \textbf{4.071e-06} & \textbf{6.348e-06} & 4.046 & 4.058 & 4.162 \\
 &  & off & 128 & 2.441e-04 & \textbf{2.264e-07} & \textbf{2.517e-07} & \textbf{3.750e-07} & 4.012 & 4.016 & 4.081 \\
 &  & off & 256 & 1.221e-04 & \textbf{1.412e-08} & \textbf{1.569e-08} & \textbf{2.279e-08} & 4.003 & 4.004 & 4.040 \\
\bottomrule
\end{tabular}
\normalsize
\end{table}

\begin{table}[H]
\centering
\small
\setlength{\tabcolsep}{3.5pt}
\caption{Region-wise relative $L^2$ errors and observed rates for LeVeque initial conditions under solid body rotation, grouped by enforced DLMP region. Within each DLMP group and each resolution $n$, the smallest error in each region column is in bold (ties in bold). Observed convergence is limited by the regularity of initial condition: the cone is $C^0$ so at best first order; the cosine hump is at most $C^1$ so at best second order; and the slotted cylinder is discontinuous ($L^\infty$ only), for which one expects, order $1/2$ in $L^2$.}
\label{tab:fv4-dlmp-relL2-rates}
\resizebox{\textwidth}{!}{%
\begin{tabular}{lrrcccccccc}
\toprule
& & &
\multicolumn{4}{c}{Relative $L^2$ errors} &
\multicolumn{4}{c}{Observed rates} \\
\cmidrule(lr){4-7}\cmidrule(lr){8-11}
Method & $n_x$ & $\Delta t$
& Cone & Cosine & Slotted & Whole
& $p_{\text{cone}}$ & $p_{\text{cos}}$ & $p_{\text{slot}}$ & $p_{\text{whole}}$ \\
\midrule

\multicolumn{11}{l}{\textbf{DLMP region: NONE (unlimited)}}\\
\addlinespace[2pt]
off  & 64  & 5.051e-04 & \textbf{3.551e-02} & \textbf{2.675e-02} & \textbf{2.640e-01} & \textbf{3.152e-01} &  &  &  &  \\
off  & 128 & 2.506e-04 & \textbf{1.358e-02} & \textbf{5.109e-03} & \textbf{1.855e-01} & \textbf{2.281e-01} & 1.386 & 2.388 & 0.510 & 0.467 \\

\addlinespace[4pt]
\hline
\addlinespace[4pt]

\multicolumn{11}{l}{\textbf{DLMP region: Global }}\\
\addlinespace[2pt]
GlobalMP-\cite{zhang2010maximum} & 64  & 5.051e-04 & \textbf{4.049e-02} & \textbf{1.634e-02} & \textbf{2.621e-01} & \textbf{3.276e-01} &  &  &  &  \\
GlobalMP-\cite{zhang2010maximum} & 128 & 2.506e-04 & \textbf{1.487e-02} & \textbf{2.966e-03} & \textbf{1.834e-01} & \textbf{2.306e-01} & 1.445 & 2.462 & 0.515 & 0.507 \\

\addlinespace[4pt]
\hline
\addlinespace[4pt]

\multicolumn{11}{l}{\textbf{DLMP region: $N(K)\cup K$}}\\
\addlinespace[2pt]
GBJ-K      & 64  & 5.051e-04 & 7.536e-01 & 7.171e-01 & 7.700e-01 & 8.223e-01 &  &  &  &  \\
GBJ-K      & 128 & 2.506e-04 & 6.289e-01 & 7.018e-01 & 6.397e-01 & 7.181e-01 & 0.261 & 0.031 & 0.267 & 0.195 \\
\addlinespace[2pt]
NKUK       & 64  & 5.051e-04 & 7.165e-01 & 7.086e-01 & 7.336e-01 & 7.931e-01 &  &  &  &  \\
NKUK       & 128 & 2.506e-04 & 5.819e-01 & 6.721e-01 & 6.037e-01 & 6.892e-01 & 0.300 & 0.076 & 0.281 & 0.203 \\
\addlinespace[2pt]
Sharp-NKUK & 64  & 5.051e-04 & \textbf{1.428e-01} & \textbf{1.719e-01} & \textbf{3.623e-01} & \textbf{4.664e-01} &  &  &  &  \\
Sharp-NKUK & 128 & 2.506e-04 & \textbf{5.629e-02} & \textbf{5.231e-02} & \textbf{2.287e-01} & \textbf{3.074e-01} & 1.343 & 1.717 & 0.664 & 0.601 \\

\addlinespace[4pt]
\hline
\addlinespace[4pt]

\multicolumn{11}{l}{\textbf{DLMP region: $N^2(K)\cup N(K)$}}\\
\addlinespace[2pt]
GBJ-NKUK     & 64  & 5.051e-04 & 7.518e-02 & 6.063e-02 & 2.922e-01 & 3.699e-01 &  &  &  &  \\
GBJ-NKUK     & 128 & 2.506e-04 & 2.871e-02 & 1.650e-02 & 1.975e-01 & 2.470e-01 & 1.389 & 1.878 & 0.566 & 0.583 \\
\addlinespace[2pt]
N2UNK        & 64  & 5.051e-04 & 6.538e-02 & 4.742e-02 & 2.812e-01 & 3.531e-01 &  &  &  &  \\
N2UNK        & 128 & 2.506e-04 & 2.453e-02 & 1.222e-02 & 1.946e-01 & 2.416e-01 & 1.414 & 1.957 & 0.531 & 0.547 \\
\addlinespace[2pt]
Sharp-N2KUNK & 64  & 5.051e-04 & \textbf{6.477e-02} & \textbf{4.709e-02} & \textbf{2.793e-01} & \textbf{3.513e-01} &  &  &  &  \\
Sharp-N2KUNK & 128 & 2.506e-04 & \textbf{2.383e-02} & \textbf{1.164e-02} & \textbf{1.935e-01} & \textbf{2.400e-01} & 1.443 & 2.016 & 0.529 & 0.550 \\

\addlinespace[4pt]
\hline
\addlinespace[4pt]

\multicolumn{11}{l}{\textbf{DLMP region: $VN(K)$}}\\
\addlinespace[2pt]
VNK       & 64  & 5.051e-04 & 6.630e-02 & \textbf{4.708e-02} & 2.876e-01 & 3.612e-01 &  &  &  &  \\
VNK       & 128 & 2.506e-04 & 2.519e-02 & \textbf{1.257e-02} & \textbf{1.957e-01} & 2.434e-01 & 1.396 & 1.905 & 0.555 & 0.569 \\
\addlinespace[2pt]
Sharp-VNK & 64  & 5.051e-04 & \textbf{6.627e-02} & 4.709e-02 & \textbf{2.874e-01} & \textbf{3.608e-01} &  &  &  &  \\
Sharp-VNK & 128 & 2.506e-04 & \textbf{2.517e-02} & \textbf{1.257e-02} & \textbf{1.957e-01} & \textbf{2.433e-01} & 1.397 & 1.905 & 0.554 & 0.568 \\

\bottomrule
\end{tabular}%
}
\end{table}

% \begingroup
% \setlength{\floatsep}{4pt}
% \setlength{\textfloatsep}{4pt}
% \setlength{\intextsep}{4pt}

% \begin{figure}[H]
% \centering
% \includegraphics[width=1.0\linewidth]{Figures/midface_joined_slices_leveque_sbr_False_SSPRK44_n128_CFL0p4_off-Sharp-NKUK-NKUK-GBJ-K.pdf}
% \caption{FV4, slice through LeVeque shapes, $128^2$ resolution. For Sharp-NKUK, NKUK, GBJ-K limiters. Including the unlimited scheme (off) for general comparison.}
% \label{fig:slice128_FV4_NKUK}
% \end{figure}
% \vspace{-10pt}

% \begin{figure}[H]
% \centering
% \includegraphics[width=1.0\linewidth]{Figures/midface_joined_slices_leveque_sbr_False_SSPRK44_n128_CFL0p4_Sharp-NVK-NVK-GlobalMP.pdf}
% \caption{FV4, slice through LeVeque shapes, $128^2$ resolution. For Sharp-VNK, VNK, and GlobalMP limiters.}
% \label{fig:slice128_FV4_NVK}
% \end{figure}

% \vspace{-10pt}

% \begin{figure}[H]
% \centering
% \includegraphics[width=1.0\linewidth]{Figures/midface_joined_slices_leveque_sbr_False_SSPRK44_n128_CFL0p4_Sharp-N2KUNK-N2UNK-GBJ-NKUK.pdf}
% \caption{FV4, slice through LeVeque shapes, $128^2$ resolution. For Sharp-N2KUNK, N2UNK, GBJ-NKUK limiters.}
% \label{fig:slice128_FV4_N2KUNK}
% \end{figure}

% \vspace{-10pt}

% \begin{figure}[H]
% \centering
% \includegraphics[width=1.0\linewidth]{Figures/midface_joined_slices_leveque_sbr_False_SSPRK44_n128_CFL0p4_Sharp-NVK-NVK-GlobalMP-Sharp-NKUK-Sharp-N2KUNK-N2UNK-GBJ-NKUK.pdf}
% \caption{Slice through LeVeque shapes, $128^2$ resolution. For many higher-order limiters.}
% \label{fig:slice128_FV4_successful}
% \end{figure}

% \endgroup
\begingroup
\setlength{\floatsep}{4pt}
\setlength{\textfloatsep}{4pt}
\setlength{\intextsep}{4pt}

\begin{figure}[p]
\centering

\includegraphics[height=0.22\textheight]{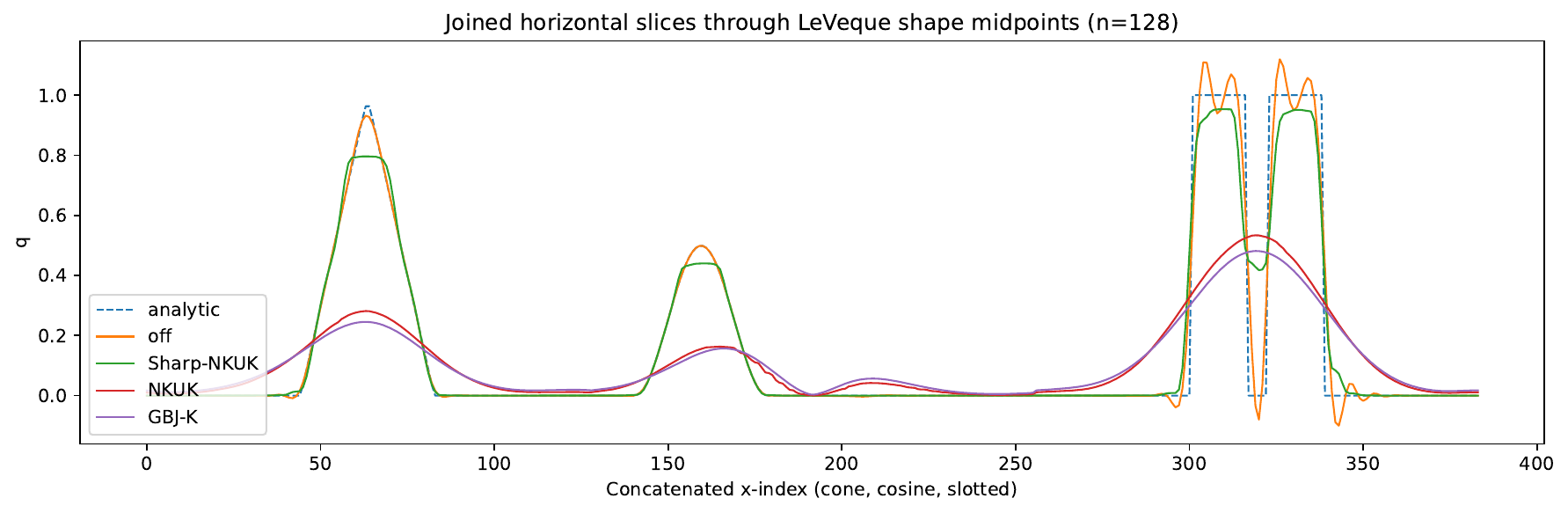}

\vspace{0.3em}

\includegraphics[height=0.22\textheight]{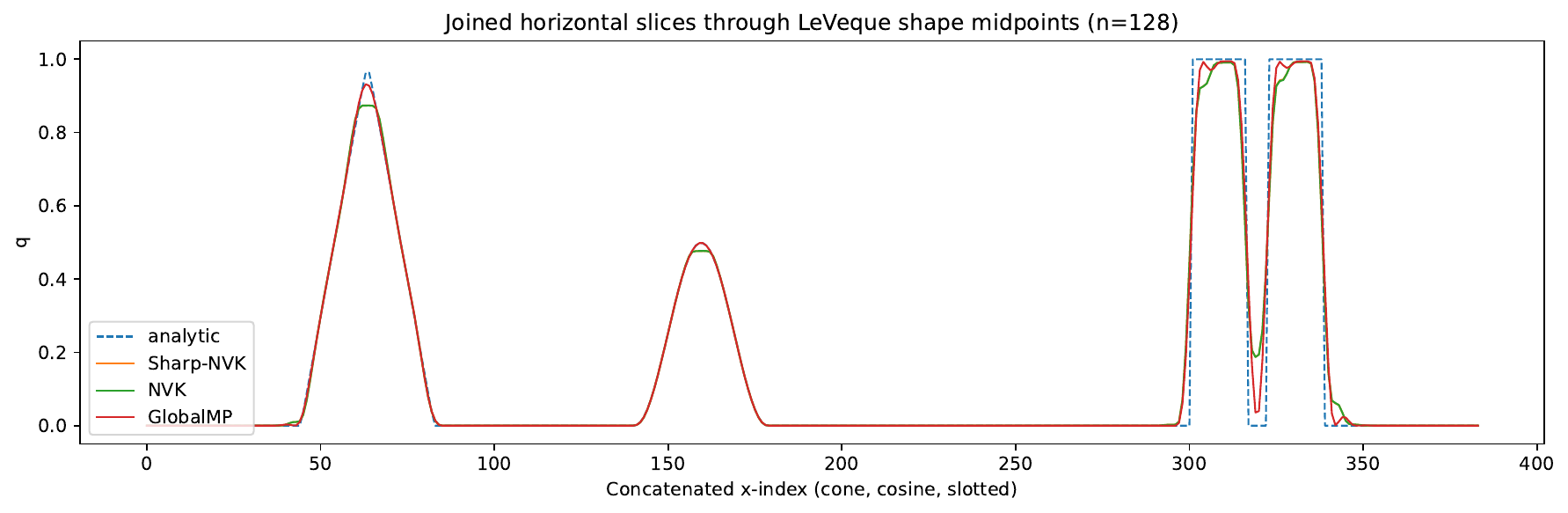}

\vspace{0.3em}

\includegraphics[height=0.22\textheight]{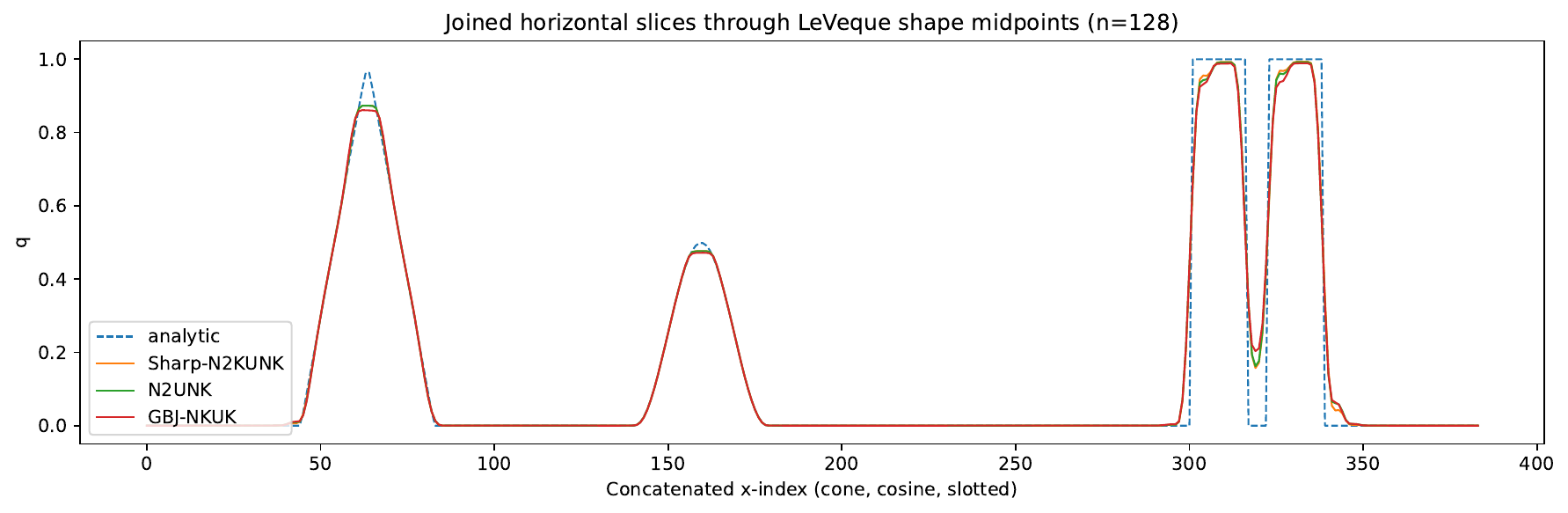}

\vspace{0.3em}

\includegraphics[height=0.22\textheight]{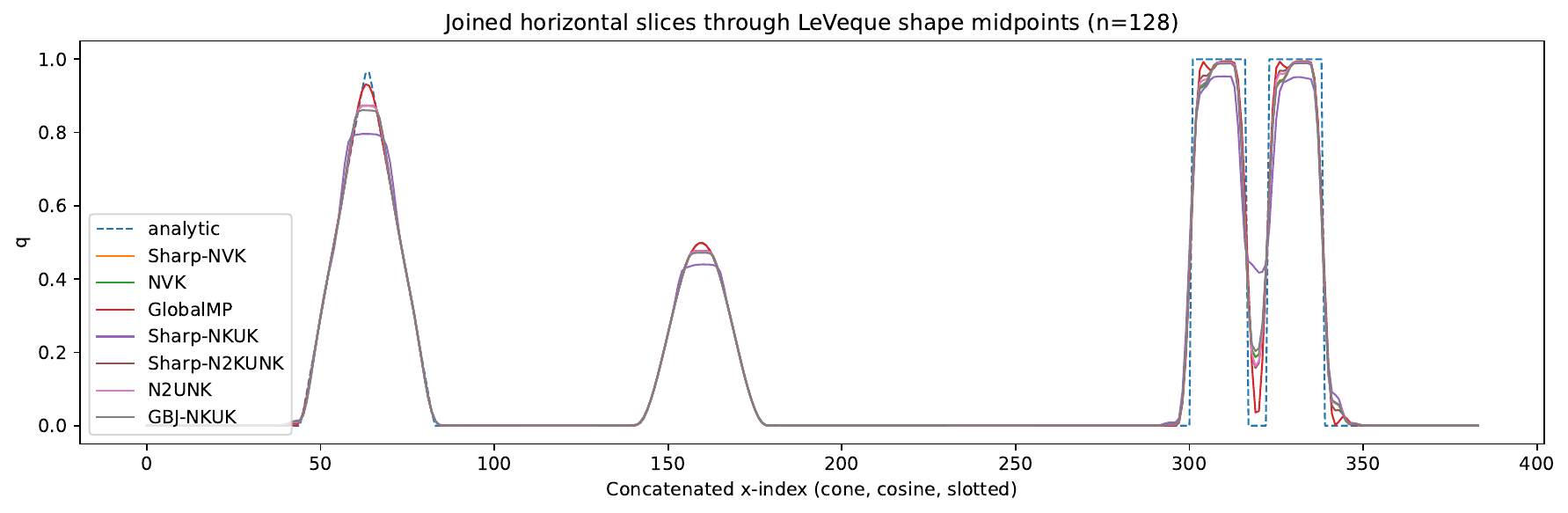}

\caption{Mid-face joined slices through the LeVeque shapes at $128^2$ resolution using the FV4 scheme.
From top to bottom: NKUK-based limiters (including the unlimited scheme), VNK-based limiters and GlobalMP,
$N^2(K)\cup N(K)$-based limiters, and a combined comparison of several higher-order limiters.}
\label{fig:slice128_FV4_all}
\end{figure}

\endgroup

\subsection{FV4: Modified initial conditions}
One might conclude from the figures in the previous experiment that the global boundedness limiter \cite{zhang2010maximum} is preferable to any of the locally bounded limiters in terms of shape preservation and peak preservation. Although the LeVeque test case is widely used, it over-represents the capability of the global boundedness limiters peak retention and shape preservation ability. This is due to two coincidental features of the initial condition. Firstly, the peak of the cone diffuses below the height of the slotted cylinder quickly, and from that point onwards remains unlimited. Secondly, the slotted cylinder contains all of its discontinuities exactly at global maxima and minima, and is therefore limited much more strongly than a general profile would be. To remove this bias, we modify the initial condition in the following experiment by halving the amplitude of the slotted cylinder and reversing the sign of the cosine bump, while leaving all other aspects of the setup unchanged.

\begin{figure}[H]
    \centering

    \begin{subfigure}[t]{0.48\linewidth}
        \centering
        \includegraphics[
            width=\linewidth,
            trim={0cm 3.5cm 0cm 4.5cm},clip
        ]{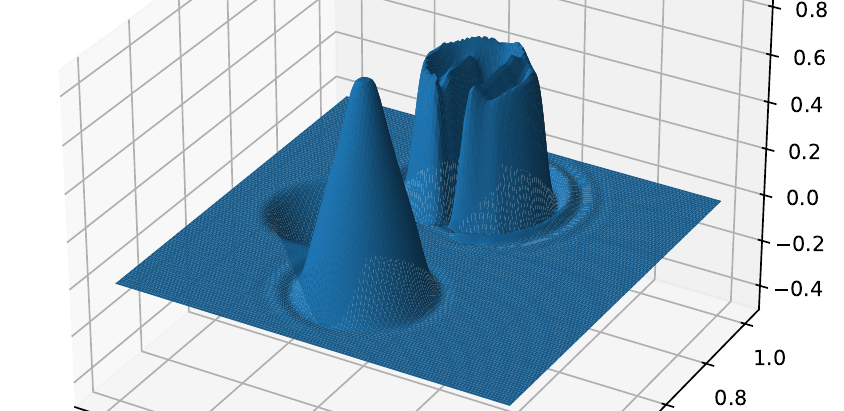}
        \caption{Global MP}
        \label{fig:gmp_new_ic}
    \end{subfigure}
    \begin{subfigure}[t]{0.48\linewidth}
        \centering
        \includegraphics[
            width=\linewidth,
            trim={0cm 3.5cm 0cm 4.5cm},clip
        ]{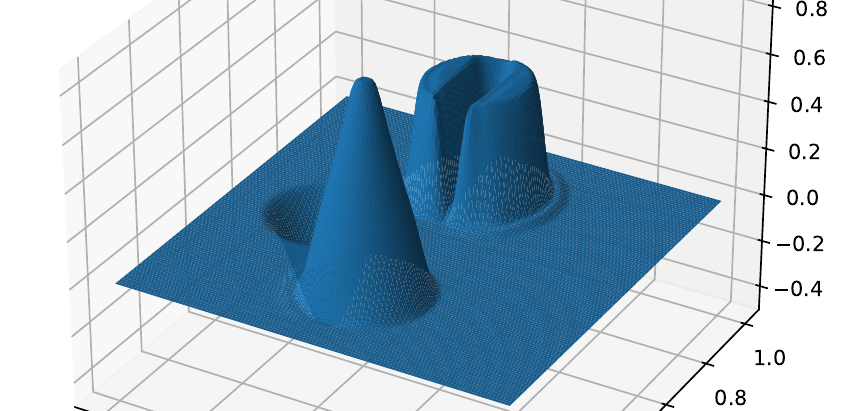}
        \caption{Sharp-N2KUNK}

    \end{subfigure}\\

    \begin{subfigure}[t]{0.48\linewidth}
        \centering
        \includegraphics[
            width=\linewidth,
            trim={0cm 3.5cm 0cm 4.5cm},clip
        ]{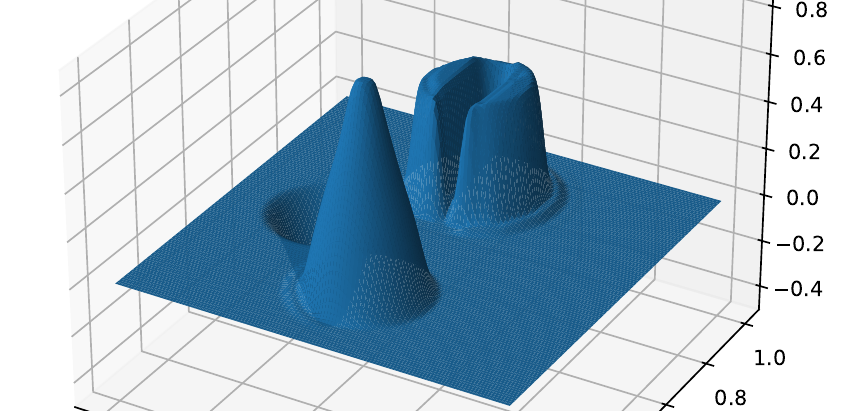}
    \caption{Sharp-VNK}
    \end{subfigure}
    \begin{subfigure}[t]{0.48\linewidth}
        \centering
        \includegraphics[
            width=\linewidth,
            trim={0cm 3.5cm 0cm 4.5cm},clip
        ]{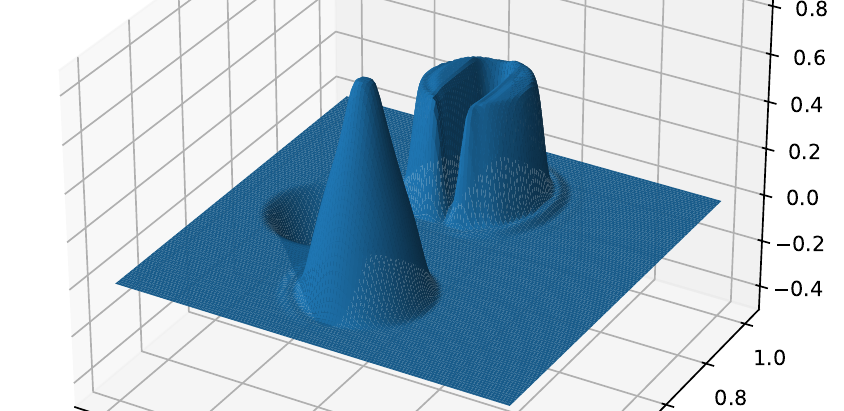}
        \caption{N2KUNK}
    \end{subfigure}\\
    \begin{subfigure}[t]{0.48\linewidth}
        \centering
        \includegraphics[
            width=\linewidth,
            trim={0cm 3.5cm 0cm 4.5cm},clip
        ]{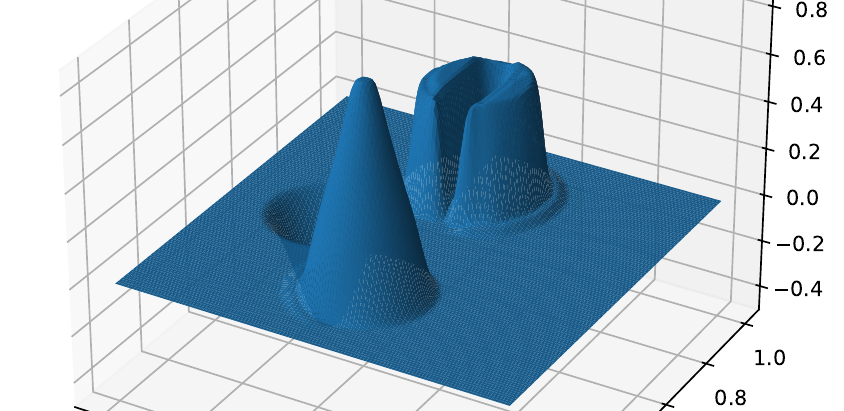}
        \caption{VNK}
    \end{subfigure}
    \begin{subfigure}[t]{0.48\linewidth}
        \centering
        \includegraphics[
            width=\linewidth,
            trim={0cm 3.5cm 0cm 4.5cm},clip
        ]{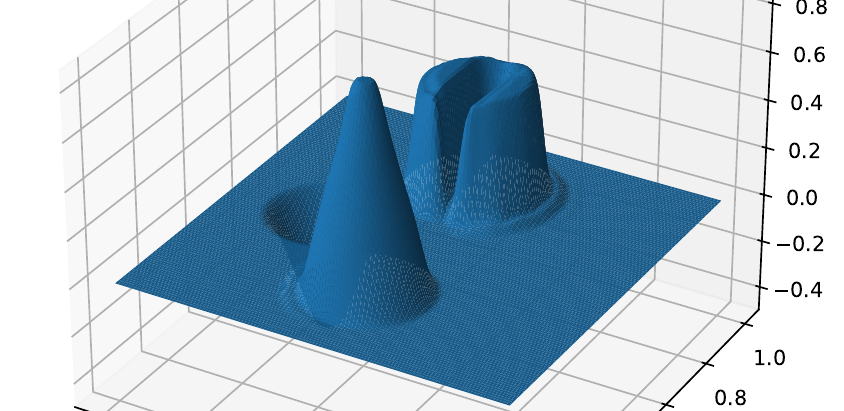}
        \caption{GBJ-NKUK}
    \end{subfigure}

    \caption{Final-time 3D surface plots for the modified LeVeque-new initial condition under solid-body rotation.
    The figure compares global, VNK-based, and $N^2(K)\cup N(K)$ DLMP limiters, highlighting differences in peak preservation and shape fidelity.}
    \label{fig:leveque-new-3d-final-n128}
\end{figure}

In \cref{fig:gmp_new_ic} we plot the Global maximum principle limiter \cite{zhang2010maximum}, which exhibits both local over and undershoots of the reduced slotted cylinder and is notably worse at shape preservation that the other limiters in \cref{fig:leveque-new-3d-final-n128}. Furthermore, in \cref{fig:gmp_new_ic} the peak retention ability of the GlobalMP limiter is near identical to the local boundedness limiters. Only the Generalised-Barth-Jespersen-$N(K)\cup K$ limiter exhibits marginally worse peak retention.

\begin{figure}[H]
    \centering
    \includegraphics[width=0.95\linewidth]{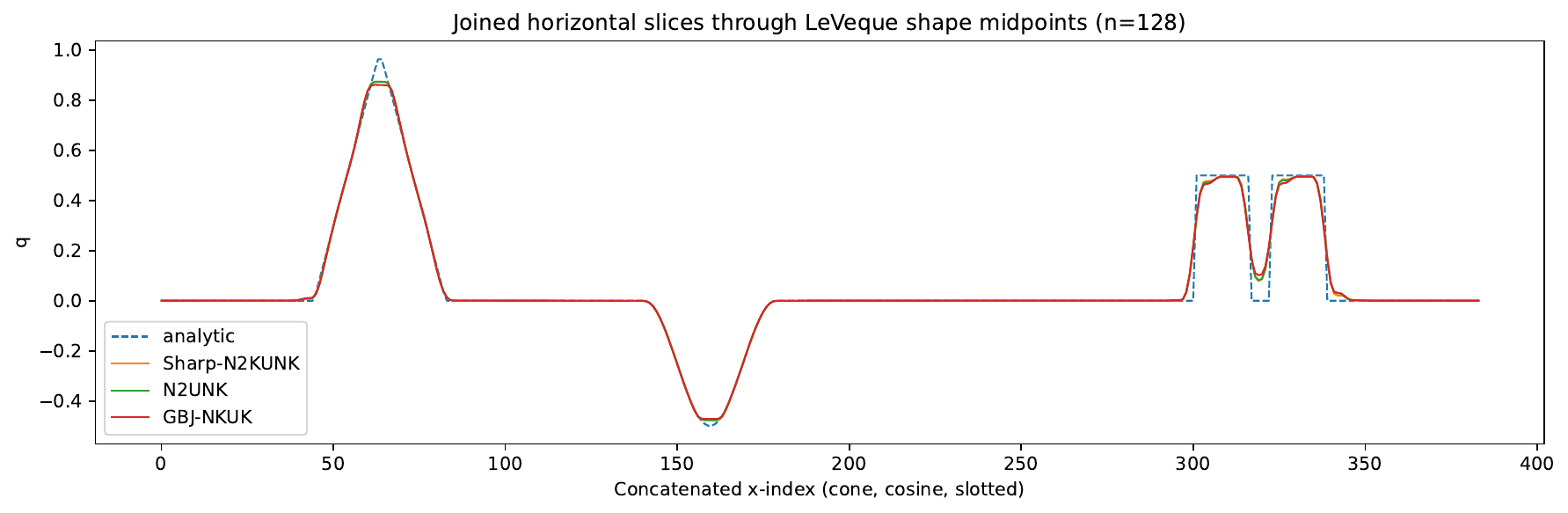}

    \vspace{0.5em}

    \includegraphics[width=0.95\linewidth]{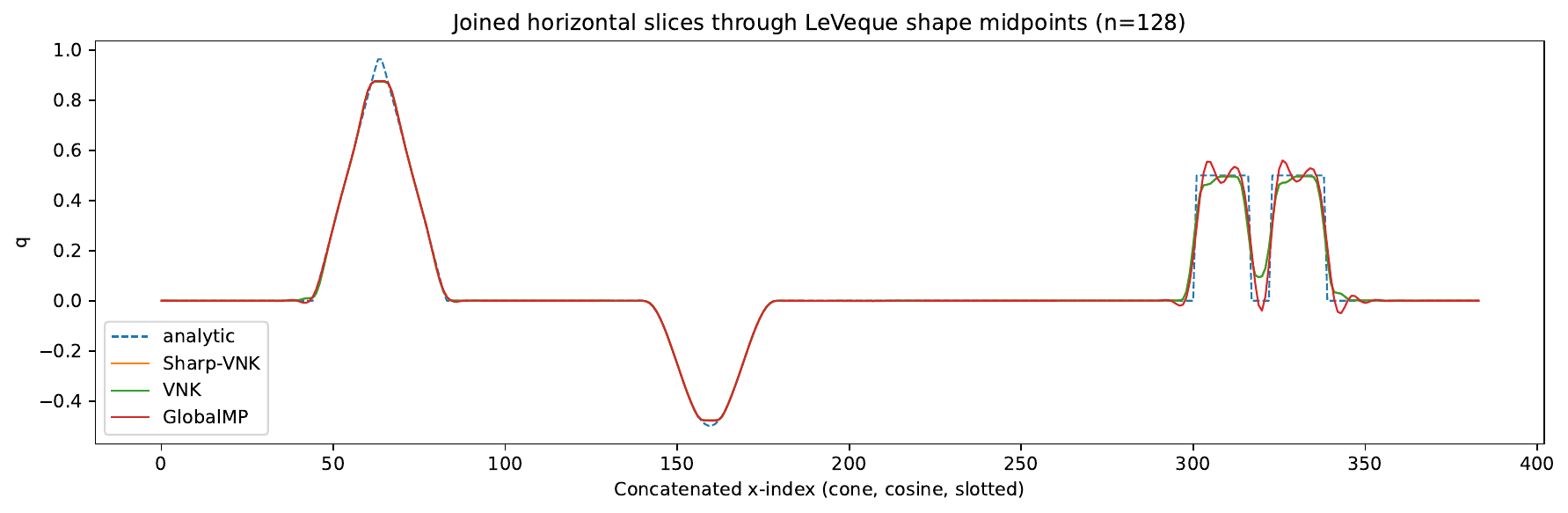}

    \caption{Mid-face joined slices for the LeVeque-new initial condition under solid-body rotation (SBR).
    Top: $128^2$ resolution comparing Sharp-N2KUNK, N2UNK, GBJ-NKUK, and NKUK limiters.
    Bottom: $128^2$ resolution comparing Sharp-VNK, VNK, and GlobalMP.}
    \label{fig:leveque-new-midface-joined-slices}
\end{figure}

This is indicated further in \cref{fig:leveque-new-midface-joined-slices} where horizontal slices through final time profiles are concatenated and plotted for several limiters. In the bottom row for both the cone, and the negative cosine, no visible differences in peak retention ability are observable between the VNK, sharp-VNK, and global maximum principle limiter. The global maximum principle limiter also exhibits undershoots and overshoots for the slotted cylinder, since the limiter never activates. In the top row of \cref{fig:leveque-new-midface-joined-slices} we see that the sharp-$N^2(K)\cup N(K)$ and $N^2(K)\cup N(K)$ limiter are marginally less diffusive than the GBJ-$N(K)\cup K$ limiter. 

\section{Conclusion}\label{sec: Application2: conclusion: chapter 2}
We have developed a local maximum principle framework for arbitrary order finite volume schemes on general unstructured meshes. Building on a localised form of the ideas in \cite{zhang2010maximum}, the main contribution is the introduction of three new neighbourhood-parameterised limiter classes: i) the Generalised-Barth-Jespersen–$BJ(K)$ family, ii) the $GN(K)$ family and iii) the Sharp--$GN(K)$ family. We proved discrete local maximum principles for each limiter on general 2D and 3D meshes.

The GBJ–$BJ(K)$ limiter proposed here differs from existing limiters in the literature by limiting non flux contributing quadrature points locally; a necessary step for positivity \cite{zhang2010maximum} and numerically having little accuracy penalty for FV4. The GBJ–$BJ(K)$ limiter is closest to existing higher order limiters and could be adopted for ease of implementation. The Sharp--$GN(K)$ limiter and the $GN(K)$ limiter variants also differ by limiting flux contributing quadrature points in accordance with face defined local maximum principle bounds, to provably reduce the severity of correction factors and allowing more of the higher order flux to be used. 

More specifically we proved on arbitrary meshes for a general class of arbitrary order FV schemes, the following structural result: Starting from any GBJ–$BJ(K)$ limiter (including the classical BJ limiter), one can construct corresponding Sharp--$GN(K)$ and $GN(K)$ limiters that enforce the same cell-mean maximum principle while requiring less severe correction factors, clarifying how local boundedness constraints can be met with less limiting.

As a secondary result, at second order for linear subcell representations, we proved in 2D and 3D for all admissible meshes that the Sharp and non-sharp variants based on $N^2(K)\cup N(K)$ and on $VN(K)$ neighbourhoods recover the same cell-mean maximum principles as BJ and Kuzmin/MLP respectively and provably apply less severe correction factors (the GBJ-$N(K)\cup
K$ reducing to the standard BJ-limiter). Although accuracy improvements are proven analytically via correction factor inequalities on all 2D and 3D meshes, and numerical experiments confirm this, in practice the numerical advantage is marginal. The real benefit to the new limiter families, is that they naturally extend to higher-order reconstructions. 

The limiting techniques provide theoretical guarantees on local boundedness principles and are likely applicable for a wide variety of schemes and DLMP regions. However, the methodology and limiting procedure require a decomposition of the cell average onto flux contributing quadrature points, which can be difficult to find. The FV4 method has a
cell mean decomposition deduced by symmetry, with only one additional point allowing the application of the proposed limiting techniques.
The non-uniqueness (\cref{remark:zhang nonuniqueness}) of such a cell mean decomposition is likely of practical consequence to the accuracy of the limiter, and may warrant further study. Modification of the present work for steady state computation could be achieved by replacing the discontinuous correction factor function (\cref{BJ_correction_factors}) with a smooth approximation following works in \cite{venkatakrishnan1993accuracy,venkatakrishnan1995convergence,gooch2006differentiability,nishikawa2022new}. Extension of local limiting procedures for discontinuous Galerkin Finite Element methods could be developed following \cite{zhang2010maximum}.

\section*{Acknowledgements}
During this work, JW has been supported by an EPSRC studentship as part of the Centre for Doctoral Training in the Mathematics of Planet Earth (grant number EP/L016613/1). Acknowledgement of Hilary Weller and Colin Cotter for valuable insights, discussions, and advice leading to the improvement of this document. 
%Rupert Klein and Peter K Sweby, for comments leading to the improvement of this document.  

\section*{Data availability}
\paragraph{Reproducibility and code availability.}
Figures, tables and schemes in this paper can be generated by scripts contained in a public repository \url{https://github.com/jameswoodfield/Local_Boundedness_limiters}, in addition to raw outputs numerical solutions, error metrics, and diagnostics such as
mass conservation and extrema as \texttt{.npz}/\texttt{.pkl}/\texttt{.csv}/\texttt{.txt} format.  

\bibliographystyle{abbrv}
\bibliography{BIBmain}

\appendix
\crefalias{section}{appendix}
\section{First order Finite volume (FV1)}\label{sec:FV1}
This is meant as a self contained introduction that outlines the main strategy of proof and frames the developments in the Historical context. 
				
We first establish the monotonicity of a forward Euler scheme in an unstructured HHLK \cite{HHL_1976} sense, and discuss how the explicit dependence on an arbitrary velocity field fits into the notion of sign preservation, and a discrete local maximum principle. We review this historical example with unstructured notation aligning with \cite{EGH_2000} but introduce additional dependence on the velocity field, rather than separate out the averaged flow through a face as in \cite{EGH_2000}, this is to ensure the later generalisation to higher order finite volume schemes in \cref{sec: Theory: theorems: chapter 2} is straightforward.

\begin{definition}[\textbf{Forward Euler Upwind}] \label{def: Forward Euler Upwind}
The forward Euler first order upwind scheme on an unstructured mesh (denoted $\mathcal{M}$), consists of approximating the compact subcell reconstruction within each cell $K$ by the constant cell mean value $\bar{u}_{K}$. The flux through a face is approximated using the midpoint of each face, and the numerical scheme takes the following form
\begin{align}
\bar{u}^{n+1}_{K} &= \bar{u}^{n}_{K} -  \Delta t\sum_{L\in N(K)} \frac{|\sigma_{KL}|}{|K|} f_{KL}( \bar{u}_{K}^{n}, \bar{u}_{L}^{n}; \b v^{n}(\b x_{KL}) \cdot \b n_{KL}), \quad \forall K \in \mathcal{M}.
\end{align}
We sketch an element of the mesh in \cref{stencil:unstructured_chap2}.
The face belonging to the boundary of cell $K$ and $L$ is denoted $\sigma_{KL}$ and assumed a subset of a hyperplane in $\mathbb{R}^d$. $N(K):= \lbrace L \in \mathcal{M} \big | |\sigma_{KL}| >0 \rbrace$ denotes the set of face-sharing neighbours of cell $K$. The midpoint of face $\sigma_{KL}$ is denoted by the position vector $\b x_{KL}$. The cell center of cell $K$ is dented by the position vector $\b x_{K}$. The face centre of $\sigma_{KL}$ and cell centre $\b x_{K}$ are defined as
\begin{align}
\b x_{KL} := \frac{1}{|\sigma_{KL}|}\int_{\b x\in \sigma_{KL}} \b x dS, \quad \b x_{K} := \frac{1}{|K|}\int_{\b x\in K} \b x dV.
\end{align}

The positive and negative superscript denotes $(\cdot)^{+} := \max (0,\cdot), (\cdot)^{-} := \min (0,\cdot)$ the positive and negative component of an input. $|K|$ denotes the volume ($d$-dimensional Lebesgue measure) of the cell $K$ and $|\sigma_{KL}|$ denotes the volume/area ($d-1$-dimensional Lebesgue measure) of the face $\sigma_{KL}$. We denote $p_{K}(\b x)$, as the subcell representation of cell $K$. $f_{KL}$ denotes the flux from cell $K$ into the cell $L$. $\b n_{KL}$ is the outward unit normal from cell $K$ into cell $L$. $\b v(\b x)$ denotes the velocity.
For the advection equation, the Riemann problem is tractable and given by the upwind/donor cell numerical flux function
\begin{align}
f_{KL} &= f_{KL}(a_{K},b_{L}; \b v \cdot \b n_{KL}) = [\b v \cdot \b n_{KL}]^+ a_{K}  +[\b v \cdot \b n_{KL}]^- b_{L}.
\end{align}
\end{definition}\smallskip

% \begin{figure}
% \centering
% 	\hspace{10pt}	\includegraphics[scale=1.2,trim={10mm 5mm 40mm 10mm},clip]{Pictures/00_geogebra_export_3.png} 
% 	\vspace{-10pt}
% 	\caption{Diagram of cell $K$, and the face $\sigma_{KL}$ of a face sharing neighbour $L \in N(K)$, with outward unit normal $\b n_{KL}$.}
% 	\label{stencil:unstructured_chap2}
% \end{figure}

More generally the definition of a consistent conservative monotone numerical flux function as defined in \cite{EGH_2000}, can be trivially extended to schemes with a faced defined velocity field as follows.
\begin{definition}\label{def:flux conservative consistent lipshitz monotone} A consistent conservative monotone numerical flux function satisfies the following properties. The numerical flux function is  consistent with the boundary flux if $f_{KL}(a,a;\b v \cdot \b n_{KL}) = f(a)\b v \cdot \b n_{KL}$. The numerical flux should also inherit the conservative properties of the continuous flux, $f_{KL}(a,b;c_{KL}) = - f_{LK}(b,a;c_{LK})$. The map defined by the numerical flux $f_{KL}$ is a monotonic flux function in the sense that it is non-decreasing with respect to the first argument and non-increasing with respect to the second argument $\partial_a f_{KL}(a,b;\b v \cdot \b n_{KL} )  \geq 0$, $\partial_b f_{KL}(a,b;\b v \cdot \b n_{KL} )\leq 0$.
\end{definition}

\begin{example} The upwind numerical flux $f_{KL}(a_{K},b_{L}; \b v \cdot \b n_{KL} ) = [\b v \cdot \b n_{KL}]^+ a_{K}  +[\b v \cdot \b n_{KL}]^- b_{L}$, is a consistent conservative monotone numerical flux function satisfying \cref{def:flux conservative consistent lipshitz monotone} for the flux form advection equation. 
\end{example}
The upwind numerical flux $f_{KL}(a_{K},b_{L}; \b v \cdot \b n_{KL} ) = [\b v \cdot \b n_{KL}]^+ a_{K}  +[\b v \cdot \b n_{KL}]^- b_{L} $ is \emph{consistent} with respect to the physical value at the boundary since it satisfies the following condition
\begin{align}
f_{KL}(a,a;\b v \cdot \b n_{KL}) =  [\b v \cdot \b n_{KL}]^+ a +[\b v \cdot \b n_{KL}]^- a= a (\b v \cdot \b n_{KL}).
\end{align}
The numerical flux is \emph{conservative} since the flux out of one face is the flux going into another
\begin{align}
f_{KL}(a,b;\b v \cdot \b n_{KL} ) &= a (\b v \cdot \b n_{KL})^+ + b(\b v \cdot \b n_{KL})^{-} \\
&= a (\b v \cdot -\b n_{LK})^+ + b(\b v \cdot -\b n_{LK})^{-} \\
&= -a (\b v \cdot \b n_{LK})^- - b(\b v \cdot \b n_{LK})^{+} \\
&= -f_{LK}(b,a,\b v \cdot \b n_{LK} ).
\end{align}
The numerical flux is \emph{monotone} by direct computation
\begin{align}
\partial_a f_{KL}(a,b;\b v \cdot \b n_{KL} ) &=  (\b v\cdot \b n_{KL})^{+} \geq 0,\\
\partial_b f_{KL}(a,b;\b v \cdot \b n_{KL} )&=  (\b v\cdot \b n_{KL})^{-}  \leq 0.
\end{align}

\begin{theorem}[Forward Euler HHLK monotone \cite{HHL_1976}] 
	Given a consistent conservative monotone numerical flux [\cref{def:flux conservative consistent lipshitz monotone}], the forward Euler scheme [\cref{def: Forward Euler Upwind}] is a monotone function of surrounding cell mean values. This is sufficient for sign preservation for compressible flow, provided the following Courant number restriction
	 \begin{align}
C_{K} = \Delta t  \sum_{L\in N(K)} \frac{|\sigma_{KL}|}{|K| } \partial_a f_{KL} \leq 1
	 \end{align}
  holds.
If in addition, the velocity field allows a discrete divergence-free condition of the following form \begin{align}
0 = \sum_{L\in N(K)} \frac{|\sigma_{KL}|}{|K|}(\b v\cdot \b n_{KL}), \label{eq:incompressible form one}
\end{align}
then the scheme has the local maximum principle \begin{align}
	    \min_{L\in N(K)\cup K} \bar{u}_{L}^{n}\leq \bar{u}^{n+1}_{K} \leq \max_{L\in N(K)\cup K} \bar{u}_{L}^{n}
	\end{align}
 with respect to neighbouring (face sharing) cell mean values.
\end{theorem}
\begin{proof} Following \cite{HHL_1976,EGH_2000}, by differentiating the function
	\begin{align}
		\bar{u}^{n+1}_{K} &= H\left(	\bar{u}_{K}, \lbrace 	\bar{u}_{L}\rbrace_{\forall L\in N(K)} ,  \lbrace \b v \cdot n_{KL}\rbrace_{\forall L\in N(K)} \right) := 	\bar{u}_{K} - \Delta t  \sum_{L\in N(K)} \frac{|\sigma_{KL}|}{|K|}f_{KL}^{n}(\bar{u}_{K},\bar{u}_{L},\bv \cdot \b n_{KL} ),\label{eq:evolution fe}
	\end{align}
	with respect to each cell mean argument 
	\begin{align}
		\frac{\partial H}{\partial 	\bar{u}_{L}} &= -\Delta t \frac{|\sigma_{KL}|}{|K|} \partial_{\bar{u}_{L}} f_{KL} \geq 0, \quad  \forall L \in N(K), \qquad
		\frac{\partial H}{\partial 	\bar{u}_{K}} =1 - \Delta t  \sum_{L\in N(K)} \frac{|\sigma_{KL}|}{|K|} \partial_{\bar{u}_{K}  }f_{KL}   \geq 0,
	\end{align}
 the scheme is a monotone function of surrounding cell mean values, provided that the following definition of a local cell defined Courant number
	\begin{align}
		C_{K} = \Delta t  \sum_{L\in N(K)} \frac{|\sigma_{KL}|}{|K|} \partial_{	\bar{u}_{K} }f_{KL} \leq 1,
	\end{align}
holds for all cells $K\in \mathcal{M}$. 
This notion of monotonicity implies sign preservation property for arbitrary velocity fields. If one additionally assumes a discrete divergence free condition of the form
\cref{eq:incompressible form one} and uses the consistency of the numerical fluxes \cref{def:flux conservative consistent lipshitz monotone} one can establish
\begin{align}
0 = \sum_{L\in N(K)} \frac{|\sigma_{KL}|}{|K|} f_{KL}(c,c;\b v\cdot \b n_{KL}), \quad \forall c \in \mathbb{R},
\end{align}
from which it can be verified that the numerical scheme is constancy preserving in the following sense
\begin{align}
c & = H\left(\bar{u}_{K}=c, \lbrace \bar{u}_{L} =c\rbrace_{\forall L\in N(K)}, \lbrace \b v \cdot \b n_{KL}\rbrace_{\forall L\in N(K)}\right),\quad \forall c \in \mathbb{R}.
\end{align}
Temporally setting local minima $m_K$ and maxima $M_K$ to be the neighbour inclusive cell mean values as follows
	\begin{align}
		m_{K} = \min_{L\in N(K)\cup K} \bar{u}_{L}, \quad M_{K} = \max_{L\in N(K)\cup K} \bar{u}_{L},
	\end{align}
	 the inclusive face sharing local maximum principle 
	\begin{align}
		m_{K} =
		H\left(m_{K}, \lbrace m_{K} \rbrace_{\forall L\in N(K)},  \lbrace \b v \cdot \b n_{KL}\rbrace_{\forall L\in N(K)}\right) \leq \bar{u}_{K}^{n+1}  \leq H\left(M_{K}, \lbrace M_{K} \rbrace_{\forall L\in N(K)}, \lbrace \b v \cdot \b n_{KL}\rbrace_{\forall L\in N(K)}\right) = M_{K},
	\end{align}
 can be established, using consistency and the monotonicity of the function $H$. 
\end{proof}

%Harten Hyman Lax and Keyfitz go further in \cite{HHL_1976} and establish a nonlinear order barrier theorem strictly stronger than the Godunov order barrier, and show that the Monotone definition forms a L1 contractive semigroup see \cite{CT_1980}. 

The HHLK-monotonicity implies sign preservation, and a divergence-free velocity field is required for a discrete local maximum principle. This background material motivates what it means in this paper for a higher-order scheme to retain monotone properties, sign preservation for compressible flow, and a discrete local maximum principle for incompressible flow. 

% In the next section, we introduce some higher-order finite volume methods and develop sufficient conditions on multidimensional slope limiters for the preservation of a local maximum principle of the following form
% \begin{align}
% 	\bar{u}^n_{K}\in [m_K,M_K], \quad \forall K.
% \end{align}
% Where $m_K,M_K$ may depend on some local quantities such as a stencil of local cell means.

\section{BJ-correction factors derivation}\label{proof:bj_derivation}
We must assume $m \leq \bar{u}_K \leq M$ otherwise no $\alpha\in[0,1]$ exists and bounds are not sensible, boundedness of $\tilde{p}_K$ at $\b x_q$ can be written as
\begin{align}
m\leq \tilde{p}_K(\b x_q) \leq M \iff m - \bar{u}_K\leq \alpha_q(p(\b x_q)-\bar{u}_{K})\leq M - \bar{u}_K\label{eq:max}
\end{align}
If $p_K(\b x_{q}) - \bar{u}_{K} > 0$ condition \cref{eq:max} is equal to 
\begin{align}
\frac{m - \bar{u}_{K}}{p_K(\b x_{q}) - \bar{u}_{K}} \leq\alpha_q \leq \frac{M - \bar{u}_{K}}{p_K(\b x_{q}) - \bar{u}_{K}},\label{eq:inequalities}
\end{align}
the LHS inequality in \cref{eq:inequalities} is automatically satisfied since $\alpha_q$ is non-negative. 
If $ p_K(\b x_{q}) - \bar{u}_{K} < 0$, then condition \cref{eq:max} is equal to 
\begin{align}
\frac{m - \bar{u}_{K}}{p_K(\b x_{q}) - \bar{u}_{K}} \geq \alpha_q \geq \frac{M - \bar{u}_{K}}{p_K(\b x_{q}) - \bar{u}_{K}}, \label{eq:inequalities2}
\end{align}
and the RHS inequality in \cref{eq:inequalities2} is automatically satisfied for any $\alpha_q\geq 0$, since the RHS is negative.
Finally, if $p_K(\b x_{q}) - \bar{u}_{K} = 0$ condition \cref{eq:max} is automatically satisfied. 
We then define the limiter function from the remaining conditions and $\alpha\in [0,1]$, giving
\begin{align}
\alpha_q := \min(1,a) \quad \text{where} \quad
a := \begin{cases}
\frac{M - \bar{u}_{K}}{p_K(\b x_{q})-\bar{u}_K}, \quad \text{if} \quad p_K(\b x_{q})-\bar{u}_K>0,\\
\frac{m - \bar{u}_{K}}{p_K(\b x_{q})-\bar{u}_K},\quad \text{if} \quad p_K(\b x_{q})-\bar{u}_K<0,\\
1, \quad \text{if} \quad p_K(\b x_{q})-\bar{u}_K=0.\\
\end{cases}
\end{align}

\section{Second order finite volume schemes FV2}\label{sec:fv2}
In this section, we:
\begin{itemize}
    \item Define a second order finite volume scheme FV2.
    \item Derive for all 2D and 3D meshes that there exists a Zhang acceptable decomposition of the cell average onto entirely flux contributing quadrature points. 
    \item Define the (Sharp)$N^2(K)\cup N(K)$ limiter, the (sharp)$VN(K)$ limiter under the FV2 scheme simplification, as well as the Barth-Jespersen limiter, and the Kuzmin/(MLP) limiter.
    \item Prove (Sharp)$VN(K)$ limiter, and Kuzmin's limiter satisfy a cell mean local maximum principle on $VN(K)$. 
    \item Prove that the Sharp-$VN(K)$ limiter does less limiting than Kuzmin limiter on all meshes.
\end{itemize}

\subsection{FV2 method}
The FV2 scheme uses the midpoint rule to approximating the flux through a face (i.e. there is one quadrature point per face chosen as the face midpoint $\b x_q = \b x_{KL}$ with quadrature weight 1.) such that the scheme \cref{eq:dg} simplifies to
\begin{align}
	\bar{u}^{n+1}_{K}   &= \bar{ u}^{n}_{K}   - \frac{ \Delta t }{ |K|} \sum_{L\in N(K)} |\sigma_{KL}|  f_{KL}( \tilde{p}_{K}(\b x_{KL}), \tilde{p}_{L}(\b x_{KL}) ; \b v(\b x_{KL}) \cdot \b n_{KL}), \label{eq:dgFV2}
\end{align}
where a linear subcell representation is chosen within each cell
\begin{align}
\tilde{p}_K(\b x) = \bar{u}_K + \alpha_K\nabla u_K \cdot (\b x-\b x_K),
\end{align}
where $\alpha_K\in [0,1]$ is attained from a slope limiting procedure.

\subsection{FV2: Zhang-acceptable decompositions of the cell mean}

\begin{theorem}[Zhang-acceptable decomposition FV2]\label{thm:ZA-FV2}
For a linear reconstruction $p_{K}(\b x)$ the cell mean $\bar{u}_{K}$, can be written entirely as a convex combination of flux contributing facet midpoint evaluations $p_{K}(\b x_{KL})$ as follows
\begin{align}
\bar{u}_{K} = p_{K}(\b x_c)  = \sum_{L\in N(K)} \gamma^{K}_L p_K(\b x_{\sigma_{KL}}), \quad \gamma^{K}_{L} >0,\quad\sum_{L\in N(K)}w^{K}_{L} = 1. 
\end{align}
Should $K$ be a convex polyhedra in 3D the weight is $\gamma^{K}_{L} = |\mathcal{P}_{KL}|/|K|$ a volume ratio defined by $\mathcal{P}_{KL}$ a pyramid constructed from joining vertices of $\sigma_{KL}$ to $\b x_c$. Should $K$ be a convex polygon in 2D, with $n_f$ faces/vertices the weight is simply given by $\gamma^{K}_{L}= \frac{1}{n_f}$. 
\end{theorem}

\begin{proof}
In 2D; let the polygon $K$ have $n_f$ vertices and $n_f$ faces, the average of an affine function is the value at the centroid. The value of an affine function at the centroid is the average of vertex values. The value at the face midpoint is the average of vertex values, and since each face has two vertices we deduce the 2D-Zhang-acceptable decomposition:
\begin{align}
\bar{u}_{K} = p_{K}(\b x_c) = \frac{1}{n_f}\sum_{v\in V(K)} p_{K}(\b x_{v}) = \frac{1}{n_f}\sum_{L\in N(K)} p_K(\b x_{\sigma_{KL}}), \quad \gamma^{K}_{L}=\frac{1}{n_f}>0,\quad\sum_{L\in N(K)}\gamma^{K}_{L} = 1. 
\end{align}

In 3D; we require a different construction.
The average of an affine function is its value at the centroid. There exists a convex combination decomposition of the cell centroid in-terms of facet centroids $\b x_{KL}$. Take face $\sigma_{KL}$ and cell centroid $\b x_{c} \in \text{Int}(K)$, and define a pyramid, $\mathcal{P}_{KL}$ by joining vertices of $\sigma_{KL}$ to $\b x_{c}$. Pyramids have positive area,  $|\mathcal{P}_{KL}|>0$ via convexity of $K$ and sum to give the volume of the cell  $|K|=\sum_{L\in N(K)}|\mathcal{P}_{KL}|$. The cell centroid can be written in terms of a convex combination of pyramid centroid values $\b x_{\mathcal{P}_{KL}}$ as follows
\begin{align}
\b x_{c} = \frac{1}{|K|}\int_{K}\b xdV = \frac{1}{|K|}\sum_{L\in N(K)} \int_{\mathcal{P}_{KL}}\b x dV = \frac{1}{|K|}\sum_{L\in N(K)} |\mathcal{P}_{KL}|\b x_{\mathcal{P}_{KL}}, 
\quad \b x_{\mathcal{P}_{KL}} :=
\frac{1}{|\mathcal{P}_{KL}|}\int_{\mathcal{P}_{KL}}\b x dV\end{align}
The centroid of each pyramid in 3D, can be expressed as cell centroid and facet centroid values via $\b x_{\mathcal{P}_{KL}} = \frac{1}{4}\b x_c +\frac{3}{4} \b x_{KL}$, substituting and rearranging gives
% \begin{align}
% \b x_{c} = \frac{1}{|K|}\sum_{L\in N(K)} |\mathcal{P}_{KL}|(\frac{1}{4}\b x_c +\frac{3}{4} \b x_{KL})= 1/4\b x_{c} + \frac{1}{|K|}\sum_{L\in N(K)}|\mathcal{P}_{KL}|3/4 \b x_{KL}
% \end{align}
% Giving
the following convex combination
\begin{align}
    \b x_{c} = \sum_{L\in N(K)} \frac{|\mathcal{P}_{KL}|}{|K|}\b x_{KL}.
\end{align}
Affine maps preserve convex combinations so that
\begin{align}
    \bar{u}_K = p_K(\b x_c) = \sum_{L\in N(K)} \gamma^{K}_{L}p_{K}(\b x_{KL}), \quad \gamma^{K}_{L} = \frac{|\mathcal{P}_{KL}|}{|K|}>0, \quad \sum_{L\in N(K)} \gamma^{K}_{L} = 1. \label{eq:FV2_zhang}
\end{align}
Here $|\mathcal{P}_{KL}|$ denotes the area of a pyramid attained from joining vertices of facet $\sigma_{KL}$ to cell center $\b x_{c}$. For a useable CFL condition, $|\mathcal{P}_{KL}|$ can be computed explicitly by the formula $|\mathcal{P}_{KL}|=1/3|\sigma_{KL}|d_{KL}$ where the perpendicular height from face to centre is given by $d_{KL} := \b n_{KL}\cdot(\b x_{\sigma_{KL}}-\b x_{c})$ (positive since $K$ is convex). 
\end{proof}

\Cref{thm:ZA-FV2} implies that in 2D we have the following monotone convex combination representation of FV2
\begin{align}
	\bar{u}^{n+1}_{K}   &= \sum_{L\in N(K)} \frac{|n_f|}{|K|} \left( \tilde{p}_K(\b x_{KL})   -  \Delta t  \frac{|\sigma_{KL}|}{|n_f|}  f_{KL}( \tilde{p}_{K}(\b x_{KL}), \tilde{p}_{L}(\b x_{KL}) ; \b v(\b x_{KL}) \cdot \b n_{KL}) \right), \label{eq:zhang_FV2_2d}
\end{align}
 and the following monotone convex combination representation of FV2 in 3D
\begin{align}
	\bar{u}^{n+1}_{K}   &= \sum_{L\in N(K)} \frac{|\mathcal{P}_{KL}|}{|K|} \left( \tilde{p}_K(\b x_{KL})   -  \Delta t  \frac{|\sigma_{KL}|}{|\mathcal{P}_{KL}|}  f_{KL}( \tilde{p}_{K}(\b x_{KL}), \tilde{p}_{L}(\b x_{KL}) ; \b v(\b x_{KL}) \cdot \b n_{KL}) \right). \label{eq:zhang_FV2}
\end{align}
We directly can infer montone-CFL restrictions of FV2 for any mesh in 2D and 3D. 
We also see why it is sufficient to limit $\tilde{p}_{K}(\b x)$, and $\tilde{p}_{L}(\b x)$, at flux contributing quadrature points $\b x_{KL}$, $\forall L\in N(K)$, since no other points arise in the above decompositions of the cell average. 

\subsection{Limiters for FV2}
In the context of second order finite volume schemes (FV2) we define: the Sharp--$N^2(K)\cup N(K)$ limiter, the $N^2(K)\cup N(K)$ limiter, the $VN(K)$ limiter, (which simplify under FV2) we also define  the Barth-Jespersen limiter, and the Kuzmin(MLP) limiter.

\begin{method}[Sharp--$N^2(K)\cup N(K)$-limiter for FV2]\label{method:Sharp-FV2N2KNK}\noindent
\begin{enumerate}
	\item Per face $\sigma_{KL}$, compute  the local face defined maximum principle bounds
	\begin{align}
	\left[ m_{\sigma_{KL}},M_{\sigma_{KL}} \right]= \left[\max\left(
    \min_{M\in N^2(K)\cup N(K)}\bar{u}_{M}^n,\min_{M\in N^2(L)\cup N(L)}\bar{u}_{M}^n\right) ,\min\left(\max_{M\in N^2(K)\cup N(K)}\bar{u}_{M}^n,\max_{M\in N^2(L)\cup N(L)}\bar{u}_{M}^n\right) \right].
		\end{align}
		\item 
		Per cell $K$ compute all the quadrature point corrections factors
        \begin{align}
            \alpha_{KL} = \alpha_{q}(m_{\sigma_{KL}},M_{\sigma_{KL}},\bar{u}_{K},p_K(\b x_{KL})), \quad \forall L \in N(K),
        \end{align} 
        via \cref{BJ_correction_factors}. %ensuring $\bar{u}_{K}+\alpha_L (p_K(\b x_{KL})-\bar{u}_{K})\in [m_{\sigma_{KL}},M_{\sigma_{KL}}]$.
        \item 
        Choose the smallest value, 
		\begin{align}
			\alpha^{N^{2}(K)\cup N(K)}_K = \min_{\forall L \in N(K)} \alpha_{KL}.
		\end{align}
	\end{enumerate}
\end{method}

\begin{method}[$N^2(K)\cup N(K)$-limiter for FV2]\label{method:FV2N2KNK}\noindent
	\begin{enumerate}
		\item Per face $\sigma_{KL}$, compute  the local face defined maximum principle bounds
		\begin{align}
			\left[ m_{\sigma_{KL}},M_{\sigma_{KL}} \right]= \left[\min_{M \in N(L) \cup N(K)} \bar{u}^n_{M},\max_{M \in N(L) \cup N(K)} \bar{u}^n_{M}\right].
		\end{align}
		\item 
		Per cell $K$ compute all the quadrature point corrections factors
        \begin{align}
            \alpha_{KL} = \alpha_{q}(m_{\sigma_{KL}},M_{\sigma_{KL}},\bar{u}_{K},p_K(\b x_{KL})), \quad \forall L \in N(K),
        \end{align} 
        via \cref{BJ_correction_factors}. %ensuring $\bar{u}_{K}+\alpha_L (p_K(\b x_{KL})-\bar{u}_{K})\in [m_{\sigma_{KL}},M_{\sigma_{KL}}]$.
        \item 
        Choose the smallest value, 
		\begin{align}
			\alpha^{N^{2}(K)\cup N(K)}_K = \min_{\forall L \in N(K)} \alpha_{KL}.
		\end{align}
	\end{enumerate}
\end{method}

\begin{method}[Barth-Jespersen limiter for FV2]\label{sec:Barth-Jespersen}
The Barth-Jespersen limiter is defined as follows,	\noindent 
	\begin{enumerate}
		\item Per cell $K$ compute the local bounds
		\begin{align}
			[m_{K},M_{K}] := \left[\min_{L\in N(K)\cup K} \bar{u}_{L}^n,\max_{L\in N(K)\cup K} \bar{u}_{L}^n\right].
		\end{align}
		\item Per cell $K$ compute all the quadrature point corrections factors
        \begin{align}
            \alpha_{KL} = \alpha_{q}(m_{K},M_{K},\bar{u}_{K},p_K(\b x_{KL})), \quad \forall L \in N(K),
        \end{align} 
        via \cref{BJ_correction_factors}.
		\item Choose
		\begin{align}
			\alpha^{BJ}_{K} = \min_{L\in N(K)} \alpha_{KL}.
		\end{align} 
	\end{enumerate}
\end{method}

\begin{method}[Sharp-$VN(K)-Limiter$ for FV2]\noindent
\begin{enumerate}
\item Per face $\sigma_{KL}$ compute the bounds
\begin{align}
[m_{\sigma_{KL}},M_{\sigma_{KL}}] = \left[\max\left(\min_{M\in VN(K)}\bar{u}_{M}, \min_{M\in VN(L)} \bar{u}_{M}\right),\min\left(\max_{M\in VN(K)}\bar{u}_{M}, \max_{M\in VN(L)} \bar{u}_{M}\right) \right].
\end{align}
\item Per cell $K$ compute all the quadrature point corrections factors
        \begin{align}
            \alpha_{KL} = \alpha_{q}(m_{\sigma_{KL}},M_{\sigma_{KL}},\bar{u}_{K},p_K(\b x_{KL})), \quad \forall L \in N(K),
        \end{align} 
        via \cref{BJ_correction_factors}.
		\item Choose
		\begin{align}
			\alpha^{\sharp,VNK}_{K} = \min_{L\in N(K)} \alpha_{KL}.
		\end{align} 
\end{enumerate}
\end{method}

\begin{method}[$VN(K)-Limiter$ for FV2]\noindent
\begin{enumerate}
\item Per face $\sigma_{KL}$ compute the bounds
\begin{align}
[m_{\sigma_{KL}},M_{\sigma_{KL}}] = \left[\min_{M\in VN(\sigma_{KL})}\bar{u}_{M}, \max_{VN(\sigma_{KL})} \bar{u}_{M} \right].
\end{align}
\item Per cell $K$ compute all the quadrature point corrections factors
        \begin{align}
            \alpha_{KL} = \alpha_{q}(m_{\sigma_{KL}},M_{\sigma_{KL}},\bar{u}_{K},p_K(\b x_{KL})), \quad \forall L \in N(K),
        \end{align} 
        via \cref{BJ_correction_factors}.
		\item Choose
		\begin{align}
			\alpha^{VNK}_{K} = \min_{L\in N(K)} \alpha_{KL}.
		\end{align} 
\end{enumerate}
\end{method}

\begin{method}[Kuzmin limiter for FV2] The Park et al.(MLP)/Kuzmin vertex-based limiter \cite{park2010multi, kuzmin2010vertex} is defined as follows:\label{sec:Kuzmin}
\noindent 
\begin{enumerate}
	\item Compute the vertex defined local bounds 
	\begin{align}
		[m_{v},M_{v}] := \left[\min_{i\in N(v)} \bar{u}_{i}^n,\max_{i\in N(v)} \bar{u}_{i}^n\right],
	\end{align}
	where $N(v)$ denotes the set of cells which share the vertex $v$.
	\item Per vertex compute a correction factor $\alpha_v$, 
	\begin{align}
		\alpha_{v} = \begin{cases}
			\min \left\lbrace 1, \frac{M_v - \bar{u}_{K}}{p_K(\b x_v) - \bar{u}_{K}}  \right\rbrace \quad \text{if} \quad  p_K(\b x_{v}) - \bar{u}_{K}>0,\\
			\min \left\lbrace 1,\frac{m_v - \bar{u}_{K}}{p_K(\b x_v) -\bar{u}_{K}} \right\rbrace \quad \text{if} \quad  p_{K}(\b x_{v}) - \bar{u}_{K}<0,\\
			1 \quad  \text{if} \quad  p_K(\b x_{v}) - \bar{u}_{K} =0.
		\end{cases}
	\end{align}
	to ensure that all vertex points are locally bounded by their local vertex sharing neighbours $[m_{v},M_{v}]$. 
	\item Then the entire subcell representation is limited based on the worst violator of the local vertex based maximum principles,
	\begin{align}
		\alpha^{Kuz}_{K} := \min_{v \in K }  \alpha_v
	\end{align} 
    to give $\tilde{p}_K(x) = \bar{u}_k + \alpha^{Kuz}_K(p(x)-\bar{u}_k)$ bounded at corner points $\b x_v$ by $[m_v,M_v]$.
\end{enumerate}
\end{method}

% \begin{remark}[Computational cost]
% At second order the only difference between the BJ limiter and the $N^2(K)\cup N(K)$ limiter, is in the face defined bounds. One can take the BJ bounds $[m^{BJ}_K,M^{BJ}_K]$, and extend them to the larger face defined bounds, by taking $[m_{\sigma_{KL}},M_{\sigma_{KL}}] = [\min\lbrace m^{BJ}_K,m^{BJ}_L\rbrace,\max\lbrace M^{BJ}_K,M^{BJ}_L\rbrace]$ for each face, this is not costly. For the Sharp variant, one takes the larger $[m_K,M_K]$ bounds defined on the larger $N^2(K)\cup N(K)$ region, and taking  $[m_{KL},M_{KL}] = [\max\lbrace m_K,m_L\rbrace,\min\lbrace M_K,M_L\rbrace]$, also not costly. Both these implementations are unstructured. 
% \end{remark}

\subsection{FV2 Maximum principles}

We prove different discrete local maximum principles for the FV2 scheme with the Kuzmin/MLP limiter. All other limiters attain maximum principles as corollaries of \cref{thm: slope limiters}.

\begin{theorem}FV2 (\cref{eq:dgFV2}) with the Kuzmin-limiter has a resulting cell mean value $\bar{u}_{K}^{n+1}$
satisfying the discrete local maximum principle
	\begin{align}
		\min_{L \in VN(K)} \bar{u}^n_{L}\leq \bar{u}_{K}^{n+1} \leq \max_{L \in VN(K)}\bar{u}^n_{L}.\label{eq:cell mean vertex max principle}
	\end{align}
\end{theorem}
\begin{proof}
In 2D on an arbitrary mesh the face midpoint of $\sigma_{KL}$ can be expressed as $\b x_{KL} = 1/2 (\b x_{v_{1}} +  \b x_{v_{2}})$ where $\b x_{v_{1}},\b x_{v_{1}}$ are the locations of the 2 vertices of that edge. Furthermore for a linear subcell representation, one can relate values at cell midpoints to values at vertices via $p_{K}(\b x_{KL}) = 1/2 p_{K}(\b x_{v_{1}})+1/2 p_{K}(\b x_{v_{2}})$.
Therefore, the Kuzmin limiter ensures that the face defined maximum principle is satisfied
\begin{align}
    p_{K}(\b x_{KL}) \in [m^{av}_{\sigma_{KL}},M^{av}_{\sigma_{KL}}]:= \left[\sum_{v_{i}\in \sigma_{KL}} \frac{1}{2}\left(\min_{j\in N(v_{i})}\bar{u}_{i}\right),\sum_{v_{i}\in \sigma_{KL}} \frac{1}{2}\left(\max_{j\in N(v_{i})}\bar{u}_{i}\right) \right],
\end{align}
and via face symmetry of such a region $p_{L}(\b x_{KL}) \in [m^{av}_{\sigma_{KL}},M^{av}_{\sigma_{KL}}]$. Therefore one can apply \cref{thm: slope limiters}, to guarantee that 
\begin{align}
\bar{u}^{n+1}_{K}\in \left[ \min_{L\in N(K)}\left( \sum_{v_{i}\in \sigma_{KL}} \frac{1}{2}\left(\min_{j\in N(v_{i})}\bar{u}^n_{i}\right)\right), \max_{L\in N(K)}\left( \sum_{v_{i}\in \sigma_{KL}} \frac{1}{2}\left(\max_{j\in N(v_{i})}\bar{u}^n_{i}\right) \right)   \right].
\end{align}
This is a stronger statement than required for the proof, proving it in 2D. 

% Curiously, this formulation of the Kuzmin limiter as a face defined limiter appears in \cite{zhang2018modified} under another name: MLP-pw limiter, where face defined bounds are taken as the average of corner defined vertex maximum principles, this appears to be identical to the Kuzmin limiter in terms of the resulting correction factor but reformulated and implemented as Barth-Jespersen type.

In 3D on an arbitrary mesh, whose faces consist of arbitrary 2D convex polygons, the face midpoint can always \footnote{One constructive proof of this fact involves triangulating the convex polygon $\sigma_{KL}$ without introducing new vertices, using that the average/integral of an affine function over a triangle equals the average of its values at vertices, then regrouping terms and swapping summands into a sum over vertices.} be expressed as a convex combination of corner defined points
\begin{align}
    \b x_{KL} = \sum_{v\in V(\sigma_{KL})} w_{v,\sigma_{KL}} \b x_{v},\quad \text{where}\quad w_{v,\sigma_{KL}} > 0, \quad \forall v \in V(\sigma_{KL}), \quad \text{and} \quad \sum_{v \in V(\sigma_{KL})} w_{v,\sigma_{KL}}=1. \label{eq:facemp_is_convex_combo_of_corner}
\end{align}

% The easiest example is a triangle since $\b x_{KL} = 1/3 \sum_{v\in V(\sigma_{KL})} \b x_{v}$. For a 4-sided face $\sigma_{KL}$ with vertices $\lbrace\b x_{i}\rbrace_{i=1}^4$, seek a triangulation into two triangles $T_1,T_2$, where the triangulation is made from a line between the $x_2,x_3$ vertex, then $\b x_{KL} = \frac{|T_1|}{3(|T_1|+|T_1|)}x_1 + \frac{1}{3}x_2+ \frac{1}{3}x_3 + \frac{|T_2|}{3(|T_1|+|T_1|)}x_4$, a convex combination. The general case follows non uniquely and gives the abstract formulation \cref{eq:facemp_is_convex_combo_of_corner}.

The linear subcell representation is affine and one may verify by a direct calculation, that the value at a midpoint of a face is a convex combination of the values at the vertices of that face
\begin{align}
    p_K(\b x_{KL}) = \sum_{v \in V(\sigma_{KL})} w_{v,\sigma_{KL}} p_{K}(\b x_{v}). \label{3D facet_vertex_relationship}
\end{align}
The Kuzmin-limiter enforces the vertex defined bounds $p_{K}(\b x_v) \in [m_v,M_v] = [\min_{j\in N(v)}\bar{u}_{j},\max_{j\in N(v)}\bar{u}_{j}]$. Therefore we attain the face defined bounds, 
\begin{align}
    p_K(\b x_{KL}) \in \left[ \sum_{v \in V(\sigma_{KL})} w_{v,\sigma_{KL}} \left(\min_{j\in N(v)}\bar{u}_{j}\right), \sum_{v \in V(\sigma_{KL})} w_{v,\sigma_{KL}} \left(\max_{j\in N(v)}\bar{u}_{j}\right)\right].
\end{align}
Where we have a face-wise convex combination of corner defined maximum principles. Since $L$ shares a face with $K$, $p_{L}(\b x_{KL})$ shares the same bounds as $p_{K}(\b x_{KL})$. As a direct consequence of \cref{thm: slope limiters}, we attain
\begin{align}
\min_{L\in N(K)}\left( \sum_{v \in V(\sigma_{KL})} w_{v,\sigma_{KL}} \left(\min_{j\in N(v)}\bar{u}_{j}\right)\right)\leq \bar{u}^{n+1}_{K}\leq \max_{L\in N(K)}\left( \sum_{v \in V(\sigma_{KL})} w_{v,\sigma_{KL}} \left(\max_{j\in N(v)}\bar{u}_{j}\right)\right).
\end{align}
This hidden maximum principle is stronger than the statement required. 
Since the face-wise convex combination is always larger than taking the face-wise minimum, and smaller that the face-wise maximum one can weaken the statement to:
\begin{align}
\min_{L\in N(K)}\left( \min_{v \in V(\sigma_{KL})} \left(\min_{j\in N(v)}\bar{u}_{j}\right)\right)\leq \bar{u}^{n+1}_{K}\leq \max_{L\in N(K)}\left( \max_{v \in V(\sigma_{KL})} \left(\max_{j\in N(v)}\bar{u}_{j}\right)\right),
\end{align}
equivalent to \cref{eq:cell mean vertex max principle} since one can simplify as follows
\begin{align}
 \max_{L\in N(K)}\left( \max_{v \in V(\sigma_{KL})} \left(\max_{j\in N(v)}\bar{u}_{j}\right)\right) =   \max_{v\in N(K)} \left(\max_{j\in N(v)}\bar{u}_{j}\right) = \max_{VN(K)}\bar{u}_{j}. 
\end{align}
\end{proof}

\subsection{Comparing 2nd order limiters, factors affecting accuracy}\label{sec:accuracy of BJ KUZ N2UN}

This section proves the $VN(K)$ limiter has less severe correction factors than the Kuzmin/MLP vertex limiter, on all meshes.

\begin{theorem}Given a arbitrary 2D mesh $\mathcal{M}$, or a arbitrary 3D mesh whose faces are flat convex polygons then the Kuzmin limiter has more severe correction factors than the $VN(K)$ limiter and the Sharp--$VN(K)$ limiter: 
\begin{align}
    \alpha_K^{Kuz} \leq \alpha^{VNK}_K\leq \alpha^{\sharp}_K.
\end{align}
\end{theorem}

\begin{proof}
In 2D where $\sigma_{KL}$ has 2 vertices, the Kuzmin limiter gives face defined bounds
\begin{align}
    p_{K}(\b x_{KL}),p_{L}(\b x_{KL}) \in [m^{av}_{\sigma_{KL}},M^{av}_{\sigma_{KL}}]:= \left[\sum_{v_{i}\in \sigma_{KL}} \frac{1}{2}\left(\min_{j\in N(v_{i})}\bar{u}_{i}\right),\sum_{v_{i}\in \sigma_{KL}} \frac{1}{2}\left(\max_{j\in N(v_{i})}\bar{u}_{i}\right) \right],
\end{align}
whereas the $VN(K)$- limiter gives the face defined bounds
\begin{align}
    p_{K}(\b x_{KL}),p_{L}(\b x_{KL}) \in [m^{*}_{\sigma_{KL}},M^{*}_{\sigma_{KL}}]:= \left[\min_{v_{i}\in \sigma_{KL}} \left(\min_{j\in N(v_{i})}\bar{u}_{i}\right),\min_{v_{i}\in \sigma_{KL}} \left(\max_{j\in N(v_{i})}\bar{u}_{i}\right) \right],
\end{align}
and the Sharp--$VN(K)$ limiter gives the face defined bounds
\begin{align}
    p_{K}(\b x_{KL}),p_{L}(\b x_{KL}) \in [m^{\sharp}_{\sigma_{KL}},M^{\sharp}_{\sigma_{KL}}]:= \left[\max\left( \min_{M\in VN(K)} \bar{u}_{M},\min_{M\in VN(L)} \bar{u}_{M} \right), \min\left( \max_{M\in VN(K)} \bar{u}_{M},\max_{M\in VN(L)} \bar{u}_{M} \right)\right]. 
\end{align}
Since
\begin{align}
[m^{av}_{\sigma_{KL}},M^{av}_{\sigma_{KL}}]
\subseteq [m^{*}_{\sigma_{KL}},M^{*}_{\sigma_{KL}}]
\subseteq [m^{\sharp}_{\sigma_{KL}},M^{\sharp}_{\sigma_{KL}}],
\end{align}
as a direct consequence of \cref{lemma:correction factor inequalities} the correction factors associated to the Sharp--$VN(K)$, the $VN(K)$ and the Kuzmin limiter at quadrature point $\b x_{KL}$ denoted $\alpha_{KL}^*$, $\alpha^{Kuz}_{KL}$, $\alpha^{\sharp}_{KL}$ respectively satisfy 
\begin{align}
   \alpha_{KL}^{Kuz}  \leq \alpha_{KL}^*\leq \alpha_{KL}^{\sharp}, \quad \forall L \in N(K). 
\end{align}
Therefore 
\begin{align}
   \alpha_K^{Kuz} =\min_{L\in N(K)}\alpha_{KL}^{Kuz}  \leq \min_{L\in N(K)} \alpha_{KL}^* = \alpha^{VNK}_{K} \leq \min_{L\in N(K)} \alpha_{KL}^{\sharp} = \alpha^{\sharp}_{K}.
\end{align}
In 3D, an entirely analogous proof holds via \cref{3D facet_vertex_relationship}.
\end{proof}

\section{Unlimited schemes}

\begin{figure}[H]
\centering
\begin{subfigure}[t]{0.45\textwidth}
  \centering
\includegraphics[width=\linewidth]{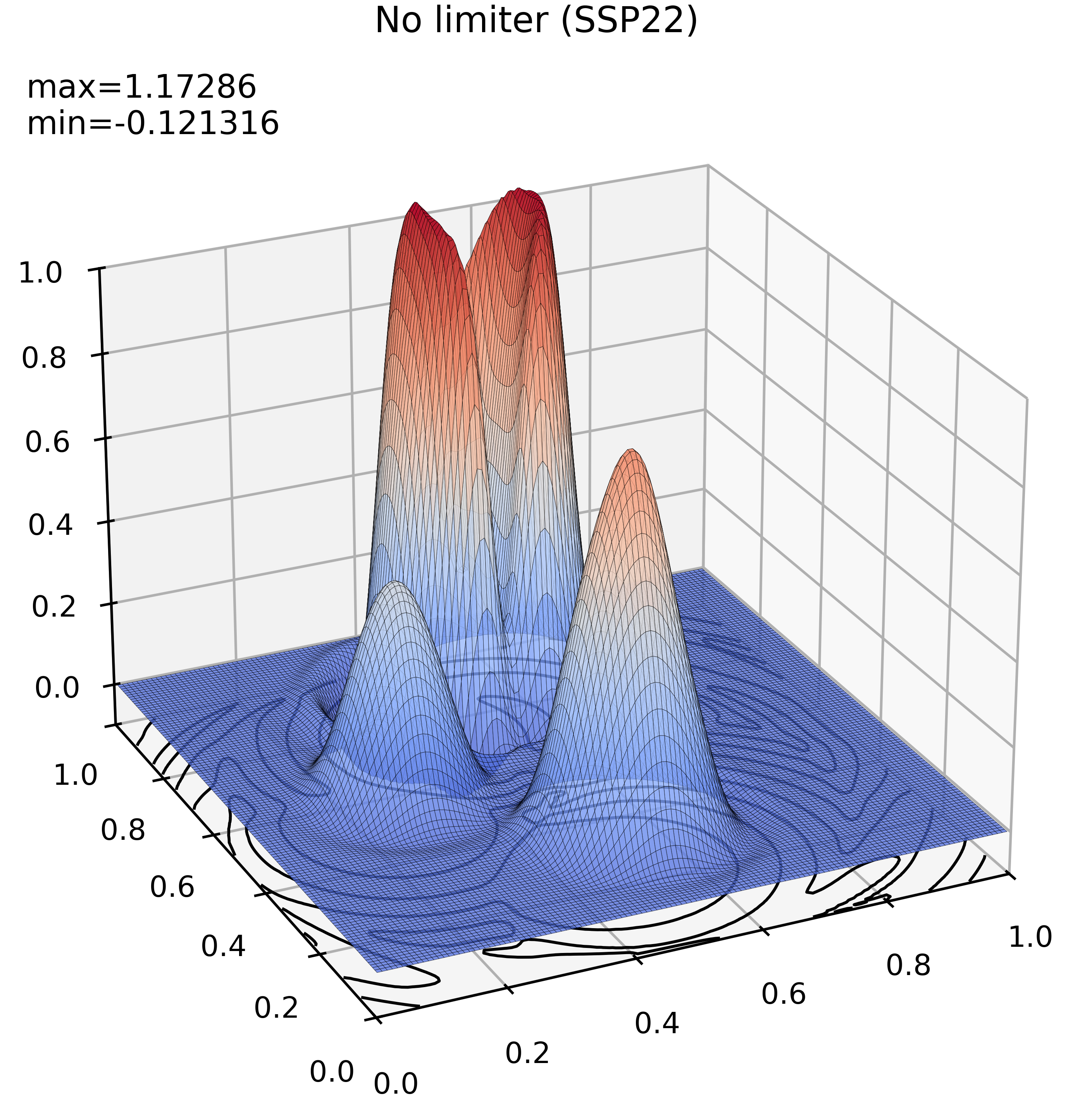}
\end{subfigure}
\begin{subfigure}
{0.45\textwidth}
  \centering
\includegraphics[width=\linewidth]{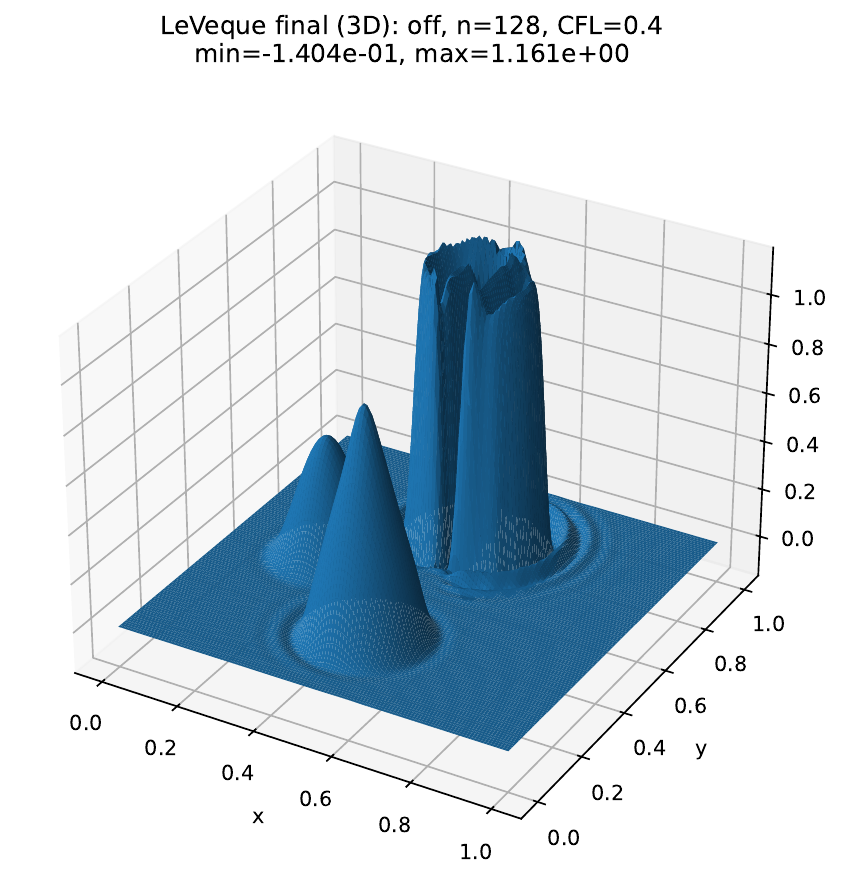}
\end{subfigure}
\caption{LeVeque solid-body rotation at $n_x=128$: FV2-SSP22 unlimited (no limiter) and FV4-SSP(ARK)44 unlimited solution at final time.}
\label{fig:leveque-3d-none-128}
\end{figure}

\end{document}